\documentclass[opre,nonblindrev]{informs4aa}
\DoubleSpacedXI % Made default 4/4/2014 at request
\usepackage{endnotes}
\let\footnote=\endnote

\usepackage{natbib}
\bibpunct[, ]{(}{)}{,}{a}{,}{,}%
\def\bibfont{\small}%
\usepackage{amsmath}
\usepackage{amssymb}
\usepackage{amsfonts}
\usepackage{amsmath}
\usepackage{bbm}
\usepackage{datetime}
\usdate
\settimeformat{ampmtime}
\usepackage[inline]{enumitem}
\usepackage{float}
\usepackage[margin=1in]{geometry}
\usepackage{graphicx}
\usepackage{latexsym}
\usepackage{mathabx}
\usepackage{mathtools}
\usepackage{multirow}
\usepackage{ifpdf}
\ifpdf
\usepackage[pdftex]{hyperref}
\else
\usepackage[hypertex]{hyperref}
\fi

\usepackage{xcolor}

\usepackage{tikz}
\usepackage{pgfplots}
\pgfplotsset{compat=1.17}

\usepackage{makecell}
\setcellgapes{3pt}

\usepackage{cancel}
\usepackage{amssymb}
\usepackage{accents}
\usepackage{verbatim}

\newcommand{\wbar}{\widebar{w}}

\newcommand{\vhat}{\widehat{v}}
\newcommand{\hhat}{\widehat{h}}
\newcommand{\xhat}{\widehat{x}}

\newcommand{\thetalb}{\underaccent{\bar}{\theta}}
\newcommand{\thetaub}{\bar{\theta}}

\newcommand{\R}{\mathbb{R}}

\newcommand{\E}[1]{{\mathbb E}\left[#1\right]}
\newcommand{\Var}[1]{{\mathrm{Var}}\left[#1\right]}

\newcommand{\1}[1]{{\mathbbm 1}\left\{ #1 \right\}}

\newcommand{\midbar}{\;\middle\vert\;}

\newcommand\cD{\mathcal D}

\newcommand\cF{\mathcal F}

\newcommand{\ceil}[1]{\left\lceil\, #1 \,\right\rceil}

\newcommand{\abs}[1]{\left\lvert\, #1 \,\right\rvert}
\newcommand{\norm}[1]{\left\lVert\, #1 \,\right\rVert}

\let\oldmarginpar\marginpar
\renewcommand\marginpar[1]{\-\oldmarginpar[\raggedleft\scriptsize \textcolor{red}{#1}]%
	{\raggedright\scriptsize \textcolor{red}{#1}}}

\renewcommand{\MSCCLASSname}{%
    \textsf{\textbf{MSC2020 subject classification:}}\enskip}

\renewcommand{\ORMSCLASSname}{%
    \textsf{\textbf{OR/MS subject classification:}}\enskip} 

\TheoremsNumberedThrough     % Preferred (Theorem 1, Lemma 1, Theorem 2)
\ECRepeatTheorems

\EquationsNumberedThrough    % Default: (1), (2), ...
\begin{document}
%%%%%%%%%%%%%%%%

%\setcounter{page}{1}
%\DOI{10.1287/xxxx.0000.0000}%

\RUNAUTHOR{A.\ Arlotto, Y.T. Kuo, and X.\ Xie}

\RUNTITLE{Reoptimization in the DSKP: from near-optimality to distributional equivalence}

% \TITLE{Bundling Information Goods of Decreasing Value}
\TITLE{Reoptimization in the dynamic and stochastic knapsack problem: from near-optimality to distributional equivalence}

\ARTICLEAUTHORS{%
	\AUTHOR{Alessandro Arlotto}
	\AFF{The Fuqua School of Business, Duke University, 100 Fuqua Drive, Durham, NC, 27708, \\ 
		\EMAIL{aa249@duke.edu}} %, \URL{}}
	\AUTHOR{Yun-Tung Kuo}
	\AFF{The Fuqua School of Business, Duke University, 100 Fuqua Drive, Durham, NC, 27708, \\ 
		\EMAIL{yun-tung.kuo@duke.edu}} %, \URL{}}

	\AUTHOR{Xinchang Xie}
	\AFF{United Parcel Service, 2311 York Road, Timonium, MD, 21093, \\  			   
		\EMAIL{xxie@ups.com}}
} % end of the block

\ABSTRACT{%
We study a dynamic and stochastic knapsack problem in which a decision maker is sequentially
presented with $T$ items with unit rewards and independent weights drawn from a 
nonnegative distribution satisfying mild regularity conditions near zero.
The decision maker seeks to maximize the expected number of selected items while satisfying a
capacity constraint and making terminal decisions as soon as each item's weight is revealed.
We study a reoptimized heuristic and compare its total number of selections with that of the
optimal dynamic programming policy.
We show that the two totals converge to the same normal limiting distribution,  after standard centering and scaling.
By contrast, another asymptotically optimal heuristic considered in the literature has strictly
larger asymptotic variance and a different, non-Gaussian limiting distribution.
These results give a finer, distributional basis for comparing policies that are otherwise
first-order indistinguishable asymptotically.
}%


\MSCCLASS{%
    Primary: 90C40, 60F05;
    secondary: 60C05, 60G42, 68W27, 90C39.
}

\ORMSCLASS{%
    Dynamic programming/optimal control: Markov;
    probability: stochastic model applications;
    analysis of algorithms: suboptimal algorithms.
}

\HISTORY{This version: November 27, 2026.
%File name: \texttt{\jobname.tex}; This version: \today, \currenttime
}
\maketitle

\section{Introduction}\label{sec:introduction}

The dynamic and stochastic knapsack problem \citep[see, e.g.,][]{PapastavrouRajagopalanKleywegt:MS1996,KleywegtPapastavrou:OR1998}
models a broad class of sequential resource-allocation problems in which a decision maker has a fixed capacity, faces a sequence of arriving items with random weights and rewards, and must decide irrevocably and in real time which items to accept, so as to maximize total expected rewards without exceeding the knapsack capacity. 
A classical example is real-time uniprocessor scheduling under conditions of overload \citep[see, e.g.,][]{BaruahHaritsaSharma:IEEE1994}, where jobs arrive sequentially and require random processing times, a single machine is available until a fixed deadline, and the objective is to maximize the number of completed jobs. 
In this scheduling context, the remaining processing time plays the role of capacity, and job durations correspond to item weights.

Dynamic programming provides a natural framework for such sequential decision problems \citep{Bellman:PUP1957}, but optimal policies are rarely available in closed form and, for the richer members of this class (multiple resource types, combinatorial state spaces, \ldots), they also become computationally intractable because of the curse of dimensionality \citep{Powell:JWS2011}.
These challenges have motivated a large literature on tractable online policies with provable performance guarantees. 
In revenue management, assortment optimization, and stochastic resource allocation, researchers have developed heuristics whose expected reward is close to that of the optimal dynamic programming policy as the planning horizon, demand, or capacity scales \citep[see, e.g.,][]{GallegoVanRyzin:MS1994,JasinKumar:MOR2012,JasinKumar:OR2013,Golrezaei:MS2014,ReimanWang:MOR2008,ArlottoGurvich:SSY2019,Vera-Banerjee:MS2021,Vera-Banerjee-Gurvich:OPRE2021,Goyal:OR2025,Jiang:arXiv2023,Vera-Arlotto-Gurvich-Levin:MOOR2025,JiangMaZhang:OR2025-degeneracy-is-ok,Zhang:arXiv2026-multi-secretary,Zhang:arXiv2026-resource-allocation,ChenChenWangZhang-WP-2026}.

Despite these successes, most existing performance guarantees remain fundamentally first-order: policies are compared through their expected reward, while the stochastic fluctuations of the reward around that expectation receive much less attention. 
Once several tractable policies achieve the same asymptotic expected reward, first-order analysis alone no longer distinguishes their quality. 
From a theoretical perspective, understanding higher-order behavior provides a finer characterization of online policies. 
From a practical perspective, however, variability is equally important, since it governs the realized performance of a policy in any given instance and plays a central role in risk assessment and scenario analysis \citep[see also][]{ArlottoGansSteele:OR2014}.
Consequently, two policies with identical asymptotic expected reward may perform quite differently in practice. 
These observations naturally raise the following question: 
can a simple online policy and the optimal dynamic programming policy exhibit the same higher-order probabilistic behavior?

To fix ideas and gain intuition, consider a dynamic and stochastic knapsack problem with
initial capacity $c = 1$ and $T$ arriving items with independent, uniformly distributed weights
$W_1, \ldots, W_T$ on $[0,1]$.
At each time $t$, the decision maker observes the weight $W_t$ of the arriving item and must
immediately accept or reject it.
Accepted items consume capacity equal to their weight and yield unit reward;
rejected items yield nothing and consume no capacity.
The objective is to maximize the expected number of accepted items subject to the capacity
constraint.
For this formulation, the existing literature provides three policies that have been studied extensively:
a nonadaptive (fixed-threshold) policy, 
an adaptive reoptimization heuristic, 
and an optimal dynamic programming policy.
First, the \emph{nonadaptive policy} was introduced by \citet{BrussRobertson:AAP1991,RheeTalagrand:JAP1991}.
If $x$ denotes the currently remaining capacity, the policy accepts the item arriving at time
$t$ if $W_t \leq \widebar h_{T-t+1}(x) \equiv \min\left\{x, \sqrt{\tfrac{2}{T}}\right\}$,
for all $t \in [T] \equiv \{1, 2, \ldots, T\}$.
For this policy, the threshold is constant over time, and the remaining capacity $x$ influences
the acceptance decision only insofar as it ensures that the selected item fits in the knapsack.
Second, the \emph{reoptimized heuristic} follows from what was proposed in \citet{Arlotto-Xie:SSY2020}.
The policy is characterized by the time- and state-dependent thresholds
$\widehat{h}_{T-t+1}(x) = \min\left\{x, \sqrt{\tfrac{2x}{T-t+1}}\right\}$
and accepts the $t$th arriving item if $W_t \leq \widehat{h}_{T-t+1}(x)$.
Third, the \emph{optimal dynamic programming policy} is characterized by a sequence of threshold
functions $\{h^*_{T}, \ldots, h^*_{T-t+1}, \ldots, h^*_{1}\}$, implicitly defined by the dynamic
programming recursion, and accepts the $t$th arriving item if $W_t \leq h^*_{T-t+1}(x)$.

If $\widebar N_T$, $\widehat N_T$, and $N_T^*$ respectively denote the number of selections made
by each of the three policies above, then prior work establishes that all three policies share
the same first-order asymptotic performance:
\[
\E{\widebar N_T}\sim \E{\widehat N_T}\sim \E{N_T^*}\sim \sqrt{2T},
\qquad \text{as } T\to\infty.
\]
Their second-order behavior, however, differs sharply.
\citet{BoshuizenKertz:AAP1999} showed that the number of selections made by the nonadaptive
policy, $\widebar N_T$, has variance
\[
\Var{\widebar N_T}\sim \frac{2(\pi-1)}{3\pi}\sqrt{2T} \approx 0.454 \sqrt{2 T},
\qquad  \text{as $T\to\infty$},
\]
and, after the standard centering and scaling, converges in distribution to the non-Gaussian random
variable $Z_1 + \min\{ Z_2, 0 \}$, where $Z_1$ and $Z_2$ are mean zero jointly normal random
variables with covariance matrix $\bigl(\begin{smallmatrix}1 & -1\\ -1 & 4/3\end{smallmatrix}\bigr)$.
By contrast, the results in this paper \citep[see also][]{ArlottoNguyenSteele:SPA2015} imply that
the number of selections made by the reoptimized heuristic, $\widehat N_T$, and by the optimal
dynamic programming policy, $N^*_T$, have variances satisfying
\[
\Var{\widehat N_T}\sim \Var{N_T^*}\sim \frac{1}{3}\sqrt{2T}, \qquad\text{as $T\to\infty$}.
\]
Furthermore, after centering by $\sqrt{2T}$ and scaling by $(2T)^{1/4}$, both $\widehat{N}_T$ and
$N^*_T$ converge to the same centered normal distribution, $N(0, \tfrac{1}{3})$.
In short, the reoptimized heuristic recovers the second-order and distributional behavior of the
optimal policy, whereas the number of selections made by the nonadaptive policy exhibits a larger
variance and obeys a different limit theorem.

In this paper, we show that this pattern extends well beyond the uniform case.
For the dynamic and stochastic knapsack problem with unit rewards and independent weights drawn from a \emph{general}
nonnegative distribution $F$ satisfying mild regularity conditions near zero, the reoptimized heuristic and
the optimal dynamic programming policy yield total numbers of selections with the same asymptotic mean, 
variance, and limiting distribution.
Methodologically, the paper complements and concludes the analysis initiated with the
logarithmic-regret bound of \citet{Arlotto-Xie:SSY2020}.
They provide a reoptimized heuristic that exhibits logarithmic regret, and we show that a version
of that same heuristic also captures the optimal policy's higher-order stochastic behavior.
In that sense, this paper identifies a heuristic that is not only first-order optimal, but also
distributionally equivalent to the optimal dynamic programming policy.

\subsection*{Organization and notation}

We review related literature in Section~\ref{sec:literature-review}
and introduce the model and reoptimized heuristic in
Section~\ref{sec:DSKP&ReoptimizedHeuristic}.
Section~\ref{sec:main-result} presents the distributional assumptions,
states the main result, and specifies the constants used in the analysis.
Sections~\ref{sec:typical-distribution-properties}
and~\ref{sec:consumption-function-properties} examine the weight
distribution and how it shapes the heuristic's acceptance thresholds.
Section~\ref{sec:vhat-derivative-bounds} establishes value-function derivative bounds. These bounds are used to determine the variance asymptotics in
Section~\ref{sec:variance-asymptotics} and establish the central limit
theorems in Section~\ref{sec:clt}.
All proofs, together with background results on regular variation,
are collected in the appendices.
Throughout the paper, we use \emph{increasing} and \emph{decreasing}
in the weak sense and state explicitly when monotonicity is strict.
For $z \in \R$, we write $(z)_+ = \max\{ z, 0 \}$ for its positive part.

\section{Literature review}
\label{sec:literature-review}

As one of the most classical problems in operations research, the knapsack problem and its variants have been extensively studied in a range of applications and under different setups.
\citet{Dantzig:OR1957} formulated a \emph{static} and \emph{deterministic} model in which the weight-reward information of a finite number of items is available upfront.
The decision maker's goal is to find a maximum-reward subset of such items under a capacity constraint.
When weight and/or reward information is random and is only revealed after inclusion, 
we obtain the static \emph{stochastic} knapsack problems, where the objective is to maximize the expected collected rewards.
For instance, a scheduling problem was studied by \citet{DermanLiebermanRoss:MS1978} through such a formulation in which items have type-dependent deterministic rewards and type-dependent exponentially distributed random weights.
Similar static and stochastic knapsack problems have also been studied by \citet{DeanGoemansVondrak:FOCS2004,DeanGoemansVondrak:SDA2005,DeanGoemansVondrak:MOR2008,BhalgatGoelKhanna:SDA2011,LiYuan:STOC2013,BladoHuToriello:SIAM2016,Ma:MOR2018,BladoToriello:SIAM2019} and \citet{BalseiroBrown:OR2019}.

When items arrive sequentially with randomness in their weight-reward pairs, we find the \emph{dynamic} and stochastic knapsack problems.
This formulation was studied by \citet{PapastavrouRajagopalanKleywegt:MS1996,KleywegtPapastavrou:OR1998,KleywegtPapastavrou:OR2001} where the weight-reward pair is revealed upon each arrival.
The decision maker seeks to maximize the expected collected rewards while making irrevocable inclusion decisions upon each arrival subject to a capacity constraint.
If the items all have unit weights and random rewards, one obtains the multi-secretary problem \citep{Cayley:ET1875,Moser:SM1956,Kleinberg:SDA2005}.
For this problem, \citet{ArlottoGurvich:SSY2019} \citep[see also][]{Vera-Banerjee:MS2021,Vera-Banerjee-Gurvich:OPRE2021,Vera-Arlotto-Gurvich-Levin:MOOR2025} proved a regret bound that is independent of the number of items and the knapsack capacity when the reward distribution is discrete.
On the other hand, if items all have unit rewards and random weights, one finds a single-machine scheduling problem of \citet{BaruahHaritsaSharma:IEEE1994}.
Such formulation is the focus of this paper, and for this formulation, \citet{Arlotto-Xie:SSY2020} \citep[see also][]{Jiang-Zhang-arXiv2020-onlineresourceallocation} proved a logarithmic regret bound under mild regularity conditions on the weight distribution.
Lastly, if both the weights and the rewards of arriving items are random, few regret bounds are available in the literature.
An exception is the special case in which the weights and the rewards are independent and are both uniformly distributed on the unit interval, and the knapsack capacity is proportional to the number of items.
For such formulation, \citet{Marchetti-SpaccamelaVercellis:MP1995} proved a $O(\log^{3/2} T)$ regret bound, and later on \citet{Lueker:JACIL1998} improved such bound to $O(\log T)$ and
showed that this order cannot be improved.

As the focus of this paper, the unit-rewards and random-weights formulation has been considered in \citet[Section 5]{KleywegtPapastavrou:OR1998}.
However, the study of such formulation dates back to the earlier work of \citet{CoffmanFlattoWeber:AAP1987}.
They showed that---under certain conditions on the weight distribution---the number of selections, $N^*_T(c)$, under the optimal dynamic programming policy satisfies the first-order asymptotics
\(
\E{N^*_T(c)} \sim T F(\epsilon_T(c)) \mbox{ as } T \rightarrow \infty,
\)
where the quantity $\epsilon_T(c)$ is defined as
\(
\epsilon_{T}(c) = \sup \left\{ \epsilon \in [0, \infty) : \int_{0}^{\epsilon} w \, d F(w) \leq \frac{c}{T} \right\}.
\)
Later on, \citet{BrussRobertson:AAP1991} generalized such result to more general weight distributions.
In addition, \citet{RheeTalagrand:JAP1991} carried a more delicate probabilistic analysis of the first-order performance of the nonadaptive heuristic.

\section{The dynamic and stochastic knapsack problem with unit rewards and a reoptimized heuristic}\label{sec:DSKP&ReoptimizedHeuristic}

In the dynamic and stochastic knapsack problem with unit rewards,
a decision maker (referred to as \emph{she}) is given a knapsack with capacity $c$
and is sequentially presented with $T$ arriving items with independent weights $W_1,\ldots,W_T$
that are drawn from a nonnegative distribution $F$.
At each period $t \in [T]$ an item arrives and its weight $W_t = w$ is revealed to the decision
maker, who must decide whether to include the item in the knapsack.
If she includes the item, she reduces the remaining knapsack capacity by $w$ and increases the
number of currently selected items by one, earning unit reward.
If she rejects the item, she uses no capacity and leaves the number of currently selected items
unchanged, earning no reward.
All of her decisions are terminal and cannot be changed at a later time.
The decision maker's objective is to maximize the expected number of items she selects,
subject to the knapsack capacity constraint.

This formulation of the dynamic and stochastic knapsack problem has been studied extensively,
and, most recently, \citet{Arlotto-Xie:SSY2020}  proposed a reoptimized heuristic with logarithmic
regret \citep[see also][]{Jiang-Zhang-arXiv2020-onlineresourceallocation}. 
In what follows, we introduce a reoptimized heuristic $\widehat \pi$
that mirrors theirs over all but a fixed number of terminal periods.
Its acceptance thresholds are determined by the optimal solution
of a fluid relaxation of the problem at an arbitrary state,
with remaining capacity $x$ and $k$ remaining periods.
Such an optimal solution, which we refer to as the
\emph{consumption function}, is defined by
\begin{equation}\label{eq:consumption-function-definition}
\epsilon_{k}(x) := \sup \left\{ \epsilon \in [0, \infty) : \int_{0}^{\epsilon} w \, d F(w)
\leq \frac{x}{k} \right\}
\qquad \text{for all $k \geq 1$ and all $x \geq 0$.}
\end{equation}
This consumption function implies that 
if every remaining item with weight at most $\epsilon_k(x)$ is accepted, then the expected capacity consumed by the $k$
remaining arrivals is exactly $x$, so that $\epsilon_k(x)$ spreads the available capacity equally (in expectation) over the remaining arrivals.

Under suitable regularity conditions on $F$, it is easy to see from  \eqref{eq:consumption-function-definition} that $\epsilon_k(x)$ is
continuous and increasing in $x$ for all $x \geq 0$, and decreasing in $k$ for all $k \geq 1$.
Next, if $\mu = \E{W_1}$  then one has the integral representation
\begin{equation}\label{eq:consumption-function-integral-representation}
\int_0^{\epsilon_k(x)} w  \, d F(w) = \frac{x}{k}
\qquad \qquad \text{for all $0 \leq x \leq k\mu$}.
\end{equation}
In turn, since for all $0 \leq x \leq k\mu$ and $j \leq k$ one has
$0 \leq \frac{j x}{k} \leq j \mu$ and
$$
\int_0^{\epsilon_{j}\left( \frac{ j x}{k} \right)} w \, dF(w) = \frac{x}{k},
$$
a comparison with \eqref{eq:consumption-function-integral-representation} gives the space-time
scaling equivalence
\begin{equation}\label{eq:consumption-function-space-time-scaling}
\epsilon_{j}\left( \frac{j x}{k} \right) = \epsilon_k(x)
\quad
\text{ for all $x \in[0, k\mu]$ and $j \leq k$}.
\end{equation}
Finally, if the weight distribution $F$ admits a density $f$ that is continuous and
nonzero at $\epsilon_k(x)$, then the implicit function
theorem applied to \eqref{eq:consumption-function-integral-representation}  tells us that
$\epsilon_k(\cdot)$ is differentiable at $x$, with derivative
\begin{equation}\label{eq:consumption-function-derivative}
\epsilon'_k(x) = \frac{1}{k\epsilon_k(x) f(\epsilon_k(x))}.
\end{equation}

We now use the consumption function \eqref{eq:consumption-function-definition} to build the
reoptimized heuristic $\widehat \pi$.
For a fixed constant $K \geq 1$, the reoptimized
heuristic $\widehat \pi$ is given by the state- and time-dependent thresholds
$\{\hhat_k : 1 \leq k < \infty \}$ defined for $x \in [0, c]$ by
\begin{equation}
\label{eq:threshold-function-definition}
    \hhat_k(x) =
    \begin{cases}
        0 & \text{if } k \leq K \\
        \min(\epsilon_k(x), x) & \text{otherwise}
    \end{cases}.
    \end{equation}
Like the consumption functions $\epsilon_k(\cdot)$, these thresholds
are increasing in the remaining capacity $x$.
For $k\leq K$, the decision maker accepts no items. 
For $k > K$, she acts greedily and selects any item that fits into the knapsack
when $x \leq \epsilon_k(x)$.
In contrast, when $x > \epsilon_k(x)$, she acts conservatively
and selects only items with weight at most $\epsilon_k(x)$.
We refer to these two cases as the \emph{greedy} and
\emph{conservative} regimes, respectively.

Under $\widehat \pi$, the item arriving at time $t \in [T]$ when the remaining capacity is $x$ is selected if and only if
$W_t \leq \hhat_{T-t+1}(x)$.
One can then keep track of how the remaining capacity evolves over time under policy
$\widehat \pi$ by setting $X_0=c$ and letting
\begin{equation}\label{eq:remaining-capacity-process}
X_t =
\begin{cases}
X_{t-1}	\quad 	&\mbox{ if } W_t > \hhat_{T-t+1}(X_{t-1})	\\
X_{t-1}-W_t 	&\mbox{ if } W_t \leq \hhat_{T-t+1}(X_{t-1}),
\end{cases}
\end{equation}
be the recursive definition of the remaining capacity at time $t \in [T]$.
In turn, the number of selections $\widehat N_T(c)$ made by policy $\widehat \pi$ is given by
$$
\widehat N_T(c) = \sum_{t=1}^T \1{W_t \leq \hhat_{T-t+1}(X_{t-1})}.
$$

Let $\vhat_k(x)$ be the value function that denotes the expected number of items selected by $\widehat\pi$ when the remaining capacity is $x$ and there are $k$ remaining periods. 
We set $\vhat_0(x)=0$ for all $x\in[0,c]$. 
For $k\geq1$ and $x\in[0,c]$, the value function satisfies the recursion 
\begin{equation}\label{eq:value-function-recursion}
\vhat_k(x) = 
\begin{cases}
0 & \text{if } k \leq K\\
\{ 1-F(\hhat_k(x)) \} \vhat_{k-1}(x) + \int_0^{\hhat_k(x)} \left\{ 1+\vhat_{k-1}(x-w) \right\}  \, dF(w)  & \text{if } k > K
\end{cases}.
\end{equation}
As usual, we have that $\vhat_T(c) = \E{\widehat N_T(c)}$, and we note that $\vhat_k$ is increasing in $x$ for each $k \geq 1$.

For $k > K$, suppose that the weight distribution $F$ has density $f$.
Under suitable regularity conditions, differentiating
\eqref{eq:value-function-recursion} gives the derivative recursion
\begin{align}
\label{eq:vhat-derivative-recursion}
    \vhat_k'(x) =
    \begin{cases}
        \left( 1 - F(x) \right) \vhat_{k-1}'(x)
        + \displaystyle \int_0^x \vhat_{k-1}'(x-w) f(w) \, dw
        + f(x) \left( 1 - \vhat_{k-1}(x) \right)
        & \text{if } \epsilon_k(x) > x \\[12pt]
        \begin{aligned}[t]
            & \left( 1 - F(\epsilon_k(x)) \right) \vhat_{k-1}'(x)
            + \int_0^{\epsilon_k(x)}
                \vhat_{k-1}'(x-w) f(w) \, dw \\
            & \quad
            + \frac{ 1 }{ k \epsilon_k(x) }
            \left\{
                1 + \vhat_{k-1}(x-\epsilon_k(x))
                - \vhat_{k-1}(x)
            \right\}
        \end{aligned}
        & \text{if } \epsilon_k(x) < x
    \end{cases}.
\end{align}
The first branch corresponds to the greedy regime, where
$\hhat_k(x) = x$, and the second to the conservative regime,
where $\hhat_k(x) = \epsilon_k(x)$.
At a switching point where $x = \epsilon_k(x)$, the derivative
need not exist, so we leave this case unspecified in
\eqref{eq:vhat-derivative-recursion}.
In Section~\ref{sec:vhat-derivative-bounds}, we first justify
this differentiation under the distributional assumptions used
in this paper and show that $\vhat_k'$ exists and satisfies
\eqref{eq:vhat-derivative-recursion} at all but finitely many
points of $(0,c)$.
We then use the derivative
recursion to show that $\vhat_k'(x)$ is well approximated by
the reciprocal consumption function $1 / \epsilon_{k+1}(x)$,
with an explicit bound on the approximation error.
This approximation allows us to control the value-function differences
caused by selecting an item, a crucial step towards determining
the variance of the total number $\widehat N_T(c)$ of selections
made by the reoptimized heuristic.

\section{Main result}\label{sec:main-result}

In our analysis, we study a knapsack regime in which the number of items $T$ grows large while
the knapsack capacity $c$ stays the same. As such, any reasonable policy becomes increasingly
selective and accepts only items with progressively smaller weights, so the behavior of the
weight distribution $F$ near zero becomes essential.
The weight distributions we consider in this paper are said to belong to the \emph{typical
class}, whose defining properties are detailed below. 
In essence, we require a monotonicity property near zero, 
together with the convergence of the inverse elasticity of $F$ to a
finite limit at a quantified rate.
These two properties carry several important implications that connect the typical class with the
theory of regularly varying functions; we develop these connections in Section
\ref{sec:typical-distribution-properties}.

\begin{definition}[Typical class of distributions]\label{def:typical-class}
A nonnegative distribution $F$ with density $f$ 
is said to belong to the \emph{typical class} if there are $0 < r, \wbar < \infty$ and $0 \leq M_F < \infty$
such that 
\begin{itemize}
\item For all $w \in (0,\wbar]$, the density $f(w)>0$, the derivative $f'(w)$ exists, and 
\begin{equation}\label{eq:typical-class-monotonicity-condition}
wf'(w) \geq -f(w).
\end{equation}
\item For all $w\in (0, \wbar]$, 
\begin{equation}\label{eq:typical-class-limit-with-rate}
	\abs{\frac{F(w)}{w f(w)}  - \frac{1}{r}} \leq M_F w^{3/4}.
\end{equation}
\end{itemize}
\end{definition}

The typical class of Definition~\ref{def:typical-class} is closely related to the classes
of weight distributions considered in earlier work on this problem.
From \eqref{eq:typical-class-limit-with-rate}, one obtains that
$ \left( \log F(w) \right)' = \frac{f(w)}{F(w)} = \frac{r}{w} + O(w^{-1/4})$ as $w \to 0$, 
which after integration and exponentiation gives $F(w) = A w^{r} \left\{ 1 + O(w^{3/4}) \right\}$ as $w \to 0$ for some constant $A>0$.
Thus, every typical distribution is asymptotic to a power law as $w \to 0$, which is the
property assumed by \citet{CoffmanFlattoWeber:AAP1987}.
This power-law behavior in turn gives $\limsup_{w \to 0} \frac{F(\lambda w)}{F(w)} < 1$ for some
$\lambda \in (0,1)$, so every typical distribution also satisfies the limsup condition
first required by \citet{BrussRobertson:AAP1991} and later by \citet{Arlotto-Xie:SSY2020}.
This power-law behavior also implies the two-sided condition assumed by
\citet{Jiang-Zhang-arXiv2020-onlineresourceallocation}.
Both \citet{Arlotto-Xie:SSY2020} and \citet{Jiang-Zhang-arXiv2020-onlineresourceallocation}
also impose a monotonicity requirement, which is implied by
\eqref{eq:typical-class-monotonicity-condition}.

The preceding discussion shows that the typical class of
Definition~\ref{def:typical-class} satisfies the near-zero conditions
used by \citet{Arlotto-Xie:SSY2020}.
Their logarithmic regret bound therefore carries over to our reoptimized
heuristic: when run on the same sequence of item weights, the two heuristics
use the same thresholds and make the same selections whenever more than
$K$ periods remain.
Since our heuristic rejects each of the last $K$ items,
it makes at most $K$ fewer selections overall.
This additional loss can be absorbed into the constant in the logarithmic
regret bound, yielding the following proposition.

\begin{proposition}[Logarithmic regret bound, {\citealp[Theorem 1]{Arlotto-Xie:SSY2020}}]
\label{prop:log-regret-bound}
Suppose that $F$ belongs to the typical class of
Definition~\ref{def:typical-class}, and 
let $N_k^{\mathrm{off}}(x)$ denote the maximum number of items
that can be selected from $k$ arrivals with capacity $x$ in hindsight.
Then there exists a constant $C < \infty$ such that
\[
    k F(\epsilon_k(x)) - C(1 + \log k)
    \leq \E{ \widehat{N}_k(x) }
    \leq \E{ N_k^{\mathrm{off}}(x) \vphantom{\widehat N} }
    \leq k F(\epsilon_k(x))
\]
for every integer $k \geq 1$ and every $x \in [0,c]$.
\end{proposition}

We are now in a position to state the main result of the paper. In a dynamic and stochastic
knapsack problem with unit rewards in which the weight distribution belongs to the typical class
of Definition~\ref{def:typical-class}, the number $\widehat N_T(c)$ of selections made by the reoptimized
heuristic $\widehat \pi$ of Section~\ref{sec:DSKP&ReoptimizedHeuristic} and the number $N^*_T(c)$ of selections made by the
optimal dynamic programming policy $\pi^*$ obey the same central limit theorem.
Both numbers of selections are centered at $T F(\epsilon_T(c))$, the leading term of their
expected performance, which is also equal to the expected number of items whose weight falls below
the initial consumption function $\epsilon_T(c)$.

\begin{theorem}[Central limit theorems]
\label{thm:main-theorem}
Consider the dynamic and stochastic knapsack problem with unit rewards, initial capacity $c$, 
weight distribution $F$ that belongs to the typical class of Definition
\ref{def:typical-class}, and set
\[
\sigma^2 \equiv \frac{2}{(r+1)(r+2)}.
\]
If $\widehat N_T(c)$ and $N^*_T(c)$ denote the numbers of selections made by the reoptimized
heuristic $\widehat \pi$ and by the optimal dynamic programming policy $\pi^*$, then as $T \rightarrow \infty$ one has the convergences in distribution
$$
\frac{\widehat N_T(c) - T F(\epsilon_T(c))}{ \sqrt{ T F(\epsilon_T(c)) } } \rightarrow N(0,\sigma^2)
\qquad \mbox{ and } \qquad
\frac{N^*_T(c) - T F(\epsilon_T(c))}{ \sqrt{ T F(\epsilon_T(c)) } } \rightarrow N(0,\sigma^2).
$$
\end{theorem}

To establish Theorem~\ref{thm:main-theorem}, we proceed in four steps.
First, in Sections~\ref{sec:typical-distribution-properties}
and~\ref{sec:consumption-function-properties}, we develop key properties
of the weight distributions in Definition~\ref{def:typical-class}
and the consumption function $\epsilon_k(x)$, respectively.
Second, in Section~\ref{sec:vhat-derivative-bounds}, we approximate
the derivative $\vhat_k'(x)$ by the reciprocal consumption function
$1 / \epsilon_{k+1}(x)$ and establish a non-asymptotic bound
on the approximation error.
Third, in Section~\ref{sec:variance-asymptotics}, we combine the resulting
estimates of value-function differences with Doob's martingale
representation of $\widehat N_T(c)$ to derive its variance asymptotics
and control the accumulated approximation errors.
Finally, in Section~\ref{sec:clt}, we verify the conditions of a martingale
central limit theorem to establish the normal limit for the number of items selected by the reoptimized heuristic.
We then use the logarithmic regret and full-information bounds of Proposition~\ref{prop:log-regret-bound}
to transfer this limit to the optimal policy.

\subsection{Choice of constants}
\label{sec:choice-of-constants}

As discussed in Section~\ref{sec:DSKP&ReoptimizedHeuristic},
the reoptimized heuristic $\widehat \pi$ rejects all items when
at most $K$ periods remain.
We now specify this cutoff and record the resulting bounds on
the consumption function.
Our choice keeps the relevant consumption thresholds within the
interval on which the typical-class assumptions apply.
With the capacity $c$ fixed, the choice of $K$ depends on the weight
distribution $F$ and an auxiliary parameter $\Theta > 1$,
whose value will be chosen sufficiently large in
Section~\ref{sec:vhat-derivative-bounds}.
Set
\begin{equation}
\label{eq:K1}
    K_1
    =
    \ceil{ \frac{ \Theta c }{ \int_0^{\wbar} w f(w) \, dw } }
    \geq \frac{ \Theta c }{ \mu },
\end{equation}
where $\wbar$ is as in Definition~\ref{def:typical-class}.
The integral representation
\eqref{eq:consumption-function-integral-representation} gives
\[
    \int_0^{\epsilon_{K_1}(\Theta c)} w f(w) \, dw
    =
    \frac{ \Theta c }{ K_1 }
    \leq
    \int_0^{\wbar} w f(w) \, dw,
\]
so $\epsilon_{K_1}(\Theta c) \leq \wbar$.
By the monotonicity of the consumption function, it follows that
\begin{equation}
\label{eq:epsilon-k-leq-wbar}
    \epsilon_k(x)
    \leq \epsilon_k(\Theta x)
    \leq \epsilon_{K_1}(\Theta c)
    \leq \wbar,
    \qquad 
    \text{for every $k \geq K_1$ and $x \in [0,c]$}.
\end{equation}
In particular, for every $k \geq K_1$, the integral representation
\eqref{eq:consumption-function-integral-representation} holds for
capacity arguments in $[0,\Theta c]$, and the derivative identity
\eqref{eq:consumption-function-derivative} holds on $(0,\Theta c)$.
We now choose $K$ to be a fixed integer satisfying
\begin{equation}
\label{eq:K}
    K \geq K_1
    \qquad \text{and} \qquad
    \epsilon_K(\Theta c) < c.
\end{equation}
Such a choice exists by the integral representation
\eqref{eq:consumption-function-integral-representation}
and the positivity of $f$ on $(0,\wbar]$.
Since the consumption function is decreasing in $k$, we also have
$\epsilon_k(\Theta c) < c$ for every $k \geq K$.
This additional restriction will allow us to partition $[0,c]$
into a small-capacity region adjacent to zero and a large-capacity
region containing $c$.
The small-capacity region shrinks as $k$ grows.
There, bounds motivated by the greedy regime remain valid even
when the heuristic makes conservative selections.
In the large-capacity region, the decision maker makes conservative
selections, and we use this property to derive the corresponding bounds.
Throughout the remainder of the paper, $K$ denotes this standing
choice, which is fixed independently of the horizon length $T$.

\section{Typical distributions and their properties}
\label{sec:typical-distribution-properties}

In this section, we develop several consequences of the typical-class
assumptions in Definition~\ref{def:typical-class}.
We first show how the constant $\sigma^2$ arises when integrating
the square of a certain linear function against a typical distribution,
a calculation used later in our variance analysis.
Then, we connect typical distributions with the theory of regularly
varying functions and derive integral estimates used repeatedly
in the rest of our analysis.
The proofs of the results in this section are given in
Appendix~\ref{apx:typical-distribution-properties-proofs}.

\subsection{Quadratic integrals of typical distributions and the constant $\sigma^2$}
\label{sec:sigma}

The following lemma gives an integral characterization of the constant
$\sigma^2$ appearing in Theorem~\ref{thm:main-theorem}.
It shows that integrating the squared linear function
$\left( 1 - \frac{ w }{ h } \right)^2$ against a typical distribution
over $[0,h]$ yields $\sigma^2 F(h)$ up to an explicitly controlled error.
In Section~\ref{sec:variance-asymptotics}, linear approximations of
value-function differences will lead to integrals of this form.
The factor $F(h)$ will then represent a per-period selection probability, and
$\sigma^2 F(h)$ will determine the per-period variance contribution.

\begin{lemma}[Bounds for quadratic integrals against typical distributions]
\label{lem:linear-approximation-second-moment-bounds}
Suppose that $F$ belongs to the typical class of
Definition~\ref{def:typical-class} with parameter $r$ and set
$\sigma^2 = \frac{ 2 }{ (r+1)(r+2) }$.
For every $h \in (0, \wbar]$,
\begin{equation}
\label{eq:linear-approximation-second-moment-bounds}
    \sigma^2 F(h) - M_F h^{3/4} F(h)
    \leq
    \int_0^h
        \left( 1 - \frac{ w }{ h } \right)^2 f(w) \, dw
    \leq
    \sigma^2 F(h) + M_F h^{3/4} F(h),
\end{equation}
where $M_F$ is the constant in Definition~\ref{def:typical-class}.
\end{lemma}

\subsection{Typical distributions as regularly varying functions}
\label{sec:typical-distributions-are-regularly-varying}

Definition~\ref{def:typical-class} implies several properties of typical distributions that will be used throughout this paper.
In particular, we find that typical densities are regularly varying at zero. 
This property will be invoked repeatedly and will simplify a number of computations.
It also allows us to draw on the theory of slowly and regularly varying functions, along with several related results that will be relevant to our analysis (see Appendix~\ref{apx:regularly-varying}).

\begin{definition}[Slowly  and regularly varying functions at zero]\label{def:slowly varying}
A continuous and positive function $\ell$, defined on an interval $(0, U) \subseteq \R_+$, 
is \emph{slowly varying} at zero if
$$
\lim_{u \rightarrow 0^+} \frac{\ell(\lambda u)}{\ell(u)} = 1 
\qquad \mbox{ for all } \lambda > 0.
$$
A continuous and positive function $g$, defined on an interval $(0,U) \subseteq \R_+$ 
is \emph{regularly varying} at zero with index $\rho$ 
if $g(u) = u^{\rho} \ell(u)$ for some slowly varying function $\ell(u)$. 
\end{definition}

Intuitively, a regularly varying function behaves like a power law near zero: 
multiplying the argument by a constant $\lambda$ multiplies the function value by approximately $\lambda^{\rho}$. 
The index $\rho$ thus captures the dominant power-law behavior, while the slowly varying factor $\ell$ collects the subtler corrections that remain once that power is factored out.

The following proposition connects the typical class to the theory of regularly varying 
functions: the two conditions of Definition~\ref{def:typical-class} force the density $f$ 
and the weight distribution $F$ to be regularly varying at zero, with indices $r-1$ and $r$ respectively.

\begin{proposition}[Regularity properties of typical distributions]
\label{prop:typical-density-properties}
If $F$ is a weight distribution that belongs to the typical class of Definition~\ref{def:typical-class} and $f$ is its density, then: 

\begin{enumerate}[label=\textnormal{(\roman*)}]

    \item \label{property:wfw-increasing}
        The map $w\mapsto wf(w)$ is increasing on $(0,\wbar]$.  
    
    \item \label{property:F-regularly-varying} 
        The weight distribution $F$ is regularly varying at zero with index $r$.
    
    \item \label{property:f-regularly-varying}
        The density $f$ is regularly varying at zero with index $r-1$; 
        that is, $f(w) = w^{r-1}\ell(w)$ for some slowly varying function $\ell$.

\end{enumerate}
\end{proposition}

Since $F$ and $f$ are regularly varying at zero, Karamata's theorem
(Lemma~\ref{lm:Karamata-RV}) gives explicit limits for normalized
integral ratios.
The following proposition quantifies this convergence for weighted
integrals of $F$, with error bounded by $M_F w^{3/4}$,
where $M_F$ is the constant in
\eqref{eq:typical-class-limit-with-rate}.

\begin{proposition}[Weighted integral bounds for typical distributions]
\label{prop:convergence-rate-bound}
Suppose that $F$ belongs to the typical class of
Definition~\ref{def:typical-class} with parameter $r$.
Then, for every real $q \geq -1$ and every
$w\in(0,\wbar]$,
\begin{equation}\label{eq:error-bound-of-limits}
\abs{
\frac{\int_0^w y^q F(y)\,dy}{w^{q+1}F(w)}
-\frac{1}{r+q+1}}
\leq M_Fw^{3/4}.
\end{equation}
\end{proposition}

The truncated first moment $\int_0^\epsilon w f(w) \, dw$
determines the consumption function through
\eqref{eq:consumption-function-integral-representation}.
The following corollary uses the preceding proposition to compare
this moment with $\epsilon^2 f(\epsilon)$ and $\epsilon F(\epsilon)$,
identifying the limiting ratios and providing bounds valid throughout
$(0,\wbar]$.
These uniform comparisons, expressed in terms of a common constant $C_F$,
will be used repeatedly to control the consumption function and its
derivative and to relate expected capacity consumption to selection
probabilities.

\begin{corollary}[Integral ratios for typical distributions]
\label{cor:karamata-ratios}
Suppose that $F$ belongs to the typical class of
Definition~\ref{def:typical-class}, with density $f$
and parameter $r$. Then
\begin{equation}
\label{eq:karamata-ratios-limit}
    \lim_{\epsilon \to 0^+}
    \frac{ \int_0^\epsilon w f(w) \, dw }{ \epsilon^2 f(\epsilon) }
    = \frac{ 1 }{ r+1 }
    \qquad \text{and} \qquad
    \lim_{\epsilon \to 0^+}
    \frac{ \int_0^\epsilon w f(w) \, dw }{ \epsilon F(\epsilon) }
    = \frac{ r }{ r+1 }.
\end{equation}
As a result, there exists a constant $C_F > 1$, depending only on $F$
and $\wbar$, such that
\begin{equation}
\label{eq:karamata-ratios-pointwise-bounds}
    \frac{ 1 }{ C_F }
    \leq
    \frac{ \int_0^\epsilon w f(w) \, dw }{ \epsilon^2 f(\epsilon) }
    \leq 1
    \qquad \text{and} \qquad
    \frac{ 1 }{ C_F }
    \leq
    \frac{ \int_0^\epsilon w f(w) \, dw }{ \epsilon F(\epsilon) }
    \leq 1,
\end{equation}
for all $0 < \epsilon \leq \wbar$.
Furthermore, for all $0 < h < \epsilon \leq \wbar$,
\begin{equation}
\label{eq:karamata-integral-ratio-bounds}
    \left( \frac{ h }{ \epsilon } \right)^{C_F}
    \leq
    \frac{ \int_0^h w f(w) \, dw }{ \int_0^\epsilon w f(w) \, dw }
    \leq
    \frac{ h }{ \epsilon }.
\end{equation}
\end{corollary}

\section{Analytical properties of the consumption function}
\label{sec:consumption-function-properties}

This section develops the properties of the consumption function
needed for the derivative analysis in
Section~\ref{sec:vhat-derivative-bounds}.
We study how $\epsilon_k(x)$ and its reciprocal change with the
remaining capacity and horizon, using the estimates
of Section~\ref{sec:typical-distribution-properties}.
We first establish basic properties of $\epsilon_k$, then derive
estimates that also involve $F$ and its density $f$.
We subsequently introduce the auxiliary boundaries $\xhat_k$,
which separate the small- and large-capacity regions used in the
derivative analysis, and establish bounds suited to each region.
The following proposition begins this journey by collecting basic analytical properties
of $\epsilon_k$, including its concavity and regular variation at zero.
The proofs of the results in this section are given in
Appendix~\ref{apx:consumption-function-properties-proofs}.

\begin{proposition}[Basic properties of the consumption function]
\label{prop:basic-properties-epsilon}
If $F$ is a weight distribution that belongs to the typical class of Definition
\ref{def:typical-class} and $\epsilon_k$ is the consumption function defined by
\eqref{eq:consumption-function-definition}, then the following properties hold.

\begin{enumerate}[label=\textnormal{(\roman*)}]

\item\label{property:consumption-function-concavity}
For all $k \geq K_1$, the map $x \mapsto \epsilon_k(x)$ is concave on $[0,c]$.
Consequently, the map $x \mapsto \epsilon_k(x)/x$ is decreasing on $(0,c]$.

\item\label{property:consumption-function-regularly-varying}
For all $k \geq K_1$, the map $x \mapsto \epsilon_k(x)$ is regularly varying at zero with index
$\frac{1}{r+1}$.  Furthermore, there is a function $L$ that is slowly varying at zero such that
\[
\epsilon_k(x) 
= \epsilon_{K_1}\left(\frac{K_1 x}{k}\right)
= \left(\frac{x}{k}\right)^{\frac{1}{r+1}} L\left(\frac{x}{k}\right)
\qquad \text{for all $k \geq K_1$ and all $x \in (0,c]$}.
\]

\item \label{property:useful-eps'/eps^2-decreasing}
For any $a > 0$ and any $k \geq K_1$, the map $x \mapsto \frac{\epsilon_k'(x)}{\epsilon_k(x)^a}$
is decreasing on $(0,c)$.

\item \label{property:ratio-eps-x-to-x-minus-eps}
For all $k \geq K_1$ and $0 \leq w \leq \epsilon_k(x) < x \leq c$, we have
    \begin{equation}\label{eq:ratio-eps-x-to-x-minus-eps}
        \frac{ \epsilon_k(x) }{ \epsilon_k(x-w) }
        \leq \frac{ x }{ x-w }
        \leq \frac{ x }{ x - \epsilon_k(x) }.
    \end{equation}
In particular,
    \begin{equation}\label{eq:ratio-eps-x-to-x-minus-eps-3/4}
        \frac{ 1 }{ (x-w)^{3/4} \epsilon_k(x-w)^{1/4} }
        \leq
        \frac{ x }{ x-w }
        \frac{ 1 }{ x^{3/4} \epsilon_k(x)^{1/4} }
        \leq
        \frac{ x }{ x - \epsilon_k(x) }
        \frac{ 1 }{ x^{3/4} \epsilon_k(x)^{1/4} }.
    \end{equation}

\item\label{property:1/epsilon-integral-bounds}
For all $k \geq K_1$ and all $w,x$ such that $0 \leq w \leq  \hhat_k(x)$,
\begin{equation}\label{eq:1/epsilon-integral-bounds}
0 \leq \int_{x-w}^x \frac{1}{\epsilon_k(y)} \, dy - \frac{w}{\epsilon_k(x)}
\leq \frac{w^2}{(x-w)\epsilon_k(x)}.
\end{equation}

\item\label{property:taylor-first-order-eps-a}  
For all $a > 0$, $k \geq K_1$, $x \in (0,c]$, and $w \in [0,x)$,
    \begin{equation}
    \label{eq:taylor-first-order-eps-a}
    \frac{\epsilon_k'(x)}{\epsilon_k(x)^{a+1}}  aw 
    \leq 
    \frac{1}{\epsilon_k(x-w)^a} - \frac{1}{\epsilon_k(x)^a}
    \leq 
    \frac{\epsilon_k'(x-w)}{\epsilon_k(x-w)^{a+1}} aw .
    \end{equation}

\end{enumerate}
\end{proposition}

We now combine the analytical properties in Proposition~\ref{prop:basic-properties-epsilon} with the distributional
estimates of Section~\ref{sec:typical-distribution-properties}.
The following proposition establishes concavity of the selection
probability $F(\epsilon_k(x))$ and gives quantitative bounds on the
derivative and reciprocal increments of the consumption function.
It also shows that the fluid benchmark grows sufficiently quickly
for the logarithmic regret to be negligible on the scales used
in the variance and central limit analyses.

\begin{proposition}[Joint properties of the weight distribution and the consumption function]
\label{prop:joint-properties-f-epsilon}
If $F$ is a weight distribution that belongs to the typical class 
of Definition~\ref{def:typical-class} and $f$ is its density,
then the following properties hold.
\begin{enumerate}[label=\textnormal{(\roman*)}]

    \item\label{property:fluid-benchmark-growth-relative-to-logT}
    For the fixed capacity $c > 0$,
    \(
        \lim_{T \to \infty}
        \frac{ \log T }{ \{ T F(\epsilon_T(c)) \}^{1/4} }
        = 0.
    \)

    \item \label{property:concavity-of-power-of-F-epsilon}
    For all
    $k \geq K_1$, the map $x \mapsto F( \epsilon_k(x) )$
    is concave on $(0,c]$.
    Moreover, the map
    $x \mapsto F( \epsilon_k(x) )/x$ is decreasing on $(0,c]$.

    \item \label{property:epsilon'/epsilon^2-delta-bounds}
    For all $a > 0$, $k \geq K_1$, and $x \in (0,c]$,
    $$
    \frac{1}{C_F} \frac{1}{x\epsilon_k(x)^a}
    \leq \frac{\epsilon_k'(x)}{\epsilon_k(x)^{a+1}}
    \leq \frac{1}{x\epsilon_k(x)^a}.
    $$

    \item \label{property:local-reciprocal-increments}
    Let
    \(
        C_\epsilon = 3 + 2 M_F c^{3/4}.
    \)
    For every $k \geq K_1$ and every $x \in (0,c]$ satisfying
    $\epsilon_k(x) \leq x/2$, the following bound holds
    for all $w \in [0, \epsilon_k(x)]$:
    \begin{equation}
    \label{eq:difference-of-consumption-function-reciprocals-linear-approximation}
        \frac{ \epsilon_k'(x) }{ \epsilon_k(x)^2 } w
        \leq
        \frac{ 1 }{ \epsilon_k(x-w) }
        - \frac{ 1 }{ \epsilon_k(x) }
        \leq
        \left(
            \frac{ \epsilon_k'(x) }{ \epsilon_k(x)^2 }
            + \frac{ C_\epsilon }{ x^{7/4} \epsilon_k(x)^{1/4} }
        \right) w.
    \end{equation}
    In particular,
    \begin{equation}
    \label{eq:difference-of-consumption-function-reciprocals-dFintegral}
            \frac{ \epsilon_k'(x) }{ \epsilon_k(x)^2 }
            \frac{ x }{ k }
            \leq
            \int_0^{\epsilon_k(x)}
            \left(
                \frac{ 1 }{ \epsilon_k(x-w) }
                - \frac{ 1 }{ \epsilon_k(x) }
            \right) f(w) \, dw 
            \leq
            \frac{ \epsilon_k'(x) }{ \epsilon_k(x)^2 }
            \frac{ x }{ k }
            + \frac{ C_\epsilon }{ k x^{3/4} \epsilon_k(x)^{1/4} }.
    \end{equation}

\end{enumerate}
\end{proposition}

For $k > K$, the reoptimized heuristic $\widehat \pi$ switches
from greedy to conservative selections as the remaining capacity
increases.
The continuation terms in both
\eqref{eq:value-function-recursion} and
\eqref{eq:vhat-derivative-recursion} involve smaller capacity levels.
Consequently, the analysis of conservative selections at small
capacities must also account for the greedy behavior at lower capacities.
We therefore introduce an auxiliary \emph{small-capacity boundary}
to distinguish two parts of the derivative analysis.
With $\Theta > 1$ as introduced in
Section~\ref{sec:choice-of-constants}, define
\begin{equation}
\label{eq:xk-definition}
    \xhat_k
    =
    \max\{ x \in [0,c]: x \leq \epsilon_k(\Theta x) \}
    \qquad \text{for all } k \geq 1.
\end{equation}
Equivalently, $\xhat_k$ is the largest $x \in [0,c]$ satisfying
$\int_0^x w f(w) \, dw \leq \Theta x/k$.
Since this integral is continuous and is $o(x)$ as $x \to 0^+$,
the maximum exists and is positive.

For $k \geq K$, the standing choice \eqref{eq:K} gives
$\epsilon_k(\Theta c) < c$, and hence $0 < \xhat_k < c$.
Continuity of the consumption function then implies that the defining
inequality holds with equality at $\xhat_k$:
\begin{equation}
\label{eq:xk-fixedpoint-property}
    \xhat_k = \epsilon_k(\Theta \xhat_k)
    \qquad \text{for all } k \geq K.
\end{equation}
For $k > K$, the small-capacity region $(0,\xhat_k)$ contains
both greedy and conservative selections, while the large-capacity
region $[\xhat_k,c]$ lies entirely in the conservative regime.
The following proposition shows that the small-capacity boundaries
decrease to zero and provides bounds on $\epsilon_k(x)/x$
on either side of $\xhat_k$.
These bounds quantify how the parameter $\Theta$ controls the
relative size of the consumption threshold in the two regions.

\begin{proposition}[Properties of $\xhat_k$]
\label{prop:properties-of-xhat}
Suppose that $F$ belongs to the typical class of
Definition~\ref{def:typical-class}, fix $\Theta > 1$, and let
$\xhat_k$ be the small-capacity boundary defined by \eqref{eq:xk-definition}.
Then the following properties hold.
\begin{enumerate}[label=\textnormal{(\roman*)}]

    \item \label{property:fixed-points-decreasing}
    The small-capacity boundaries $\xhat_k$ are decreasing in $k$ for all
    $k \geq K_1$, and $\xhat_k \to 0$ as $k \to \infty$.

    \item \label{property:consumption-bound-below-xhat}
    For every $k \geq K_1$ and every $x \in [0, \xhat_k]$,
    \begin{equation}
    \label{eq:consumption-bound-below-xhat}
        \int_0^x w f(w) \, dw
        \leq \frac{ \Theta x }{ k }
        \qquad \text{and} \qquad
        x \leq \epsilon_k(\Theta x).
    \end{equation}
    
    \item \label{property:epsilon/x-ratio-bounds}
    Let $C_F > 1$ be the constant in Corollary~\ref{cor:karamata-ratios},
    and set
    \[
        \thetalb = \frac{ 1 }{ \Theta }
        \qquad \text{and} \qquad
        \thetaub = \left( \frac{ 1 }{ \Theta } \right)^{1 / C_F}.
    \]
    Then, for every $k \geq K$,
    \[
        \frac{ \epsilon_k(x) }{ x } \geq \thetalb
        \quad \text{for } x \in (0, \xhat_k),
        \qquad \text{and} \qquad
        \frac{ \epsilon_k(x) }{ x } \leq \thetaub
        \quad \text{for } x \in [\xhat_k, c].
    \]
    Furthermore, $\thetalb < \thetaub < 1/2$ whenever $\Theta > 2^{C_F}$.

\end{enumerate}
\end{proposition}

In the large-capacity region $x \geq \xhat_k$, the bound
$\epsilon_k(x)/x \leq \thetaub$ allows us to control the coefficient
$f(w) - 1 / (k \epsilon_k(x))$.
The next two lemmas first show that this coefficient is nonnegative
outside an interval adjacent to zero, and then bound a weighted
integral of its absolute value.
This integral estimate will control the propagation of approximation
errors through the derivative recursion in
Section~\ref{sec:vhat-derivative-bounds}.

\begin{lemma}[Density lower bound with large capacity]
\label{lemma:density-lower-bound-large-capacity}
Suppose that $F$ belongs to the typical class of
Definition~\ref{def:typical-class} with parameter $r$.
Fix $\Theta > 1$, 
%and let $K$, $C_F$, and $\thetaub$ be as in
% Proposition~\ref{prop:properties-of-xhat}--\ref{property:epsilon/x-ratio-bounds}.
and set
\(
    \lambda_\Theta
    = \thetaub^{1/(C_F - 1)}.
\)
Then, for every $k \geq K$ and $x \in [\xhat_k, c]$,
\begin{equation}
\label{eq:density-lower-bound-large-capacity}
    \frac{ 1 }{ k \epsilon_k(x) }
    \leq f(w),
    \qquad \text{for all }
    \lambda_\Theta \epsilon_k(x)
    \leq w \leq \epsilon_k(x).
\end{equation}
Consequently,
\begin{equation}
\label{eq:positive-part-integral-bound}
    \int_0^{\epsilon_k(x)}
    \left(
        \frac{ 1 }{ k \epsilon_k(x) } - f(w)
    \right)_+
    \, dw
    \leq
    \frac{ \lambda_\Theta }{ k }.
\end{equation}
\end{lemma}

\begin{lemma}[Absolute-coefficient estimate]
\label{lemma:absolute-coefficient-estimate}
Suppose that $F$ belongs to the typical class of
Definition~\ref{def:typical-class} with parameter $r$.
Fix $\Theta$ sufficiently large that $\thetaub < 1/2$,
and set
\[
    C_\Theta
    =
    \frac{ \thetaub + 2 \lambda_\Theta }{ 1 - \thetaub }
    % =
    % \frac{ \thetaub + 2 \thetaub^{1/(C_F - 1)} }{ 1 - \thetaub }.
\]
Then, for every $k \geq K$ and $x \in (\xhat_k, c)$,
\begin{equation}
\label{eq:absolute-coefficient-estimate}
    \int_0^{\epsilon_k(x)}
    \frac{ 1 }{ (x - w)^{3/4} \epsilon_k(x - w)^{1/4} }
    \abs{ \frac{1}{k \epsilon_k(x)} - f(w) }
    \, dw
    \leq
    \left(F(\epsilon_k(x)) + \frac{ C_\Theta }{ k } \right) 
    \frac{ 1 }{ x^{3/4} \epsilon_k(x)^{1/4} }.
\end{equation}
\end{lemma}

We now turn to the small-capacity region $x \leq \xhat_k$.
Since every selected item must have weight at most $x$, the expected
number of selections is bounded by $k F(x)$.
The following lemma uses the defining property of $\xhat_k$ to bound
this quantity uniformly, providing the value-function estimate needed
for the small-capacity derivative analysis.

\begin{lemma}
\label{lemma:vhat-upperbound-small-capacity}
Suppose that $F$ belongs to the typical class of
Definition~\ref{def:typical-class}.
Fix $\Theta > 1$ and set $M_v = C_F \Theta$, where $C_F$
is the constant in Corollary~\ref{cor:karamata-ratios}.
Then, for every $k \geq K_1$,
\begin{equation}
\label{eq:vhat-upperbound-small-capacity}
    0
    \leq \vhat_k(x)
    \leq k F(x)
    \leq M_v,
    \qquad \text{for all } 0 \leq x \leq \xhat_k.
\end{equation}
\end{lemma}

\section{Non-asymptotic approximation of the value-function derivative}
\label{sec:vhat-derivative-bounds}

This section establishes a quantitative approximation for the derivative
of the value function associated with the reoptimized heuristic $\widehat \pi$.
Proposition~\ref{prop:vhat-derivative-bounds} shows that the marginal value
of remaining capacity, $\vhat_k'(x)$, is approximated by the reciprocal
consumption function $1 / \epsilon_{k+1}(x)$.
It gives a non-asymptotic bound on the approximation error, with an explicit
rate expressed in terms of the remaining capacity $x$ and the consumption
function $\epsilon_{k+1}(x)$.
This estimate will allow us to control the value-function increments
in the martingale decomposition of the total reward collected by policy
$\widehat \pi$, and thereby derive its variance asymptotics
(Section~\ref{sec:variance-asymptotics}) and prove the central limit theorem
(Section~\ref{sec:clt}).
The proofs of the results in this section are given in
Appendix~\ref{apx:vhat-derivative-bounds-proofs}.

We begin by establishing the regularity needed to justify the derivative
recursion \eqref{eq:vhat-derivative-recursion}.
The lemma below records the continuity and differentiability
properties of the value functions $\{ \vhat_k: 0 \leq k < \infty \}$.

\begin{lemma}[Absolute continuity and differentiability of the value functions]
\label{lem:vhat-continuity&differentiability}
Suppose that the weight distribution $F$ belongs to the typical class
of Definition~\ref{def:typical-class}.
Then, for every $k \geq 0$, $\vhat_k$ is absolutely continuous on $[0,c]$.
Furthermore, for $k \leq K$, $\vhat_k$ is continuously differentiable
on $(0,c)$ and $\vhat_k' \equiv 0$.
For $k > K$, there is a finite set $\cD_k \subseteq (0,c)$ such that
$\vhat_k$ is continuously differentiable on $(0,c) \setminus \cD_k$,
with derivative given by the recursion
\eqref{eq:vhat-derivative-recursion}.
\end{lemma}

We now collect the key constants that enter the derivative bound.
With the capacity $c$ fixed, there are two main inputs.
The first is the weight distribution $F$, with the constants $M_F$
from the rate condition \eqref{eq:typical-class-limit-with-rate}
and $C_F$ from Corollary~\ref{cor:karamata-ratios}.
Together with $c$, the distribution also determines the constant
$C_\epsilon$ in
Proposition~\ref{prop:joint-properties-f-epsilon}--\ref{property:local-reciprocal-increments}.
The second input is the parameter $\Theta > 1$, which defines the
small-capacity boundary $\xhat_k$ and the interval $(0, \xhat_k)$.
The same parameter $\Theta$ determines $\thetalb$ and $\thetaub$, which
bound the ratio $\epsilon_k(x)/x$ from below for $x < \xhat_k$
and from above for $x \geq \xhat_k$, respectively.
Recall the constant
$\lambda_\Theta = \thetaub^{1/(C_F - 1)}$
from Lemma~\ref{lemma:density-lower-bound-large-capacity}.
When $\Theta$ is sufficiently large so that $\thetaub < 1/2$,
the constant
\(
    C_\Theta
    =
    \frac{ \thetaub + 2 \lambda_\Theta }{ 1 - \thetaub }
\)
from Lemma~\ref{lemma:absolute-coefficient-estimate}
controls the amplification of the propagated error in
Proposition~\ref{prop:Pk-absolute-upper-bound}.
Finally, $F$ and $\Theta$ together determine
$M_v = C_F \Theta$ in
Lemma~\ref{lemma:vhat-upperbound-small-capacity},
and the constant $C_v = (1 + M_v)e^{M_v}$ in
Proposition~\ref{prop:vhat-intermediate-derivative-bounds-small-capacity}.
As $\Theta \to \infty$, both $\thetaub$ and $C_\Theta$ converge to zero.
In addition to the preceding requirements, we choose $\Theta$ sufficiently large so that 
\[
    d_\Theta
    \equiv
    \frac{ 1 }{ 8 C_F } - C_\Theta
    >
    0.
\]
This positive gap allows a finite choice of $\Delta$ in the
proposition below.

\begin{proposition}[Non-asymptotic derivative bounds]
\label{prop:vhat-derivative-bounds}
Suppose that $F$ belongs to the typical class of
Definition~\ref{def:typical-class} with parameter $r$.
Choose $\Theta$ sufficiently large that $\thetaub < 1/2$ and 
\(
    d_\Theta
    =
    \frac{ 1 }{ 8 C_F } - C_\Theta
    >
    0.
\)
Let $\Delta$ be any constant such that 
\[
    \Delta
    \geq
    \max\left\{
        \Theta C_F C_v,
        (2\Theta)^{3/4},
        \left( \frac{ c }{ \epsilon_{K+1}(c) } \right)^{3/4},
        \frac{ 3 C_\epsilon + 1 }{ 2 d_\Theta }
    \right\}.
\]
Then, for every $k \geq K$ and every $x \in (0,c) \setminus \cD_k$, 
one has 
\[
    \abs{ \vhat_k'(x) - \frac{ 1 }{ \epsilon_{k+1}(x) } }
    \leq
    \frac{ \Delta }{ x^{3/4} \epsilon_{k+1}(x)^{1/4} }.
\]
For every $x \in (0, c]$ and $w \in [0, \hhat_{k+1}(x)]$,
one also has
\begin{equation}
\label{eq:intermediate-integrand-overall-bound}
    \abs{
        \vhat_k(x)
        - \vhat_k(x-w)
        - \frac{ w }{ \epsilon_{k+1}(x) }
    }
    \leq
    \int_{x-w}^x
        \frac{ \Delta }{ y^{3/4} \epsilon_{k+1}(y)^{1/4} }
    \, dy
    +
    \frac{ w^2 }{ (x-w) \epsilon_{k+1}(x) }.
\end{equation}
\end{proposition}

Proposition~\ref{prop:vhat-derivative-bounds} connects the marginal value
of capacity, $\vhat_k'(x)$, with the reciprocal consumption function
$1 / \epsilon_{k+1}(x)$. This relationship has a natural economic
interpretation. 
The consumption function $\epsilon_{k+1}(x)$ is the optimal acceptance threshold in the residual fluid knapsack problem with $k+1$ remaining periods and capacity $x$ (Section~\ref{sec:DSKP&ReoptimizedHeuristic}).
The reciprocal of the consumption function, $1 / \epsilon_{k+1}(x)$, is the optimal dual multiplier associated with the capacity constraint in the fluid relaxation. It therefore represents the shadow price of capacity, measuring the marginal value of an
additional unit of capacity in the fluid problem.

The proof uses the auxiliary boundary $\xhat_k$ to divide the state space.
In Section~\ref{subsec:vhat-derivative-bounds-small-capacity}, we consider \emph{small} capacity
$x \in (0, \xhat_k)$, an interval which encompasses greedy as well as conservative selections.
We first establish an inductive bound that implies that $\abs{ \vhat_k'(x) } \leq C_v k f(x)$ for some $C_v > 0$, 
and then convert it into the form of the bound asserted.
In Section~\ref{subsec:vhat-derivative-bounds-large-capacity}, we consider \emph{large} capacity $x \geq \xhat_k$, a region where the reoptimized heuristic acts conservatively. 
In this region, we decompose the approximation error using the derivative recursion into two terms:
the propagated error $P_k(x)$, which carries over the error inherited from $\vhat_{k-1}'$,
and the local error $R_k(x)$, which does not involve the value functions.

Proposition~\ref{prop:vhat-derivative-bounds} has two important consequences.
The first is that, over the selection interval
$w \in [0, \hhat_{k+1}(x)]$, the difference
$\vhat_k(x) - \vhat_k(x-w)$ is bounded in absolute value by a constant
independent of $k$, $x$, and $w$.
This parallels a basic property of the optimal dynamic programming policy:
on its acceptance interval, the corresponding decrement in the optimal
value function lies between zero and one.

\begin{corollary}[Bounded difference in value functions within selection interval]
\label{cor:uniform-bound-on-difference-in-value-function}
Under the assumptions of Proposition~\ref{prop:vhat-derivative-bounds},
there exists a finite constant $\Delta_v$ such that, for every $k \geq 1$ and $x \in [0,c]$,
\[
    \abs{ \vhat_k(x) - \vhat_k(x-w) }
    \leq \Delta_v
    \qquad
    \text{for all } w \in [0, \hhat_{k+1}(x)].
\]
\end{corollary}

The second consequence of Proposition~\ref{prop:vhat-derivative-bounds} is more technical and will be useful in the
variance analysis of Section~\ref{sec:variance-asymptotics}.
It provides a quantitative bound on the error obtained when replacing
$\left( 1 + \vhat_k(x-w) - \vhat_k(x) \right)^2$
by its squared linear approximation
$\left( 1 - \frac{ w }{ \hhat_{k+1}(x) } \right)^2$
and integrating against the weight distribution over the selection interval.

\begin{corollary}[Squared linear approximation error]
\label{cor:error-bound-of-linear-integrand-approximation-k}
Under the assumptions of Proposition~\ref{prop:vhat-derivative-bounds},
there is a constant $\Delta_2 < \infty$ such that, for all $k\geq K$ and $x\in(0,c]$,
\begin{equation}
\label{eq:error-bound-of-linear-integrand-approximation-square-k}
\abs{
    \int_0^{\hhat_{k+1}(x)}
    \left\{
    \left(1+\vhat_k(x-w)-\vhat_k(x)\right)^2
    -\left(1-\frac{w}{\hhat_{k+1}(x)}\right)^2
    \right\} f(w)\,dw
    }
    \leq
    \Delta_2 
    \left\{
    \frac{F^{1/4}(\hhat_{k+1}(x))}{(k+1)^{3/4}}
    +
    \frac{1}{k+1}
    \right\}.
    \end{equation}
\end{corollary}

\subsection{Derivative bound with small capacity}
\label{subsec:vhat-derivative-bounds-small-capacity}

When the remaining capacity $x$ is small, specifically below $\xhat_k$,
the reoptimized heuristic operates in the greedy regime ($x \leq x_k$) as well as in the
conservative regime ($x_k < x < \xhat_k$).
Since every selected item has weight at most $x$, the expected reward is at most $k F(x)$,
the expected number of arrivals that fit, and in the greedy regime this bound describes the
value function $\vhat_k$ well.
Informally, the derivative $k f(x)$ of this upper bound then also describes $\vhat_k'(x)$.
Proposition~\ref{prop:vhat-intermediate-derivative-bounds-small-capacity} makes this
precise and shows that the same scale persists on all of $(0, \xhat_k)$, that is,
$\abs{ \vhat_k'(x) } \leq C_v k f(x)$.
In other words, when capacity is small, the derivative of the value function in the conservative regime behaves very much like the same derivative in the greedy regime.
This is partly because the derivative recursion
\eqref{eq:vhat-derivative-recursion} samples $\vhat_{k-1}'$ over
$[x - \epsilon_k(x), x]$, which reaches deep into the greedy regime since $\epsilon_k(x)$ is
comparable to $x$ on this interval.

\begin{proposition}[Intermediate derivative bound with small capacity]
\label{prop:vhat-intermediate-derivative-bounds-small-capacity}
Suppose that the weight distribution $F$ belongs to the typical class of
Definition~\ref{def:typical-class}.
% Let $M_v$ be the constant in
% Lemma~\ref{lemma:vhat-upperbound-small-capacity}, which depends
% only on $F$ and $\Theta$.
Then, for every $k \geq K$ and every $x \in (0, \xhat_k) \setminus \cD_k$,
\begin{equation}
\label{eq:vhat-derivative-bounds-small-capacity-induction}
    \abs{ \vhat_k'(x) }
    \leq
    (1 + M_v) \{ 1 + F(x) \}^{k-1} k f(x).
\end{equation}
Consequently, if $C_v = (1 + M_v) e^{M_v}$, we also have for the same $k$ and $x$ that 
\begin{equation}
\label{eq:vhat-derivative-bounds-small-capacity}
    \abs{ \vhat_k'(x) } \leq C_v k f(x).
\end{equation}
\end{proposition}

With a bit more work, the bound \eqref{eq:vhat-derivative-bounds-small-capacity} can be expressed in terms of the reciprocal of the consumption function, just like Proposition~\ref{prop:vhat-derivative-bounds} requires. 
This quantification is settled by the next result. 

\begin{proposition}[Non-asymptotic derivative bound with small capacity]
\label{prop:vhat-derivative-bounds-small-capacity}
Suppose that the weight distribution $F$ belongs to the typical class of
Definition~\ref{def:typical-class}.
Then, with
\(
    \Delta
    \geq \max\left\{
        \Theta C_F C_v,
        (2\Theta)^{3/4}
    \right\},
\)
one has
\[
    \abs{ \vhat_k'(x) - \frac{ 1 }{ \epsilon_{k+1}(x) } }
    \leq
    \frac{ \Delta }{ x^{3/4} \epsilon_{k+1}(x)^{1/4} }
\]
for every $k \geq K$ and every $x \in (0, \xhat_k) \setminus \cD_k$.
\end{proposition}

\subsection{Derivative bound with large capacity}
\label{subsec:vhat-derivative-bounds-large-capacity}

When the remaining capacity $x$ is large, specifically at least $\xhat_k$,
the reoptimized heuristic operates only in the conservative regime, with threshold
$\hhat_k(x) = \epsilon_k(x)$.
Here, every accepted item consumes only a small fraction of the remaining capacity, 
so the fluid relaxation describes the heuristic well, and the marginal value of capacity is
close to the shadow price $1 / \epsilon_{k+1}(x)$.

To make this precise, fix $k > K$ and $x \in [\xhat_k, c) \setminus \cD_k$.
By absolute continuity of $\vhat_{k-1}$,
\[
    \vhat_{k-1}(x) - \vhat_{k-1}(x - \epsilon_k(x))
    =
    \int_0^{\epsilon_k(x)} \vhat_{k-1}'(x - w) \, dw.
\]
Substituting this identity into the conservative branch of the derivative recursion
\eqref{eq:vhat-derivative-recursion} and subtracting $1 / \epsilon_{k+1}(x)$,
we decompose the approximation error into the sum of propagated and local errors as
\begin{equation}
\label{eq:vhat-derivative-decomposition-large-capacity}
    \vhat_k'(x) - \frac{1}{\epsilon_{k+1}(x)}
    =
    P_k(x) + R_k(x).
\end{equation}
In this decomposition, $P_k(x)$ is the \emph{propagated error},
which carries over the error inherited from $\vhat_{k-1}'$,
\begin{align}
    \label{eq:Pk}
    P_k(x)
    = 
    & \{ 1 - F(\epsilon_k(x)) \}
    \left(
        \vhat_{k-1}'(x) - \frac{1}{\epsilon_k(x)}
    \right) \\
    & +
    \int_0^{\epsilon_k(x)}
    \left(
        \vhat_{k-1}'(x - w)
        - \frac{1}{\epsilon_k(x - w)}
    \right)
    \left(
        f(w) - \frac{1}{k \epsilon_k(x)}
    \right)
    \, dw,
    \notag
\end{align}
and $R_k(x)$ is the \emph{local error}, the error that the recursion would produce in one step if
$\vhat_{k-1}'$ were exactly $1 / \epsilon_k$, and which does not involve the value functions:
\begin{align}
    \label{eq:Rk}
    R_k(x)
    = 
    & \frac{1}{\epsilon_k(x)}
    - \frac{1}{\epsilon_{k+1}(x)} \\
    & +
    \int_0^{\epsilon_k(x)}
    \left(
        \frac{1}{\epsilon_k(x - w)}
        - \frac{1}{\epsilon_k(x)}
    \right)
    f(w) \, dw 
    -
    \frac{1}{k \epsilon_k(x)}
    \int_0^{\epsilon_k(x)}
    \left(
        \frac{1}{\epsilon_k(x - w)}
        - \frac{1}{\epsilon_k(x)}
    \right)
    \, dw.
    \notag
\end{align}

To bound the absolute value of the approximation error
\eqref{eq:vhat-derivative-decomposition-large-capacity}, the rest of this section studies
$\abs{ P_k(x) }$ and $\abs{ R_k(x) }$ separately.
The bound for $\abs{ P_k(x) }$ is given in Proposition~\ref{prop:Pk-absolute-upper-bound}
under a suitable assumption on the approximation error of $\vhat_{k-1}'$
while the bound for $\abs{ R_k(x) }$ follows in Proposition~\ref{prop:Rk-absolute-upper-bound}.

\begin{proposition}[Propagated-error bound]
\label{prop:Pk-absolute-upper-bound}
Suppose that $F$ belongs to the typical class of
Definition~\ref{def:typical-class} with parameter $r$.
Fix $\Theta$ sufficiently large that $\thetaub < 1/2$.
Let $k > K$ and let $M > 0$ be any constant such that
\begin{equation}
\label{eq:vhat-derivative-approximation-error-assumption}
    \abs{ \vhat_{k-1}'(x) - \frac{ 1 }{ \epsilon_k(x) } }
    \leq
    \frac{ M }{ x^{3/4} \epsilon_k(x)^{1/4} }
\end{equation}
for every $x \in (0, c) \setminus \cD_{k-1}$.
Then, for every $x \in (\xhat_k, c) \setminus \cD_k$,
\[
    \abs{ P_k(x) }
    \leq
    \left( 1 + \frac{ C_\Theta }{ k } \right)
    \frac{ M }{ x^{3/4} \epsilon_k(x)^{1/4} }.
\]
\end{proposition}

\begin{proposition}[Local-error bound]
\label{prop:Rk-absolute-upper-bound}
Suppose that $F$ belongs to the typical class of
Definition~\ref{def:typical-class}.
Fix $\Theta$ sufficiently large that $\thetaub < 1/2$.
% Let $C_\epsilon$ be the constant in
% Proposition~\ref{prop:joint-properties-f-epsilon}--\ref{property:local-reciprocal-increments}.
Then, for every $k \geq K$ and $x \in (\xhat_k, c)$,
\[
    \abs{ R_k(x) }
    \leq
    \left(\frac{ 3 C_\epsilon + 1 }{ 2 }\right)
    \frac{ 1 }{ k x^{3/4} \epsilon_k(x)^{1/4} }.
\]
\end{proposition}

The last two estimates hold when the remaining capacity is large, complementing the small-capacity bound in Proposition~\ref{prop:vhat-derivative-bounds-small-capacity}.
Together, these bounds will be used in the induction argument that proves Proposition~\ref{prop:vhat-derivative-bounds}.

\section{Variance asymptotics}
\label{sec:variance-asymptotics}

This section establishes the variance asymptotics for the number
$\widehat N_T(c)$ of selections made by the reoptimized heuristic
$\widehat \pi$, with $T$ periods and fixed initial knapsack capacity $c$.
The proofs of the results in this section are given in
Appendix~\ref{apx:variance-asymptotics-proofs}.

\begin{proposition}[Variance asymptotics]
\label{prop:variance-asymptotics}
Suppose that the weight distribution $F$ belongs to the typical class
of Definition~\ref{def:typical-class} with parameter $r$.
Set
\(
    \sigma^2 = \frac{ 2 }{ (r+1)(r+2) }.
\)
Then
\(
    \Var{ \widehat N_T(c) } \sim \sigma^2 \E{ \widehat N_T(c) },
\)
as $T \to \infty$.
\end{proposition}

The proof of Proposition~\ref{prop:variance-asymptotics} uses a classical martingale
argument.
Let $\{\cF_t\}_{t=0}^T$, with
$\cF_t = \sigma(W_1, \ldots, W_t)$, be the natural filtration of the
item weights, and consider the Doob martingale associated with
$\widehat N_T(c)$:
\[
    M_t
    =
    \E{ \widehat N_T(c) \midbar \cF_t }
    =
    \sum_{j=1}^t
        \1{ W_j \leq \hhat_{T-j+1}(X_{j-1}) }
    + \vhat_{T-t}(X_t),
    \qquad \text{for $t = 0, \ldots, T$.}
\]

For $t \in [T]$, let
\(
    G_t = \{ W_t \leq \hhat_{T-t+1}(X_{t-1}) \}
\)
denote the event that policy $\widehat \pi$ selects the $t$th arriving
item. 
By the remaining capacity recursion
\eqref{eq:remaining-capacity-process}, the martingale difference
\[
d_t 
= M_t - M_{t-1}
= A_t + B_t,
\]
where
\begin{equation}
\label{eq:At}
    A_t
    =
    \left\{
        1 + \vhat_{T-t}(X_{t-1} - W_t)
        - \vhat_{T-t}(X_{t-1})
    \right\}
    \1{ G_t }
\end{equation}
and
\begin{equation*}
%\label{eq:Bt}
    B_t
    =
    \vhat_{T-t}(X_{t-1}) - \vhat_{T-t+1}(X_{t-1}).
\end{equation*}
Here, $A_t$ is the \emph{selection difference}: on $G_t$ and with $T-t$
periods remaining, it compares the unit reward plus the continuation value
after selection with the continuation value after rejection.
In \eqref{eq:At}, we evaluate the expression in braces only on $G_t$ and set $A_t = 0$ otherwise.
The term $B_t$ is the \emph{time difference}, representing the change in
the value function when the number of remaining periods decreases by one
at fixed capacity.

Since $B_t$ is $\cF_{t-1}$-measurable and
$\E{ d_t \midbar \cF_{t-1} } = 0$, we have
$B_t = -\E{ A_t \midbar \cF_{t-1} }$.
Hence,
\begin{equation}
\label{eq:dt-squared-decomposition}
    \E{ d_t^2 \midbar \cF_{t-1} }
    =
    \E{ A_t^2 \midbar \cF_{t-1} }
    + B_t^2
    + 2 \E{ A_t \midbar \cF_{t-1} } B_t
    =
    \E{ A_t^2 \midbar \cF_{t-1} } - B_t^2.
\end{equation}

The independence of future item weights and the structure of the reoptimized heuristic 
imply that the conditional variance of the remaining martingale
differences depends on the observed history only through the remaining
capacity and horizon.
For $t = 0, \ldots, T$, we therefore write
\begin{equation*}
    \omega_{T-t}(X_t)
    =
    \Var{ \sum_{j=t+1}^T d_j \midbar \cF_t },
\end{equation*}
so setting $t = 0$ gives
$\omega_T(c) = \Var{ \widehat N_T(c) }$.
Conditional orthogonality and the tower property give
\begin{equation}
\label{eq:conditional-variance-decomposition}
    \omega_{T-t}(X_t)
    =
    \sum_{j=t+1}^T \E{ d_j^2 \midbar \cF_t }
    =
    \sum_{j=t+1}^T \E{ A_j^2 \midbar \cF_t }
    - \sum_{j=t+1}^T \E{ B_j^2 \midbar \cF_t }.
\end{equation}

The remainder of the proof of
Proposition~\ref{prop:variance-asymptotics} is organized as follows.
Section~\ref{sec:At-squared-contribution} controls the approximation error
in the conditional second moments of $A_t$.
Section~\ref{sec:Bt-square-contribution} bounds the accumulated contribution
of $B_t^2$.

\subsection{The contribution of $A_t^2$}
\label{sec:At-squared-contribution}

Here we study the conditional second moments
$\E{ A_t^2 \midbar \cF_{t-1} }$, for $t \in [T]$, and quantify
their contribution to the first sum in the conditional variance
decomposition \eqref{eq:conditional-variance-decomposition}.
By the definition of $A_t$ and the independence of $W_t$ from $\cF_{t-1}$,
we have the integral representation
\begin{equation}
\label{eq:At-conditional-second-moment}
    \E{ A_t^2 \midbar \cF_{t-1} }
    =
    \int_0^{\hhat_{T-t+1}(X_{t-1})}
        \left\{
            1 + \vhat_{T-t}(X_{t-1}-w)
            - \vhat_{T-t}(X_{t-1})
        \right\}^2
        f(w) \, dw.
\end{equation}
We analyze this integral by replacing the squared expression with
its squared linear approximation
$\left( 1 - \frac{ w }{ \hhat_{T-t+1}(X_{t-1}) } \right)^2$
and controlling the resulting error.
Lemma~\ref{lem:linear-approximation-second-moment-bounds}
identifies $\sigma^2 F(\hhat_{T-t+1}(X_{t-1}))$ as the leading term
of the squared linear integral over
$[0, \hhat_{T-t+1}(X_{t-1})]$.
Using this result, we show
that the accumulated error in approximating the conditional second
moments by $\sigma^2 F(\hhat_{T-t+1}(X_{t-1}))$ is asymptotically
negligible on the scale of $\vhat_T(c)$.

\begin{proposition}[Asymptotic negligibility of $A_t^2$'s accumulated approximation error]
\label{prop:A_t-square-accumulated-approximation-error-negligible}
Suppose that $F$ belongs to the typical class of Definition~\ref{def:typical-class}. Then, as $T\to \infty$,
\begin{equation*}
\frac{ 1 }{ \vhat_T(c) }
\sum_{t=1}^T
\E{
\abs{
\E{ A_t^2\midbar\cF_{t-1} }
-\sigma^2F(\hhat_{T-t+1}(X_{t-1}))
}
}
\longrightarrow 0.
\end{equation*}
\end{proposition}

As an intermediate step to the proof of Proposition~\ref{prop:A_t-square-accumulated-approximation-error-negligible}, 
we now bound weighted sums of future acceptance probabilities in terms of
$\vhat_{T-t}(X_t)$, the expected number of remaining selections given $\cF_t$.

\begin{lemma}[Expected weighted sum of future acceptance probabilities]
\label{lem:series-upper-bound}
Let the remaining capacity process $\{X_t: t=0,1,\ldots,T\}$ under the adaptive heuristic $\widehat \pi$ be defined as in \eqref{eq:remaining-capacity-process}. For any fixed $a\in (0,1]$ and for all $t \in \{0,\ldots,T-1\}$,
\[
\sum_{j=t+1}^T
\E{
\frac{ F^a(\hhat_{T-j+1}(X_{j-1})) }
     { (T-j+1)^{1-a} }
\midbar \cF_t
}
\leq
\vhat_{T-t}(X_t)^{a}
\{1+\log(T-t)\}^{1-a}.
\]
\end{lemma}

\subsection{The contribution of $B_t^2$}
\label{sec:Bt-square-contribution}

Since $B_t = -\E{ A_t \midbar \cF_{t-1} }$, the definition
\eqref{eq:At} gives the integral representation
\begin{equation}
\label{eq:Bt-integral-representation}
    B_t
    =
    - \int_0^{\hhat_{T-t+1}(X_{t-1})}
        \left\{
            1 + \vhat_{T-t}(X_{t-1}-w)
            - \vhat_{T-t}(X_{t-1})
        \right\}
        f(w) \, dw.
\end{equation}
The bounded-difference estimate from
Corollary~\ref{cor:uniform-bound-on-difference-in-value-function}
bounds the integrand in absolute value by $(1 + \Delta_v) f(w)$,
so that
$\abs{ B_t } \leq (1 + \Delta_v) F(\hhat_{T-t+1}(X_{t-1}))$.
It therefore suffices to show that the sum of the expected squared
acceptance probabilities, $\sum_{t=1}^T \E{ F(\hhat_{T-t+1}(X_{t-1}))^2 }$, is negligible on the scale of $\vhat_T(c)$
to obtain the same conclusion for $\sum_{t=1}^T \E{ B_t^2 }$.
The following proposition establishes both conclusions.

\begin{proposition}[Negligible sums of squared acceptance probabilities and time differences]
\label{prop:B_t-square-total-contribution-negligible}
Suppose that $F$ belongs to the typical class of
Definition~\ref{def:typical-class}, and let the remaining capacity process
$\{X_t: t = 0, 1, \ldots, T\}$ under the reoptimized heuristic
$\widehat \pi$ be defined as in \eqref{eq:remaining-capacity-process}.
Then
$$
    \lim_{T \to \infty}
    \frac{ 1 }{ \vhat_T(c) }
    \sum_{t=1}^T
        \E{ F(\hhat_{T-t+1}(X_{t-1}))^2 }
    = 0
    \qquad \text{and} \qquad
    \lim_{T \to \infty}
    \frac{ 1 }{ \vhat_T(c) }
    \sum_{t=1}^T \E{ B_t^2 }
    = 0.
$$
\end{proposition}

\section{Distributional equivalence of the reoptimized heuristic: two central limit theorems}
\label{sec:clt}

In this section, we establish the central limit theorem for the number
of items selected by the reoptimized heuristic and then transfer the
limit to the optimal policy.
Our argument uses the following version of the martingale central limit
theorem \citep[cf.][]{McLeish:AP1974}.
The proofs of the results in this section and of Theorem~\ref{thm:main-theorem} are given in Appendix~\ref{apx:clt-proofs}.

\begin{proposition}[Martingale central limit theorem]
\label{prop:martingale-clt}
For each $T \geq 1$, let $\{Z_{T,t}: 1 \leq t \leq T\}$ be a martingale
difference sequence with respect to the filtration
$\{\cF_{T,t}: 0 \leq t \leq T\}$.
If
\begin{equation}
\label{eq:negligibility-condition-of-martingale-clt}
    \max_{1 \leq t \leq T} \norm{ Z_{T,t} }_\infty
    \longrightarrow 0
    \qquad \mbox{as } T \to \infty,
\end{equation}
and
\begin{equation}
\label{eq:convergence-in-probability-condition-of-martingale-clt}
    \sum_{t=1}^T \E{ Z_{T,t}^2 \midbar \cF_{T,t-1} }
    \xrightarrow{p} 1
    \qquad \mbox{as } T \to \infty,
\end{equation}
then we have the convergence in distribution
\[
    \sum_{t=1}^T Z_{T,t}
    \xrightarrow{d} N(0,1)
    \qquad \mbox{as } T \to \infty.
\]
\end{proposition}

To apply Proposition~\ref{prop:martingale-clt} to our problem, recall
the martingale differences $d_t = M_t - M_{t-1}$ and the variance
$\Var{ \widehat N_T(c) } = \omega_T(c)$ from
Section~\ref{sec:variance-asymptotics}.
For all sufficiently large $T$, we scale each $d_t$ by
$\sqrt{ \omega_T(c) }$ and set
\begin{equation}
\label{eq:Zt}
    Z_t := \frac{ d_t }{ \sqrt{ \omega_T(c) } }
    \qquad \text{for } t = 1, \ldots, T.
\end{equation}
Then $\{Z_t: 1 \leq t \leq T\}$ is a martingale difference sequence
with respect to $\{\cF_t: 0 \leq t \leq T\}$, with the dependence on $T$
implicit in the notation.
With this specialization, we will show that the sum
\[
    \sum_{t=1}^T Z_t
    =
    \frac{ \widehat N_T(c) - \vhat_T(c) }{ \sqrt{ \omega_T(c) } }
\]
converges in distribution to a standard normal random variable.
As an intermediate step,
Section~\ref{subsec:quadratic-variation-concentration} studies the
predictable quadratic variation
\begin{equation}\label{eq:QT}
    Q_T := \sum_{t=1}^T \E{ d_t^2 \midbar \cF_{t-1} }
\end{equation}
of the martingale $\{M_t: 0 \leq t \leq T\}$ and proves that
\[
 \sum_{t=1}^T \E{ Z_t^2 \midbar \cF_{t-1} } = \frac{ Q_T }{ \omega_T(c) } \xrightarrow{p} 1,
\]
as required by
\eqref{eq:convergence-in-probability-condition-of-martingale-clt}.

\subsection{Convergence in probability of the predictable quadratic variation}
\label{subsec:quadratic-variation-concentration}

Section~\ref{sec:variance-asymptotics} identifies the asymptotic
behavior of $\E{ Q_T } = \omega_T(c)$.
To verify \eqref{eq:convergence-in-probability-condition-of-martingale-clt},
we must also control the fluctuations of the predictable quadratic variation $Q_T$ around its mean.
The following proposition establishes the required convergence
in probability by proving the stronger convergence of
$Q_T/\omega_T(c)$ in $L^1$.

\begin{proposition}[$L^1$ convergence of the normalized predictable quadratic variation]
\label{prop:quadratic-variation-convergence}
Suppose that $F$ belongs to the typical class of
Definition~\ref{def:typical-class}. Then
\begin{equation*}
    \E{ \abs{ \frac{ Q_T }{ \omega_T(c) } - 1 } }
    \longrightarrow 0
    \qquad \mbox{as } T \to \infty.
\end{equation*}
In particular, $Q_T/\omega_T(c)$ converges to one in probability.
\end{proposition}

As an intermediate step towards proving
Proposition~\ref{prop:quadratic-variation-convergence}, we study the
sum of the conditional selection probabilities.
Recall that $\widehat N_T(c) = \sum_{t=1}^T \1{ W_t \leq \hhat_{T-t+1}(X_{t-1}) }$, and replace each summand by its conditional expectation given the past. 
This yields the selection compensator
\begin{equation}
\label{eq:LambdaT}
    \Lambda_T
    =
    \sum_{t=1}^T F(\hhat_{T-t+1}(X_{t-1})).
\end{equation}
Each summand in $\Lambda_T$ is the probability of selection given
the remaining capacity before item $t$ arrives.
By the tower property,
$\E{ \Lambda_T } = \E{ \widehat N_T(c) } = \vhat_T(c)$.
The following lemma compares $\Lambda_T$ with the realized selection
count $\widehat N_T(c)$ to show that $\Lambda_T/\vhat_T(c)$ converges
to one in mean square.

\begin{lemma}[Mean-square convergence of the normalized selection compensator]
\label{lemma:compensator-L2-convergence}
Suppose that the weight distribution $F$ belongs to the typical class of
Definition~\ref{def:typical-class}. Then
\begin{equation*}
    \E{ \left( \frac{ \Lambda_T }{ \vhat_T(c) } - 1 \right)^2 }
    \longrightarrow 0,
    \qquad \mbox{as } T \to \infty.
\end{equation*}
\end{lemma}

\section{Closing remarks}\label{sec:closing-remarks}

Different heuristics can achieve the same asymptotic expected reward
while exhibiting different variances and different limiting distributions.
In this paper, we study the dynamic and stochastic knapsack problem
with unit rewards and show that, after the same centering and scaling,
the number of selections made by a reoptimized heuristic has the same
normal limiting distribution as the number of selections made by the
optimal dynamic programming policy.
The expected-performance guarantee of the reoptimized heuristic is
therefore accompanied by a stronger distributional equivalence.

Our analysis highlights a connection between value functions and
the stochastic fluctuations of realized rewards.
While the value function describes the expected number of selections,
estimates of its derivative allow us to determine the corresponding
variance and limiting distribution.
This connection suggests that value-function sensitivity can also be
useful for studying the higher-order behavior of other online policies.

From a decision-making perspective, we demonstrate that expected
rewards alone may provide limited guidance when choosing among
asymptotically optimal heuristics.
Understanding higher moments, tail probabilities, and limiting
distributions provides a fuller basis for comparing their realized
performance and assessing risk.

\ACKNOWLEDGMENT{%
This material is based upon work supported by the National Science Foundation under CAREER Award No. 1553274.
}

{
    \SingleSpacedXI
    \bibliography{KnapsackCLT}

@book{BinghamGoldieTeugels:CUP1989,
	author = {Bingham, Nicholas H and Goldie, Charles M and Teugels, Jef L},
	publisher = {Cambridge University Press},
	title = {Regular Variation},
	volume = {27},
	year = {1989}}

@article{ChenChenWangZhang-WP-2026,
	author = {Chen, Chen and Chen, Yilun and Wang, Wenjia and Zhang, Jiawei},
	journal = {SSRN},
	pages = {},
	title = {Online Resource Allocation: Logarithmic Corrections and Constant Optimality Gaps},
	url = {http://dx.doi.org/10.2139/ssrn.7515319},
	volume = {},
	year = {2026},
    }

@article{Lamperti:TAMS1958,
	author = {Lamperti, John},
	doi = {10.2307/1993222},
	fjournal = {Transactions of the American Mathematical Society},
	issn = {0002-9947},
	journal = {Trans. Amer. Math. Soc.},
	mrclass = {60.00},
	mrnumber = {94863},
	mrreviewer = {J. Feldman},
	pages = {380--387},
	title = {An occupation time theorem for a class of stochastic processes},
	url = {https://doi.org/10.2307/1993222},
	volume = {88},
	year = {1958}}

@article{Arlotto-Xie:SSY2020,
	author = {Arlotto, Alessandro and Xie, Xinchang},
	doi = {10.1287/stsy.2019.0055},
	fjournal = {Stochastic Systems},
	journal = {Stoch. Syst.},
	mrclass = {90C27 (60K30 62L99 90C15)},
	mrnumber = {4119249},
	number = {2},
	pages = {170--191},
	title = {Logarithmic regret in the dynamic and stochastic knapsack problem with equal rewards},
	url = {https://doi.org/10.1287/stsy.2019.0055},
	volume = {10},
	year = {2020}}

@inproceedings{BhalgatGoelKhanna:SDA2011,
	author = {Bhalgat, Anand and Goel, Ashish and Khanna, Sanjeev},
	booktitle = {Proceedings of the {T}wenty-{S}econd {A}nnual {ACM}-{SIAM} {S}ymposium on {D}iscrete {A}lgorithms},
	mrclass = {90C15 (68W25 90C27)},
	mrnumber = {2858426},
	pages = {1647--1665},
	publisher = {SIAM, Philadelphia, PA},
	title = {Improved approximation results for stochastic knapsack problems},
	year = {2011}}

@article{BladoToriello:SIAM2019,
	author = {Blado, Daniel and Toriello, Alejandro},
	doi = {10.1137/16M1101209},
	fjournal = {SIAM Journal on Optimization},
	issn = {1052-6234},
	journal = {SIAM J. Optim.},
	mrclass = {90C34 (90C40)},
	mrnumber = {3895343},
	number = {1},
	pages = {1--30},
	title = {Relaxation analysis for the dynamic knapsack problem with stochastic item sizes},
	url = {https://doi.org/10.1137/16M1101209},
	volume = {29},
	year = {2019}}

@article{BladoHuToriello:SIAM2016,
	author = {Blado, Daniel and Hu, Weihong and Toriello, Alejandro},
	doi = {10.1137/15M1036233},
	fjournal = {SIAM Journal on Optimization},
	issn = {1052-6234},
	journal = {SIAM J. Optim.},
	mrclass = {90C39 (90C10 90C27 90C34)},
	mrnumber = {3537890},
	mrreviewer = {Charles E. Blair},
	number = {3},
	pages = {1625--1648},
	title = {Semi-infinite relaxations for the dynamic knapsack problem with stochastic item sizes},
	url = {https://doi.org/10.1137/15M1036233},
	volume = {26},
	year = {2016}}

@article{DeanGoemansVondrak:MOR2008,
	author = {Dean, Brian C. and Goemans, Michel X. and Vondr\'{a}k, Jan},
	doi = {10.1287/moor.1080.0330},
	fjournal = {Mathematics of Operations Research},
	issn = {0364-765X},
	journal = {Math. Oper. Res.},
	mrclass = {90C27 (68W25 90C15 90C59)},
	mrnumber = {2464653},
	number = {4},
	pages = {945--964},
	title = {Approximating the stochastic knapsack problem: the benefit of adaptivity},
	url = {https://doi.org/10.1287/moor.1080.0330},
	volume = {33},
	year = {2008}}

@inproceedings{DeanGoemansVondrak:FOCS2004,
	author = {Dean, Brian C. and Goemans, Michel X. and Vondr\'{a}k, Jan},
	booktitle = {Proceedings of the 45th {A}nnual {IEEE} {S}ymposium on {F}oundations of {C}omputer {S}cience},
	doi = {10.1109/FOCS.2004.15},
	pages = {208--217},
	publisher = {IEEE Press, Piscataway, NJ},
	title = {Approximating the stochastic knapsack problem: the benefit of adaptivity},
	year = {2004}}

@inproceedings{DeanGoemansVondrak:SDA2005,
	author = {Dean, Brian C. and Goemans, Michel X. and Vondr\'{a}k, Jan},
	booktitle = {Proceedings of the {S}ixteenth {A}nnual {ACM}-{SIAM} {S}ymposium on {D}iscrete {A}lgorithms},
	mrclass = {90C15 (68W25 90B60)},
	mrnumber = {2298288},
	pages = {395--404},
	publisher = {ACM, New York, NY},
	title = {Adaptivity and approximation for stochastic packing problems},
	year = {2005}}

@inproceedings{BaruahHaritsaSharma:IEEE1994,
	author = {Baruah, Sanjoy and Haritsa, Jayant and Sharma, Nitin},
	booktitle = {1994 Proceedings Real-Time Systems Symposium},
	pages = {228--236},
	publisher = {IEEE, San Juan, PR},
	title = {On-line scheduling to maximize task completions},
	year = {1994}}

@article{ArlottoGurvich:SSY2019,
	author = {Arlotto, Alessandro and Gurvich, Itai},
	fjournal = {Stochastic Systems},
	journal = {Stochastic Systems},
	number = {3},
	pages = {231-260},
	title = {Uniformly bounded regret in the multi-secretary problem},
	volume = {9},
	year = {2019}}

@article{Vera-Banerjee:MS2021,
author = {Vera, Alberto and Banerjee, Siddhartha},
title = {The Bayesian Prophet: A Low-Regret Framework for Online Decision Making},
journal = {Management Science},
volume = {67},
number = {3},
pages = {1368-1391},
year = {2021},
doi = {10.1287/mnsc.2020.3624},
}

@article{ArlottoGansSteele:OR2014,
	author = {Arlotto, Alessandro and Gans, Noah and Steele, J. Michael},
	doi = {10.1287/opre.2014.1281},
	fjournal = {Operations Research},
	issn = {0030-364X},
	journal = {Oper. Res.},
	mrclass = {90C40},
	mrnumber = {3246083},
	mrreviewer = {Enrique Lemus-Rodr\'{\i}guez},
	number = {4},
	pages = {864--875},
	title = {Markov decision problems where means bound variances},
	url = {https://doi.org/10.1287/opre.2014.1281},
	volume = {62},
	year = {2014}}

@article{PapastavrouRajagopalanKleywegt:MS1996,
	author = {Papastavrou, Jason D. and Rajagopalan, Srikanth and Kleywegt, Anton J.},
	fjournal = {Management Science},
	issn = {00251909},
	journal = {Management Sci.},
	number = {12},
	pages = {1706--1718},
	title = {The Dynamic and Stochastic Knapsack Problem with Deadlines},
	volume = {42},
	year = {1996}}

@article{DermanLiebermanRoss:MS1978,
	author = {Derman, C. and Lieberman, G. J. and Ross, S. M.},
	doi = {10.1287/mnsc.24.5.554},
	fjournal = {Management Science. Journal of the Institute of Management Science. Application and Theory Series},
	issn = {0025-1909},
	journal = {Management Sci.},
	mrclass = {90B25},
	mrnumber = {521664},
	number = {5},
	pages = {554--561},
	title = {A renewal decision problem},
	url = {https://doi.org/10.1287/mnsc.24.5.554},
	volume = {24},
	year = {1978}}

@article{BalseiroBrown:OR2019,
	author = {Balseiro, Santiago R. and Brown, David B.},
	fjournal = {Operations Research},
	journal = {Oper. Res.},
	number = {2},
	pages = {577--597},
	title = {Approximations to Stochastic Dynamic Programs via Information Relaxation Duality},
	volume = {67},
	year = {2019}}

@article{ArlottoNguyenSteele:SPA2015,
	author = {Arlotto, Alessandro and Nguyen, Vinh V. and Steele, J. Michael},
	doi = {10.1016/j.spa.2015.03.009},
	fjournal = {Stochastic Processes and their Applications},
	issn = {0304-4149},
	journal = {Stochastic Process. Appl.},
	mrclass = {60C05 (60F05 60G40 60G42 90C27 90C39 90C40)},
	mrnumber = {3357621},
	mrreviewer = {Aleksey B. Piunovski\u{\i}},
	number = {9},
	pages = {3596--3622},
	title = {Optimal online selection of a monotone subsequence: a central limit theorem},
	url = {https://doi.org/10.1016/j.spa.2015.03.009},
	volume = {125},
	year = {2015}}

@article{BrussRobertson:AAP1991,
	author = {Bruss, F. Thomas and Robertson, James B.},
	coden = {AAPBBD},
	doi = {10.2307/1427625},
	fjournal = {Advances in Applied Probability},
	issn = {0001-8678},
	journal = {Adv. in Appl. Probab.},
	mrclass = {60G40 (62L15)},
	mrnumber = {1122878 (92i:60083)},
	mrreviewer = {S. M. Samuels},
	number = {3},
	pages = {612--623},
	title = {``{W}ald's lemma'' for sums of order statistics of i.i.d.\ random variables},
	volume = {23},
	year = {1991}}

@article{RheeTalagrand:JAP1991,
	author = {Rhee, WanSoo and Talagrand, Michel},
	coden = {JPRBAM},
	fjournal = {Journal of Applied Probability},
	issn = {0021-9002},
	journal = {J. Appl. Probab.},
	mrclass = {60G40 (90A05)},
	mrnumber = {1133802 (93d:60079)},
	mrreviewer = {Yoshio Ohtsubo},
	number = {4},
	pages = {919--923},
	title = {A note on the selection of random variables under a sum constraint},
	volume = {28},
	year = {1991}}

@article{CoffmanFlattoWeber:AAP1987,
	author = {Coffman, Jr., E. G. and Flatto, L. and Weber, R. R.},
	coden = {AAPBBD},
	doi = {10.2307/1427427},
	fjournal = {Advances in Applied Probability},
	issn = {0001-8678},
	journal = {Adv. in Appl. Probab.},
	mrclass = {60G40 (90B05 90C40)},
	mrnumber = {MR889945 (88i:60083)},
	mrreviewer = {S. M. Samuels},
	number = {2},
	pages = {454--473},
	title = {Optimal selection of stochastic intervals under a sum constraint},
	volume = {19},
	year = {1987}}

@article{BoshuizenKertz:AAP1999,
	author = {Boshuizen, Frans A and Kertz, Robert P},
	fjournal = {Advances in Applied Probability},
	journal = {Adv. in Appl. Probab.},
	number = {1},
	pages = {178--198},
	publisher = {Applied Probability Trust},
	title = {Smallest-fit selection of random sizes under a sum constraint: weak convergence and moment comparisons},
	volume = {31},
	year = {1999}}

@article{Cayley:ET1875,
	author = {Cayley, Arthur},
	journal = {Educational Times},
	note = {See \emph{The {C}ollected {M}athematical {P}apers of {A}rthur {C}ayley}, 10, 587--588, (1986). Cambridge University Press, Cambridge},
	pages = {18--19},
	title = {Mathematical questions and their solutions},
	volume = {22},
	year = {1875}}

@article{Moser:SM1956,
	author = {Moser, Leo},
	journal = {Scripta Mathematica},
	pages = {289--292},
	title = {On a problem of {C}ayley},
	volume = {22},
	year = {1956}}

@inproceedings{Kleinberg:SDA2005,
	author = {Kleinberg, Robert},
	booktitle = {Proceedings of the {S}ixteenth {A}nnual {ACM}-{SIAM} {S}ymposium on {D}iscrete {A}lgorithms},
	mrclass = {91B26 (68W40)},
	mrnumber = {2298313},
	pages = {630--631},
	publisher = {ACM, New York, NY},
	title = {A multiple-choice secretary algorithm with applications to online auctions},
	year = {2005}}

@article{Ma:MOR2018,
	author = {Ma, Will},
	doi = {10.1287/moor.2017.0884},
	fjournal = {Mathematics of Operations Research},
	issn = {0364-765X},
	journal = {Math. Oper. Res.},
	mrclass = {90C15 (90C27 90C40)},
	mrnumber = {3846073},
	mrreviewer = {Alireza Ghaffari-Hadigheh},
	number = {3},
	pages = {789--812},
	title = {Improvements and generalizations of stochastic knapsack and {M}arkovian bandits approximation algorithms},
	url = {https://doi.org/10.1287/moor.2017.0884},
	volume = {43},
	year = {2018}}

@article{JasinKumar:OR2013,
	author = {Jasin, Stefanus and Kumar, Sunil},
	doi = {10.1287/opre.2013.1216},
	fjournal = {Operations Research},
	issn = {0030-364X},
	journal = {Oper. Res.},
	mrclass = {91B26 (90C15 90C90)},
	mrnumber = {3158248},
	number = {6},
	pages = {1312--1320},
	title = {Analysis of deterministic {LP}-based booking limit and bid price controls for revenue management},
	url = {https://doi.org/10.1287/opre.2013.1216},
	volume = {61},
	year = {2013}}

@article{ReimanWang:MOR2008,
	author = {Reiman, Martin I. and Wang, Qiong},
	doi = {10.1287/moor.1070.0288},
	fjournal = {Mathematics of Operations Research},
	issn = {0364-765X},
	journal = {Math. Oper. Res.},
	mrclass = {90B15 (60F05 60F17 60K30 62P20 90B05)},
	mrnumber = {2415991},
	number = {2},
	pages = {257--282},
	title = {An asymptotically optimal policy for a quantity-based network revenue management problem},
	url = {https://doi.org/10.1287/moor.1070.0288},
	volume = {33},
	year = {2008}}

@article{Jiang-Zhang-arXiv2020-onlineresourceallocation,
      title={Online Resource Allocation with Stochastic Resource Consumption}, 
      author={Jiashuo Jiang and Jiawei Zhang},
      year={2020},
      journal={arXiv e-print 2012.07933},
      eprint={2012.07933},
      archivePrefix={arXiv},
      primaryClass={math.OC},
      url={https://arxiv.org/abs/2012.07933}, 
}

@article{GallegoVanRyzin:MS1994,
	author = {Gallego, Guillermo and van Ryzin, Garrett},
	doi = {10.1287/mnsc.40.8.999},
	fjournal = {Management Science},
	journal = {Management Sci.},
	number = {8},
	pages = {999--1020},
	title = {Optimal Dynamic Pricing of Inventories with Stochastic Demand over Finite Horizons},
	url = {https://doi.org/10.1287/mnsc.40.8.999},
	volume = {40},
	year = {1994}}

@incollection{LiYuan:STOC2013,
	author = {Li, Jian and Yuan, Wen},
	booktitle = {S{TOC}'13---{P}roceedings of the 2013 {ACM} {S}ymposium on {T}heory of {C}omputing},
	doi = {10.1145/2488608.2488731},
	mrclass = {68W25 (60C05 62E17 90C15 90C27)},
	mrnumber = {3210859},
	pages = {971--980},
	publisher = {ACM, New York, NY},
	title = {Stochastic combinatorial optimization via {P}oisson approximation},
	url = {https://doi.org/10.1145/2488608.2488731},
	year = {2013}}

@article{Dantzig:OR1957,
	author = {Dantzig, George B.},
	doi = {10.1287/opre.5.2.266},
	fjournal = {Operations Research},
	issn = {0030-364X},
	journal = {Oper. Res.},
	mrclass = {90.0X},
	mrnumber = {0089098},
	mrreviewer = {P. Wolfe},
	number = {2},
	pages = {266--277},
	title = {Discrete-variable extremum problems},
	url = {https://doi.org/10.1287/opre.5.2.266},
	volume = {5},
	year = {1957}}

@article{Marchetti-SpaccamelaVercellis:MP1995,
	author = {Marchetti-Spaccamela, A. and Vercellis, C.},
	doi = {10.1007/BF01585758},
	fjournal = {Mathematical Programming},
	issn = {0025-5610},
	journal = {Math. Program.},
	mrclass = {90C15 (90C27)},
	mrnumber = {1312106},
	mrreviewer = {Anton Volgenant},
	number = {1, Ser. A},
	pages = {73--104},
	title = {Stochastic on-line knapsack problems},
	url = {https://doi.org/10.1007/BF01585758},
	volume = {68},
	year = {1995}}

@article{Lueker:JACIL1998,
	author = {Lueker, George S.},
	doi = {10.1006/jagm.1998.0954},
	fjournal = {Journal of Algorithms. Cognition, Informatics and Logic},
	issn = {0196-6774},
	journal = {J. Algorithms},
	mrclass = {90C27 (68Q25 90C10)},
	mrnumber = {1661604},
	mrreviewer = {Charles E. Blair},
	number = {2},
	pages = {277--305},
	title = {Average-case analysis of off-line and on-line knapsack problems},
	url = {https://doi.org/10.1006/jagm.1998.0954},
	volume = {29},
	year = {1998}}

@article {Vera-Arlotto-Gurvich-Levin:MOOR2025,
    AUTHOR = {Vera, Alberto and Arlotto, Alessandro and Gurvich, Itai and
              Levin, Eli},
     TITLE = {Dynamic resource allocation: the geometry and robustness of
              constant regret},
   JOURNAL = {Math. Oper. Res.},
  FJOURNAL = {Mathematics of Operations Research},
    VOLUME = {50},
      YEAR = {2025},
    NUMBER = {4},
     PAGES = {2834--2872},
      ISSN = {0364-765X,1526-5471},
   MRCLASS = {91B32 (90C05 90C39 90C90)},
  MRNUMBER = {4989848},
       DOI = {10.1287/moor.2020.0334},
       URL = {https://doi.org/10.1287/moor.2020.0334},
}

@article {Vera-Banerjee-Gurvich:OPRE2021,
    AUTHOR = {Vera, Alberto and Banerjee, Siddhartha and Gurvich, Itai},
     TITLE = {Online allocation and pricing: constant regret via {B}ellman
              inequalities},
   JOURNAL = {Oper. Res.},
  FJOURNAL = {Operations Research},
    VOLUME = {69},
      YEAR = {2021},
    NUMBER = {3},
     PAGES = {821--840},
      ISSN = {0030-364X,1526-5463},
   MRCLASS = {90C15 (90B80 90C10 90C39 90C40)},
  MRNUMBER = {4280419},
MRREVIEWER = {Krishna\ Kalyanam},
       DOI = {10.1287/opre.2020.2061},
       URL = {https://doi.org/10.1287/opre.2020.2061},
}

@article{KleywegtPapastavrou:OR2001,
	author = {Kleywegt, Anton J. and Papastavrou, Jason D.},
	coden = {OPREAI},
	doi = {10.1287/opre.49.1.26.11185},
	fjournal = {Operations Research},
	issn = {0030-364X},
	journal = {Oper. Res.},
	mrclass = {90C27 (90C15 90C40)},
	mrnumber = {1814607 (2001k:90061)},
	mrreviewer = {Ignacy Kaliszewski},
	number = {1},
	pages = {26--41},
	publisher = {INFORMS},
	sjournal = {Oper. Res.},
	title = {The dynamic and stochastic knapsack problem with random sized items},
	volume = {49},
	year = {2001}}

@article{KleywegtPapastavrou:OR1998,
	author = {Kleywegt, Anton J. and Papastavrou, Jason D.},
	coden = {OPREAI},
	doi = {10.1287/opre.46.1.17},
	fjournal = {Operations Research},
	issn = {0030-364X},
	journal = {Oper. Res.},
	mrclass = {90C27 (90C15)},
	mrnumber = {1658657 (99k:90134)},
	mrreviewer = {Guntram Scheithauer},
	number = {1},
	pages = {17--35},
	publisher = {INFORMS},
	sjournal = {Oper. Res.},
	title = {The dynamic and stochastic knapsack problem},
	volume = {46},
	year = {1998}}

@article{McLeish:AP1974,
	author = {McLeish, D. L.},
	doi = {10.1214/aop/1176996608},
	fjournal = {The Annals of Probability},
	journal = {Ann. Probab.},
	month = {08},
	number = {4},
	pages = {620--628},
	publisher = {The Institute of Mathematical Statistics},
	title = {Dependent Central Limit Theorems and Invariance Principles},
	url = {https://doi.org/10.1214/aop/1176996608},
	volume = {2},
	year = {1974}}

@book{Bellman:PUP1957,
	author = {Bellman, Richard},
	mrclass = {90.0X},
	mrnumber = {0090477},
	mrreviewer = {J. Kiefer},
	pages = {xxv+342},
	publisher = {Princeton University Press, Princeton, N. J.},
	title = {Dynamic programming},
	year = {1957}}

@book{Powell:JWS2011,
	author = {Powell, Warren B.},
	doi = {10.1002/9781118029176},
	edition = {Second},
	isbn = {978-0-470-60445-8},
	mrclass = {90-01 (90C15 90C39)},
	mrnumber = {2839330},
	pages = {xviii+627},
	publisher = {John Wiley \& Sons, Inc., Hoboken, NJ},
	series = {Wiley Series in Probability and Statistics},
	title = {Approximate Dynamic Programming: Solving the Curses of Dimensionality},
	url = {https://doi.org/10.1002/9781118029176},
	year = {2011}}

@book{SidneyResnick2007,
	author = {Resnick, Sidney},
	isbn = {978-0-387-24272-9},
	pages = {XIX+404},
	publisher = {Springer New York},
	series = {Springer Series in Operations Research and Financial Engineering},
	title = {Heavy-Tail Phenomena: Probabilistic and Statistical Modeling},
	year = {2007},
    doi = {10.1007/9780387450247},
	edition = {first}}

@article{JasinKumar:MOR2012,
author = {Jasin, Stefanus and Kumar, Sunil},
title = {A Re-Solving Heuristic with Bounded Revenue Loss for Network Revenue Management with Customer Choice},
journal = {Mathematics of Operations Research},
volume = {37},
number = {2},
pages = {313-345},
year = {2012},
doi = {10.1287/moor.1120.0537}
}

@article{Golrezaei:MS2014,
author = {Golrezaei, Negin and Nazerzadeh, Hamid and Rusmevichientong, Paat},
title = {Real-Time Optimization of Personalized Assortments},
journal = {Management Science},
volume = {60},
number = {6},
pages = {1532-1551},
year = {2014},
doi = {10.1287/mnsc.2014.1939}
}

@article{Goyal:OR2025,
author = {Goyal, Vineet and Iyengar, Garud and Udwani, Rajan},
title = {Asymptotically Optimal Competitive Ratio for Online Allocation of Reusable Resources},
journal = {Operations Research},
volume = {73},
number = {4},
pages = {1897-1915},
year = {2025},
doi = {10.1287/opre.2021.0695}
}

@article{JiangMaZhang:OR2025-degeneracy-is-ok,
author = {Jiang, Jiashuo and Ma, Will and Zhang, Jiawei},
title = {Degeneracy Is OK: Logarithmic Regret for Network Revenue Management with Indiscrete Distributions},
journal = {Operations Research},
volume = {73},
number = {6},
pages = {3405-3420},
year = {2025},
doi = {10.1287/opre.2022.0641}
}

@article{Zhang:arXiv2026-resource-allocation,
      title={Online Resource Allocation with Continuous Random Consumption: Regret under Degeneracy}, 
      author={Jiawei Zhang},
      journal={arXiv e-print 2607.02196},
      year={2026},
      eprint={2607.02196},
      archivePrefix={arXiv},
      primaryClass={cs.LG},
      url={https://arxiv.org/abs/2607.02196}
}

@article{Jiang:arXiv2023,
      title={Constant Approximation for Network Revenue Management with Markovian-Correlated Customer Arrivals}, 
      author={Jiashuo Jiang},
      journal={arXiv e-print 2305.05829},
      year={2023},
      eprint={2305.05829},
      archivePrefix={arXiv},
      primaryClass={cs.DS},
      url={https://arxiv.org/abs/2305.05829}
}

@article{Zhang:arXiv2026-multi-secretary,
      title={Tight Lower Bounds for the Multi-Secretary Problem via Bellman Certificates}, 
      author={Jiawei Zhang},
      journal={arXiv e-print 2607.02150},
      year={2026},
      eprint={2607.02150},
      archivePrefix={arXiv},
      primaryClass={cs.DS},
      url={https://arxiv.org/abs/2607.02150}
}
    \bibliographystyle{apalike}
    \par
}

% \begin{APPENDICES}
% \section{<Title of Section A>}
% \section{<Title of Section B>}
% etc
% \end{APPENDICES}

\newpage 

\begin{APPENDICES}

\section{Preliminaries on slowly  and regularly varying functions}
\label{apx:regularly-varying}

In this appendix, we review results from the theory of slowly  and regularly varying 
functions that we use in our analysis \citep[see, e.g.,][]{BinghamGoldieTeugels:CUP1989}.
The connection with this literature comes from the fact that---as discussed in 
Section~\ref{sec:typical-distribution-properties}---the density of a typical distribution 
(see Definition~\ref{def:typical-class}) is regularly varying at zero.
Accordingly, we state all results in that form rather than the equivalent form 
at infinity that is common in the literature.

The following lemma provides a convenient test for regular variation when $g$ is positive 
and differentiable near zero.

\begin{lemma}[{\citealp[Theorem 2]{Lamperti:TAMS1958}}]\label{lm:limit-condition-regularly-varying}
If $g$ is a positive and differentiable function on $(0, U) \subseteq \R_+$ such that
\begin{equation*}
\lim_{u \rightarrow 0^+} \frac{ug'(u)}{g(u)} = \rho,
\end{equation*}
then $g$ is regularly varying at zero with index $\rho$.
That is, $g(u)=u^{\rho} \ell(u)$ for some slowly varying function $\ell$.
\end{lemma}

% AA: this is also \citealp[Exercise 1.11.13]{BinghamGoldieTeugels:CUP1989} but with continuously differentiable as an assumption}

The next result will play a crucial role in our analysis, especially in the 
context of typical densities. 
It shows that the ratio of the integral of a regularly 
varying function to a suitable rescaling of the function itself converges to an explicit 
constant.

% AA: the result below is also in {see also \citealp[Proposition 1.5.8]{BinghamGoldieTeugels:CUP1989}}, \citealp[Lemma 1]{Csorgo-Mason:AP1986}}

\begin{lemma}[Karamata's Theorem, {\citealp[see, e.g.,][Theorem 2.1]{SidneyResnick2007}}]\label{lm:Karamata-RV}
 Suppose $g$ is regularly varying at zero with index $\rho > -1$. 
Then $G(u)=\int_0^u g(z)dz$ is regularly varying at zero with index $\rho+1$ and 
$$
\lim_{u \rightarrow 0^+} \frac{\int_0^u g(z) \, dz }{u g(u)} = \frac{1}{\rho + 1}.
$$
\end{lemma}

\section{Proofs of the results in Section~\ref{sec:typical-distribution-properties}}
\label{apx:typical-distribution-properties-proofs}

\proof{Proof of Lemma~\ref{lem:linear-approximation-second-moment-bounds}.}
Fix $h \in (0, \wbar]$ and set
\(
    D(w) = F(w) - \frac{ wf(w) }{ r }.
\)
For $p = 1, 2$, integration by parts and the substitution of
$F(w) = wf(w)/r + D(w)$ give
\[
    \int_0^h w^p f(w) \, dw
    =
    \frac{ r }{ r+p } h^p F(h)
    -
    \frac{ pr }{ r+p }
    \int_0^h w^{p-1} D(w) \, dw.
\]
Thus, expanding the square yields
\begin{align*}
    \int_0^h
        \left( 1 - \frac{ w }{ h } \right)^2 f(w) \, dw
    & =
    F(h)
    - \frac{ 2 }{ h }
        \int_0^h wf(w) \, dw
    + \frac{ 1 }{ h^2 }
        \int_0^h w^2 f(w) \, dw \\
    & = 
    \left( 1 - \frac{ 2r }{ r+1 } + \frac{ r }{ r+2 } \right)
    F(h) 
    + \rho_h,
\end{align*}
where 
\[
    1 - \frac{ 2r }{ r+1 } + \frac{ r }{ r+2 }  = \sigma^2
    \qquad \text{and} \qquad
    \rho_h
    =
    \frac{ 2r }{ h }
    \int_0^h
        \left(
            \frac{ 1 }{ r+1 }
            - \frac{ w }{ h(r+2) }
        \right)
        D(w) \, dw.
\]

It remains to bound $\rho_h$.
For $0 < w \leq h \leq \wbar$, the rate condition
\eqref{eq:typical-class-limit-with-rate} gives
\[
    \abs{ D(w) }
    \leq M_F w^{3/4} wf(w).
\]
Moreover, since $r/(r+1) \leq 1$ and
$(r+1)/(r+2) \geq 1/2$, we have
\[
    0
    \leq
    \frac{ 2rw }{ h }
    \left(
        \frac{ 1 }{ r+1 }
        - \frac{ w }{ h(r+2) }
    \right)
    \leq
    \frac{ 2w }{ h }
    \left( 1 - \frac{ w }{ 2 h } \right)
    \leq 1.
\]
Consequently,
\[
    \abs{ \rho_h }
    \leq
    M_F \int_0^h
        \frac{ 2rw }{ h }
        \left(
            \frac{ 1 }{ r+1 }
            - \frac{ w }{ h(r+2) }
        \right)
        w^{3/4} f(w) \, dw 
     \leq
    M_F \int_0^h w^{3/4} f(w) \, dw
    \leq
    M_F h^{3/4} F(h).
\]
Combining this bound with the preceding identity proves both
inequalities in \eqref{eq:linear-approximation-second-moment-bounds}.
\halmos
\endproof

\proof{Proof of Proposition~\ref{prop:typical-density-properties}.}
\ref{property:wfw-increasing}
        This property is an immediate consequence of the monotonicity condition \eqref{eq:typical-class-monotonicity-condition}.
        From \eqref{eq:typical-class-monotonicity-condition}, we have 
        $wf'(w) + f(w) \geq 0$ for all $w \in (0,\wbar]$, so 
        $\left( wf(w) \right)' = wf'(w) + f(w) \geq 0$ on $(0,\wbar]$
        and the map $w \mapsto wf(w)$ is increasing on the same interval.

\ref{property:F-regularly-varying}
        Because $F$ is a weight distribution with density $f$, the rate condition \eqref{eq:typical-class-limit-with-rate} gives us the limit
        \begin{equation}\label{eq:intermediate-limits-weightdist}
        \lim_{w \to 0^+} \frac{w F'(w)}{F(w)} = \lim_{w \to 0^+} \frac{w f(w)}{F(w)} = r.
        \end{equation}
        Both $f, F > 0$ on $(0,\wbar]$, so Lemma~\ref{lm:limit-condition-regularly-varying} implies 
        that $F$ is regularly varying at zero with index $r$.

\ref{property:f-regularly-varying}
        By part~\ref{property:F-regularly-varying}, we may write $F(w) = w^{r} L(w)$ 
        for some slowly varying $L$. Thus, 
        $$
        f(w) = \frac{F(w)}{w} \cdot \frac{w f(w)}{F(w)} = w^{r-1} \ell(w),
        \qquad \mbox{where} \qquad 
        \ell(w) = L(w) \cdot \frac{w f(w)}{F(w)}.
        $$
        The map $\ell$ is continuous and positive on $(0,\wbar]$, and it is slowly varying 
        because it is the product of the slowly varying function $L$ 
        and a function converging to the positive constant $r$ by 
        \eqref{eq:intermediate-limits-weightdist}.
\halmos\endproof

\proof{Proof of Proposition~\ref{prop:convergence-rate-bound}.}
Fix $q \geq -1$ and $w\in(0,\wbar]$.
Multiplying the rate condition
\eqref{eq:typical-class-limit-with-rate} by $y^{q+1}f(y)$ gives
\[
\left|
y^q F(y) - \frac{1}{r} y^{q + 1}f(y)
\right|
\leq M_F y^{3/4} y^{q+1} f(y).
\]
Since $q + 1\geq0$, the functions $y^{q + 1}f(y)$ and $y^q F(y)$ are both integrable over $(0,w)$.
Integrating and using $y^{3/4}\leq w^{3/4}$ yields
\[
\abs{
\int_0^w y^q F(y)\,dy
-\frac{1}{r}\int_0^w y^{q+1}f(y)\,dy}
\leq
M_Fw^{3/4}\int_0^w y^{q+1}f(y)\,dy.
\]
Integration by parts gives 
\[
\int_0^w y^{q+1}f(y)\,dy
=
w^{q+1}F(w)
-(q+1)\int_0^w y^q F(y)\,dy.
\]
Substituting this identity and rearranging, we obtain
\begin{align*}
\abs{
\int_0^w y^q F(y)\,dy
-\frac{1}{r}\int_0^w y^{q+1}f(y)\,dy}
& = 
\abs{
\int_0^w y^q F(y)\,dy
-\frac{1}{r}\left[
w^{q+1}F(w)
-(q+1)\int_0^w y^q F(y)\,dy
\right]}\\
& = 
\abs{
\frac{r+q+1}{r}\int_0^w y^q F(y)\,dy
-\frac{1}{r}w^{q+1}F(w)} \\
& =
\frac{r+q+1}{r}
\abs{
\int_0^w y^q F(y)\,dy
-\frac{w^{q+1}F(w)}{r+q+1}}.
\end{align*}
The preceding inequality therefore implies
\[
\abs{
\frac{\int_0^w y^q F(y)\,dy}{w^{q+1}F(w)}
-\frac{1}{r+q+1}}
\leq
\frac{r}{r+q+1}M_F w^{3/4}
\frac{\int_0^w y^{q+1}f(y)\,dy}
     {w^{q+1}F(w)}
\leq
\frac{r}{r+q+1} M_Fw^{3/4},
\]
since $y^{q+1}\leq w^{q+1}$ on $(0,w)$. Finally, $q + 1 \geq 0$, so $\frac{r}{r+q+1} \leq 1$, and the desired result follows.
\halmos\endproof

\proof{Proof of Corollary~\ref{cor:karamata-ratios}.}
For $0 < \epsilon \leq \wbar$, let
\(
    I(\epsilon) = \int_0^\epsilon w f(w) \, dw.
\)
Integration by parts gives
\[
    \frac{ I(\epsilon) }{ \epsilon F(\epsilon) }
    =
    1 - \frac{ \int_0^\epsilon F(w) \, dw }{ \epsilon F(\epsilon) }.
\]
By \eqref{eq:error-bound-of-limits} with $q = 0$, the right-hand side
converges to $1 - \frac{ 1 }{ r+1 } = \frac{ r }{ r+1 }$,
proving the second limit.
For the first limit, write
\[
    \frac{ I(\epsilon) }{ \epsilon^2 f(\epsilon) }
    =
    \frac{ I(\epsilon) }{ \epsilon F(\epsilon) }
    \frac{ F(\epsilon) }{ \epsilon f(\epsilon) }.
\]
The first factor converges to $r/(r+1)$, while the second converges
to $1/r$ by \eqref{eq:typical-class-limit-with-rate}.
Their product therefore converges to $1/(r+1)$.

To establish the lower bounds in
\eqref{eq:karamata-ratios-pointwise-bounds},
note that both ratios are positive and continuous on $(0,\wbar]$,
with positive finite limits at zero.
They are therefore bounded away from zero on $(0,\wbar]$,
which gives a common constant $C_F > 1$ for the lower bounds.
For the upper bounds, simply note that the monotonicity of the map
$w \mapsto w f(w)$ on $(0,\wbar]$
(Proposition~\ref{prop:typical-density-properties}--\ref{property:wfw-increasing})
gives $I(\epsilon) \leq \epsilon^2 f(\epsilon)$,
while the inequality $w \leq \epsilon$ gives
$I(\epsilon) \leq \epsilon F(\epsilon)$.

To establish \eqref{eq:karamata-integral-ratio-bounds}, the left estimate in
\eqref{eq:karamata-ratios-pointwise-bounds} gives
\[
    \frac{ 1 }{ \epsilon }
    \leq
    \frac{ I'(\epsilon) }{ I(\epsilon) }
    =
    \frac{ \epsilon f(\epsilon) }{ I(\epsilon) }
    \leq
    \frac{ C_F }{ \epsilon }.
\]
Integrating the inequalities above between $h$ and $\epsilon$, where
$0 < h < \epsilon \leq \wbar$, and exponentiating yields
\[
    \left( \frac{ h }{ \epsilon } \right)^{C_F}
    \leq
    \frac{ I(h) }{ I(\epsilon) }
    \leq
    \frac{ h }{ \epsilon },
\]
just as needed.
\halmos\endproof

\section{Proofs of the results in Section~\ref{sec:consumption-function-properties}}
\label{apx:consumption-function-properties-proofs}

\proof{Proof of Proposition~\ref{prop:basic-properties-epsilon}.}
\ref{property:consumption-function-concavity}
If $k\geq K_1$, the definition of $K_1$ \eqref{eq:K1} tells us 
that $\epsilon_k(x)\leq\wbar$ for all $x\in[0,c]$. 
By Proposition~\ref{prop:typical-density-properties}--\ref{property:wfw-increasing} and the monotonicity of the map $x \mapsto \epsilon_k(x)$,
it follows that the composition $\epsilon_k(x)f(\epsilon_k(x))$ is increasing in $x$. 
Consequently, the first derivative \eqref{eq:consumption-function-derivative}, given by 
$\epsilon_k'(x) = \frac{1}{k\epsilon_k(x)f(\epsilon_k(x))}$,
is decreasing in $x$ for all $k \geq K_1$.
Hence, $\epsilon_k(x)$ is concave for $x \in [0,c]$ when $k\geq K_1$. 
Finally, since $\epsilon_k(x) > 0$ for $x\in(0,c]$, the  convexity of $1/\epsilon_k(x)$ for $k \geq K_1$ follows immediately.
For $0<y\leq x\leq c$, concavity and $\epsilon_k(0)=0$
give
\[
\epsilon_k(y)
\geq
\frac{y}{x}\epsilon_k(x)
+\left(1-\frac{y}{x}\right)\epsilon_k(0)
=
\frac{y}{x}\epsilon_k(x).
\]
Dividing by $y$ on both sides proves that $\epsilon_k(x)/x$
is decreasing.

\ref{property:consumption-function-regularly-varying}
For $k \geq K_1$, the definition of $K_1$ in \eqref{eq:K1} ensures that the integral
representation \eqref{eq:consumption-function-integral-representation} and the derivative \eqref{eq:consumption-function-derivative} hold for all
$x \in (0,c]$, where $\epsilon_k$ is also positive.
Since $\epsilon_k$ is continuous with $\epsilon_k(0) = 0$, 
we have $\epsilon_k(x) \to 0$ as $x \to 0^+$, and therefore
$$
\lim_{x \to 0^+} \frac{x \epsilon_k'(x)}{\epsilon_k(x)}
= \lim_{x \to 0^+} \frac{x}{k \epsilon_k(x)^2 f(\epsilon_k(x))}
= \lim_{x \to 0^+} \frac{\int_0^{\epsilon_k(x)} y f(y) \, dy}{\epsilon_k(x)^2 f(\epsilon_k(x))}
= \lim_{w \to 0^+} \frac{\int_0^{w} y f(y) \, dy}{w^2 f(w)}
= \frac{1}{r+1},
$$
where the third equality is the change of variables $w = \epsilon_k(x)$ and the last one is the
left limit in \eqref{eq:karamata-ratios-limit}. Lemma~\ref{lm:limit-condition-regularly-varying} then implies that $\epsilon_k$ is regularly varying at zero with index $\frac{1}{r+1}$ for all $k \geq K_1$.
Finally, the space-time scaling equivalence \eqref{eq:consumption-function-space-time-scaling} and the regularly varying property of $\epsilon_{K_1}(x)$ 
together imply the existence of a function $L$ that is slowly varying at zero 
such that
$$
\epsilon_k(x) 
= \epsilon_{K_1}\left(\frac{x}{k} K_1\right) 
= \left(\frac{x}{k}\right)^{\frac{1}{r+1}} L\left( \frac{x}{k} \right) 
\qquad \text{for all  $k \geq K_1$ and $x \in [0,c]$}. 
$$

\ref{property:useful-eps'/eps^2-decreasing}
By Proposition~\ref{prop:basic-properties-epsilon}--\ref{property:consumption-function-concavity}, $\epsilon_k(x)$ is concave in $x$, so $\epsilon_k'(x)$ is nonnegative and decreasing in $x$.
Since $\epsilon_k(x)$ is positive and increasing on $(0,c]$, $\epsilon_k^{-a}(x)$ is positive and decreasing.
Therefore, their product $\epsilon_k'(x)/\epsilon_k^a(x)$ is decreasing, as expected.

\ref{property:ratio-eps-x-to-x-minus-eps}
Fix $k \geq K_1$ and $0 \leq w \leq \epsilon_k(x) < x \leq c$.
By Proposition~\ref{prop:basic-properties-epsilon}--\ref{property:consumption-function-concavity},
the map $y \mapsto \epsilon_k(y)/y$ is decreasing on $(0,c]$,
so
\(
    \frac{ \epsilon_k(x) }{ x }
    \leq
    \frac{ \epsilon_k(x-w) }{ x-w }.
\)
All quantities in this inequality are positive, so rearranging gives
the first inequality of \eqref{eq:ratio-eps-x-to-x-minus-eps}.
The second inequality of \eqref{eq:ratio-eps-x-to-x-minus-eps}
then follows because $w \leq \epsilon_k(x)$.

Next, note that 
\[
    \frac{ 1 }{ (x-w)^{3/4} \epsilon_k(x-w)^{1/4} }
     =
    \left( \frac{ x }{ x-w } \right)^{3/4}
    \left( \frac{ \epsilon_k(x) }{ \epsilon_k(x-w) } \right)^{1/4}
    \frac{ 1 }{ x^{3/4} \epsilon_k(x)^{1/4} }.
\]
Applying the first inequality of
\eqref{eq:ratio-eps-x-to-x-minus-eps} and then using
$w \leq \epsilon_k(x)$ gives
\[
    \frac{ 1 }{ (x-w)^{3/4} \epsilon_k(x-w)^{1/4} }
    \leq
    \frac{ x }{ x-w }
    \frac{ 1 }{ x^{3/4} \epsilon_k(x)^{1/4} } 
    \leq
    \frac{ x }{ x-\epsilon_k(x) }
    \frac{ 1 }{ x^{3/4} \epsilon_k(x)^{1/4} },
\]
which proves \eqref{eq:ratio-eps-x-to-x-minus-eps-3/4}.

\ref{property:1/epsilon-integral-bounds} 
Since $\epsilon_k$ is increasing and positive on $(0,c]$, the map $y \mapsto 1/\epsilon_k(y)$ is
decreasing there, so bounding the integrand by its values at the endpoints of $[x-w,x]$ gives
$$
\frac{w}{\epsilon_k(x)} \leq \int_{x-w}^{x} \frac{1}{\epsilon_k(y)} \, dy
\leq \frac{w}{\epsilon_k(x-w)}.
$$
The left inequality of \eqref{eq:1/epsilon-integral-bounds} is immediate from the
first of these bounds.
For the right inequality, we subtract $\frac{w}{\epsilon_k(x)}$ from the second bound to obtain
$$
\int_{x-w}^{x} \frac{1}{\epsilon_k(y)} \, dy - \frac{w}{\epsilon_k(x)}
\leq w \left( \frac{1}{\epsilon_k(x-w)} - \frac{1}{\epsilon_k(x)} \right).
$$
For $k \geq K_1$, \eqref{eq:ratio-eps-x-to-x-minus-eps} implies that for $0 \leq w \leq \hhat_k(x)$, 
$\frac{1}{\epsilon_k(x-w)} \leq \frac{x}{x-w} \frac{1}{\epsilon_k(x)}$ after rearrangement.
Using this last inequality in the preceding display, we find
$$
\int_{x-w}^{x} \frac{1}{\epsilon_k(y)} \, dy - \frac{w}{\epsilon_k(x)}
\leq w \left( \frac{x}{x-w} - 1 \right) \frac{1}{\epsilon_k(x)}
= \frac{w^2}{(x-w)\epsilon_k(x)},
$$
completing the proof.

\ref{property:taylor-first-order-eps-a}
Fix $a > 0$, $k \geq K_1$, $x \in (0,c]$, and
$w \in [0,x)$.
If $w = 0$, all three expressions in
\eqref{eq:taylor-first-order-eps-a} are zero, so the result
is immediate.
Next, suppose that $w > 0$.
Since $\epsilon_k$ is positive and continuous on $[x-w,x]$
and differentiable on $(x-w,x)$, the mean value theorem
applied to $y \mapsto 1/\epsilon_k(y)^a$ gives some
$z \in (0,w)$ such that
\[
    \frac{ 1 }{ \epsilon_k(x-w)^a }
    - \frac{ 1 }{ \epsilon_k(x)^a }
    =
    a w \cdot
    \frac{ \epsilon_k'(x-z) }{ \epsilon_k(x-z)^{a+1} }.
\]
By Proposition~\ref{prop:basic-properties-epsilon}--\ref{property:useful-eps'/eps^2-decreasing},
applied with exponent $a+1$, the map
$y \mapsto \epsilon_k'(y)/\epsilon_k(y)^{a+1}$ is decreasing.
Since $x-w \leq x-z \leq x$, we therefore have
\[
    \frac{ \epsilon_k'(x) }{ \epsilon_k(x)^{a+1} }
    \leq
    \frac{ \epsilon_k'(x-z) }{ \epsilon_k(x-z)^{a+1} }
    \leq
    \frac{ \epsilon_k'(x-w) }{ \epsilon_k(x-w)^{a+1} }.
\]
Multiplying by $a w$ and using the preceding identity proves
\eqref{eq:taylor-first-order-eps-a}.
\halmos\endproof

\proof{Proof of Proposition~\ref{prop:joint-properties-f-epsilon}.}
\ref{property:fluid-benchmark-growth-relative-to-logT}
For $T\geq K_1$, monotonicity of the consumption function and \eqref{eq:epsilon-k-leq-wbar} give $0<\epsilon_T(c)\leq\epsilon_{K_1}(c)\leq\wbar$.
We use the lower bound in \eqref{eq:karamata-integral-ratio-bounds} and the integral representation \eqref{eq:consumption-function-integral-representation} to obtain
\[
\left(\frac{\epsilon_T(c)}{\epsilon_{K_1}(c)}\right)^{C_F}
\leq \frac{\int_0^{\epsilon_T(c)} w f(w) \, dw}{\int_0^{\epsilon_{K_1}(c)} w f(w) \, dw} = \frac{c/T}{c/K_1}
= \frac{K_1}{T}.
\]
Combining this inequality with the bound
$c / \epsilon_T(c) \leq T F(\epsilon_T(c))$ from the upper bound
on the second ratio in
\eqref{eq:karamata-ratios-pointwise-bounds} gives
\[
\{TF(\epsilon_T(c))\}^{1/4}
\geq \left(\frac{c}{\epsilon_T(c)}\right)^{1/4}
\geq \left(\frac{c}{\epsilon_{K_1}(c)}\right)^{1/4}
     \left(\frac{T}{K_1}\right)^{1/(4C_F)}.
\]
Since $\log T=o(T^\alpha)$ for any $\alpha>0$, the desired limit follows.

\ref{property:concavity-of-power-of-F-epsilon}
For $k \geq K_1$ and $x \in (0,c)$, the derivative identity
\eqref{eq:consumption-function-derivative} gives
\(
\frac{d}{dx} F(\epsilon_k(x)) = f(\epsilon_k(x)) \epsilon_k'(x) = \frac{1}{k \epsilon_k(x)}.
\)
Differentiating once more yields
\(
\frac{d^2}{dx^2} F(\epsilon_k(x)) = - \frac{1}{k^2 \epsilon_k(x)^3 f(\epsilon_k(x))} \leq 0, 
\)
which proves concavity.
Note that $F(\epsilon_k(0)) = 0$, so the concavity of $F(\epsilon_k(x))$ implies that $F(\epsilon_k(x))/x$ is decreasing, by the same argument used in Proposition~\ref{prop:basic-properties-epsilon}--\ref{property:consumption-function-concavity}.

\ref{property:epsilon'/epsilon^2-delta-bounds}
Recall from \eqref{eq:consumption-function-derivative} that
$\epsilon'_k(x) = \frac{1}{k\epsilon_k(x) f(\epsilon_k(x))}$, so that for any $a > 0$ and
$x \in (0,c]$,
$$
\frac{\epsilon_k'(x)}{\epsilon_k(x)^{a+1}}
= \frac{x}{k \epsilon_k(x)^2 f(\epsilon_k(x))} \cdot \frac{1}{x \epsilon_k(x)^a}.
$$
Note that the left estimate in
\eqref{eq:karamata-ratios-pointwise-bounds} at $\epsilon = \epsilon_k(x) \leq \wbar$ gives the lower and upper bounds for the first factor for all $k \geq K_1$ and $x\in(0,c]$.
Multiplying that two-sided bound by the second factor gives the claim.

\ref{property:local-reciprocal-increments}
Fix $k \geq K_1$ and $x \in (0, c]$ satisfying
$\epsilon_k(x) \leq x/2$.
The lower bound in \eqref{eq:difference-of-consumption-function-reciprocals-linear-approximation}
follows immediately from \eqref{eq:taylor-first-order-eps-a}
with exponent $a = 1$.
For the upper bound, we fix $w \in [0, \epsilon_k(x)]$, let $s \in [0,w]$, and obtain an upper bound for
$\epsilon_k'(x-s) / \epsilon_k(x-s)^2$ and a corresponding
lower bound for $\epsilon_k'(x) / \epsilon_k(x)^2$.
We then subtract these bounds and integrate to estimate
the reciprocal increment.
Since $k \geq K_1$ and $\epsilon_k$ is increasing, we have
\(
    0 < \epsilon_k(x-s) \leq \epsilon_k(x) \leq \wbar.
\)
By \eqref{eq:ratio-eps-x-to-x-minus-eps} and the monotonicity of $y\mapsto F(\epsilon_k(y))/y$ from part~\ref{property:concavity-of-power-of-F-epsilon}, we have
\begin{equation}\label{eq:concavity-lower-bounds}
    \left(1-\frac{s}{x}\right)\epsilon_k(x)
    \leq 
    \epsilon_k(x-s)
    \qquad \text{and} \qquad
    \left(1-\frac{s}{x}\right)F(\epsilon_k(x))
    \leq
    F(\epsilon_k(x-s)).
\end{equation}
If we now use the upper bound in part~\ref{property:epsilon'/epsilon^2-delta-bounds} of this proposition with $a = 1$ at $x-s$, together with the left
inequality above, we find
\begin{equation}
\label{eq:epsilon'/epsilon^2-bound}
    \frac{ \epsilon_k'(x-s) }{ \epsilon_k(x-s)^2 }
    \leq
    \frac{ 1 }{ (x-s) \epsilon_k(x-s) }
    \leq
    \frac{ x }{ (x-s)^2 \epsilon_k(x) }
    =
    \frac{ 1 }{ x \epsilon_k(x) }
    \left( 1 - \frac{ s }{ x } \right)^{-2}.
\end{equation}

To obtain the corresponding lower bound at $x$, we first note from the derivative identity \eqref{eq:consumption-function-derivative} that 
\[
    F( \epsilon_k(y) ) \epsilon_k'(y)
    =
    \frac{ 1 }{ k }
    \frac{ F( \epsilon_k(y) ) }{ \epsilon_k(y) f( \epsilon_k(y) ) },
\]
and we bound the quantity on the left by using the rate condition \eqref{eq:typical-class-limit-with-rate}.
For $y = x -  s$ and $y = x $ we find 
\[
    F( \epsilon_k(x-s) ) \epsilon_k'(x-s)
    \leq
    \frac{ 1 }{ kr }
    + \frac{ M_F }{ k } \epsilon_k(x-s)^{3/4}
    \qquad \text{and} \qquad
    \frac{ 1 }{ kr }
    - \frac{ M_F }{ k } \epsilon_k(x)^{3/4}
    \leq 
    F( \epsilon_k(x) ) \epsilon_k'(x).
\]
Since $\epsilon_k$ is increasing, we can now replace $\epsilon_k(x-s)$ with $\epsilon_k(x)$ in the first inequality and combine them to obtain that 
\[
    F( \epsilon_k(x-s) ) \epsilon_k'(x-s)
    \leq
    F( \epsilon_k(x) ) \epsilon_k'(x)
    + \frac{ 2 M_F }{ k } \epsilon_k(x)^{3/4}.
\]
If we now divide this inequality by
$\epsilon_k(x)^2 F( \epsilon_k(x) )$,
set 
\[
\Gamma_k(s,x)
= 
\frac{ F( \epsilon_k(x-s) ) \epsilon_k'(x-s) }{ \epsilon_k(x)^2 F( \epsilon_k(x) ) },
\]
and note that  \eqref{eq:consumption-function-integral-representation} implies
\(
    x
    =
    k \int_0^{\epsilon_k(x)} y f(y) \, dy
    \leq
    k \epsilon_k(x) F( \epsilon_k(x) ),
\)
we then obtain 
\[
    \Gamma_k(s,x)
    \leq
    \frac{ \epsilon_k'(x) }{ \epsilon_k(x)^2 }
    +
    \frac{ 2 M_F \epsilon_k(x)^{3/4} }{ k \epsilon_k(x)^2 F( \epsilon_k(x) ) }
    \leq
    \frac{ \epsilon_k'(x) }{ \epsilon_k(x)^2 }
    +
    \frac{ 2 M_F }{ x \epsilon_k(x)^{1/4} }.
\]
Furthermore, the concavity bounds in
\eqref{eq:concavity-lower-bounds} imply
\[
        \Gamma_k(s,x)
         =
        \frac{ F( \epsilon_k(x-s) ) }{ F( \epsilon_k(x) ) }
        \left( \frac{ \epsilon_k(x-s) }{ \epsilon_k(x) } \right)^2
        \frac{ \epsilon_k'(x-s) }{ \epsilon_k(x-s)^2 }
        \geq
        \left( 1 - \frac{ s }{ x } \right)^3
        \frac{ \epsilon_k'(x-s) }{ \epsilon_k(x-s)^2 }.
\]
The upper and lower bounds that we just derived for $\Gamma_k(s,x)$ can now be combined to obtain that 
\[
    \left( 1 - \frac{ s }{ x } \right)^3
    \frac{ \epsilon_k'(x-s) }{ \epsilon_k(x-s)^2 }
    -
    \frac{ 2 M_F }{ x \epsilon_k(x)^{1/4} }
    \leq
    \frac{ \epsilon_k'(x) }{ \epsilon_k(x)^2 }.
\]
If we now subtract $\epsilon_k'(x-s) / \epsilon_k(x-s)^2$ from both sides, 
change sign, and use \eqref{eq:epsilon'/epsilon^2-bound}, we find
\begin{align*}
    \frac{ \epsilon_k'(x-s) }{ \epsilon_k(x-s)^2 }
    - \frac{ \epsilon_k'(x) }{ \epsilon_k(x)^2 }
    & \leq
        \left\{
            1 - \left( 1 - \frac{ s }{ x } \right)^3
        \right\}
        \frac{ \epsilon_k'(x-s) }{ \epsilon_k(x-s)^2 }
        +
        \frac{ 2 M_F }{ x \epsilon_k(x)^{1/4} }
        \notag\\
    & \leq
        \frac{ 1 }{ x \epsilon_k(x) }
        \left\{
            \left( 1 - \frac{ s }{ x } \right)^{-2}
            -
            \left( 1 - \frac{ s }{ x } \right)
        \right\}
        +
        \frac{ 2 M_F }{ x \epsilon_k(x)^{1/4} }.
\end{align*}

If we now integrate the bound above for $s \in [0,w]$ and use the fundamental theorem of calculus, 
we find
\begin{equation}
\label{eq:consumption-function-reciprocal-difference-intermediate}
        \frac{ 1 }{ \epsilon_k(x-w) }
        - \frac{ 1 }{ \epsilon_k(x) }
        - \frac{ \epsilon_k'(x) }{ \epsilon_k(x)^2 } w
        \leq
        \frac{ 1 }{ x \epsilon_k(x) }
        \int_0^w
        \left\{
            \left( 1 - \frac{ s }{ x } \right)^{-2}
            -
            \left( 1 - \frac{ s }{ x } \right)
        \right\} \, ds
        +
        \frac{ 2 M_F w }{ x \epsilon_k(x)^{1/4} }.
\end{equation}
For the integral above, we factor $s/x$ out of the integrand
and bound the remaining factor using $0 \leq s \leq w \leq x/2$.
Hence, we obtain
\[
        \int_0^w
        \left\{
            \left( 1 - \frac{ s }{ x } \right)^{-2}
            -
            \left( 1 - \frac{ s }{ x } \right)
        \right\} \, ds
        =
        \int_0^w
        \left\{
        \left(1-\frac{s}{x}\right)^{-2}
        +\left(1-\frac{s}{x}\right)^{-1}
        +1
        \right\}
        \frac{s}{x}
        \, ds
        \leq 
        \int_0^w
        7 \frac{s}{x}
        \, ds
        =
        \frac{7 w^2}{2x}.
        \]
If we now substitute this bound into
\eqref{eq:consumption-function-reciprocal-difference-intermediate},
use the bound $w / \epsilon_k(x) \leq 1$, rearrange,
and recall that $\epsilon_k(x) \leq x/2$, $x \leq c$,
and $C_\epsilon = 3 + 2 M_F c^{3/4}$, then we find
\[
    \frac{ 1 }{ \epsilon_k(x-w) }
    - \frac{ 1 }{ \epsilon_k(x) }
    - \frac{ \epsilon_k'(x) }{ \epsilon_k(x)^2 } w
    \leq
    \left(
    \frac{7 w}{ 2 x^{1/4} \epsilon_k(x)^{3/4}}
    +
    2 M_F x^{3/4}
    \right)
    \frac{ w }{ x^{7/4} \epsilon_k(x)^{1/4} }
    \leq 
    \frac{ C_\epsilon }{ x^{7/4} \epsilon_k(x)^{1/4} } w.
\]
This proves the upper bound in
\eqref{eq:difference-of-consumption-function-reciprocals-linear-approximation}.
Finally, integrating both bounds in
\eqref{eq:difference-of-consumption-function-reciprocals-linear-approximation} against $f(w) \, dw$
over $[0, \epsilon_k(x)]$ and using the integral representation \eqref{eq:consumption-function-integral-representation} of the consumption function 
gives \eqref{eq:difference-of-consumption-function-reciprocals-dFintegral}.
\halmos\endproof

\proof{Proof of Proposition~\ref{prop:properties-of-xhat}.}
\ref{property:fixed-points-decreasing}
Since $\epsilon_k(x)$ is decreasing in $k$, the sets over which
the maximum is taken in \eqref{eq:xk-definition} are decreasing
in $k$.
Consequently, $\xhat_{k+1} \leq \xhat_k$ for every $k \geq K_1$.
Moreover, the definition of $\xhat_k$, the monotonicity of
$\epsilon_k$, and \eqref{eq:epsilon-k-leq-wbar} give
\[
    0 \leq \xhat_k
    \leq \epsilon_k(\Theta \xhat_k)
    \leq \epsilon_k(\Theta c)
    \leq \wbar
    \qquad \text{for all $k \geq K_1$}.
\]
The integral representation
\eqref{eq:consumption-function-integral-representation} gives
\[
    \int_0^{\epsilon_k(\Theta c)} w f(w) \, dw
    =
    \frac{ \Theta c }{ k }
    \longrightarrow 0.
\]
Since $f$ is positive on $(0,\wbar]$, it follows that
$\epsilon_k(\Theta c) \to 0$, and hence $\xhat_k \to 0$,
just as needed.

\ref{property:consumption-bound-below-xhat}
Fix $k \geq K_1$.
Both inequalities are immediate when $x = 0$.
Note that for $0 < x \leq \xhat_k$, \eqref{eq:epsilon-k-leq-wbar} and the definition \eqref{eq:xk-definition} give $\xhat_k \leq \epsilon_k(\Theta \xhat_k) \leq \wbar$.
Consequently, the integral representation
\eqref{eq:consumption-function-integral-representation}
implies
\[
    \int_0^{\xhat_k} w f(w) \, dw
    \leq
    \int_0^{\epsilon_k(\Theta \xhat_k)} w f(w) \, dw
    =
    \frac{ \Theta \xhat_k }{ k }.
\]
Since $w \mapsto w f(w)$ is increasing on $(0, \wbar]$,
a change of variable gives
\[
    \frac{ 1 }{ x } \int_0^x w f(w) \, dw
    =
    \int_0^1 xt f(xt) \, dt
    \leq
    \int_0^1 \xhat_k t f(\xhat_k t) \, dt
    =
    \frac{ 1 }{ \xhat_k }
    \int_0^{\xhat_k} w f(w) \, dw
    \leq
    \frac{ \Theta }{ k }.
\]
Multiplying both sides by $x$ proves the first inequality of \eqref{eq:consumption-bound-below-xhat}.
The second inequality then  follows from the definition
\eqref{eq:consumption-function-definition} of $\epsilon_k$.

\ref{property:epsilon/x-ratio-bounds}
Fix $k \geq K$ and write
$I(\epsilon) = \int_0^\epsilon w f(w) \, dw$.
The fixed-point identity \eqref{eq:xk-fixedpoint-property}
and the integral representation
\eqref{eq:consumption-function-integral-representation} give
\[
    \frac{ I(\epsilon_k(\xhat_k)) }{ I(\xhat_k) }
    =
    \frac{ I(\epsilon_k(\xhat_k)) }{ I(\epsilon_k(\Theta \xhat_k)) }
    =
    \frac{ \xhat_k / k }{ \Theta \xhat_k / k }
    =
    \frac{ 1 }{ \Theta }.
\]
Consequently, when we use the bounds
\eqref{eq:karamata-integral-ratio-bounds} from
Corollary~\ref{cor:karamata-ratios} with
$h = \epsilon_k(\xhat_k)$ and $\epsilon = \xhat_k$,
we find
\[
    \left(
        \frac{ \epsilon_k(\xhat_k) }{ \xhat_k }
    \right)^{C_F}
    \leq
    \frac{ 1 }{ \Theta }
    \leq
    \frac{ \epsilon_k(\xhat_k) }{ \xhat_k },
    \qquad \text{so} \qquad
    \frac{ 1 }{ \Theta }
    \leq
    \frac{ \epsilon_k(\xhat_k) }{ \xhat_k }
    \leq
    \left( \frac{ 1 }{ \Theta } \right)^{1/C_F}.
\]
By Proposition~\ref{prop:basic-properties-epsilon}--\ref{property:consumption-function-concavity},
the map $x \mapsto \epsilon_k(x)/x$ is decreasing on $(0,c]$.
For $x \in (0,\xhat_k)$, we thus obtain
\[
    \frac{ \epsilon_k(x) }{ x }
    \geq
    \frac{ \epsilon_k(\xhat_k) }{ \xhat_k }
    \geq
    \frac{ 1 }{ \Theta }
    = \thetalb.
\]
Similarly, for $x \in [\xhat_k,c]$,
\[
    \frac{ \epsilon_k(x) }{ x }
    \leq
    \frac{ \epsilon_k(\xhat_k) }{ \xhat_k }
    \leq
    \left(\frac{ 1 }{ \Theta }\right)^{1/C_F}
    =
    \thetaub.
\]
Finally, $C_F > 1$ and $\Theta > 1$ imply that
$0 < \thetalb < \thetaub < 1$, while
$\Theta > 2^{C_F}$ implies that $\thetaub < 1/2$.
\halmos\endproof

\proof{Proof of Lemma~\ref{lemma:density-lower-bound-large-capacity}.}
First, we note that $\lambda_\Theta \in (0,1)$ since
$\Theta > 1$, $C_F > 1$, and
\(
    \lambda_\Theta
    = \thetaub^{1/(C_F - 1)}
    = \Theta^{-\frac{ 1 }{ C_F(C_F - 1) }}.
\)
Next, fix $k \geq K$, $x \in [\xhat_k, c]$, and
$w \in [\lambda_\Theta \epsilon_k(x), \epsilon_k(x)]$.
Since $\epsilon_k(x) \leq \wbar$, the integral representation
\eqref{eq:consumption-function-integral-representation},
the first inequality in \eqref{eq:karamata-integral-ratio-bounds},
and the monotonicity of $y \mapsto y f(y)$ give
\[
    \left( \frac{ w }{ \epsilon_k(x) } \right)^{C_F}
    \frac{ x }{ k }
    =
    \left( \frac{ w }{ \epsilon_k(x) } \right)^{C_F}
    \int_0^{\epsilon_k(x)} y f(y) \, dy
    \leq
    \int_0^w y f(y) \, dy
    \leq
    w^2 f(w).
\]
Hence, \eqref{eq:density-lower-bound-large-capacity} follows if we can show that
\[
    1
    \leq
    \left( \frac{ w }{ \epsilon_k(x) } \right)^{C_F - 1}
    \frac{ x }{ w }
    \qquad
    \text{for all $w \in [\lambda_\Theta \epsilon_k(x), \epsilon_k(x)]$}.
\]
This is an immediate implication of the inequality
$\epsilon_k(x)/x \leq \thetaub$ from
Proposition~\ref{prop:properties-of-xhat}--\ref{property:epsilon/x-ratio-bounds},
the identity $\lambda_\Theta^{C_F - 1} = \thetaub$,
and the bounds
$\lambda_\Theta \epsilon_k(x) \leq w \leq \epsilon_k(x)$.
We have
\[
    1
    \leq
    \thetaub \frac{ x }{ \epsilon_k(x) }
    \leq
    \lambda_\Theta^{C_F - 1} \frac{ x }{ w }
    \leq
    \left( \frac{ w }{ \epsilon_k(x) } \right)^{C_F - 1}
    \frac{ x }{ w },
\]
just as needed to complete the proof of
\eqref{eq:density-lower-bound-large-capacity}.

The integrand in \eqref{eq:positive-part-integral-bound}
therefore vanishes for
$w \in [\lambda_\Theta \epsilon_k(x), \epsilon_k(x)]$,
while, on the remaining interval, it is at most
$1 / (k \epsilon_k(x))$.
Therefore,
\[
    \int_0^{\epsilon_k(x)}
    \left(
        \frac{ 1 }{ k \epsilon_k(x) } - f(w)
    \right)_+
    \, dw
    \leq
    \int_0^{\lambda_\Theta \epsilon_k(x)}
    \frac{ 1 }{ k \epsilon_k(x) } \, dw
    =
    \frac{ \lambda_\Theta }{ k },
\]
concluding the proof.
\halmos
\endproof

\proof{Proof of Lemma~\ref{lemma:absolute-coefficient-estimate}.}
Fix $k \geq K$ and $x \in (\xhat_k, c)$,
so that Proposition~\ref{prop:properties-of-xhat}--\ref{property:epsilon/x-ratio-bounds}
implies
\(
    \epsilon_k(x) \leq \thetaub x.
\)
Using $\abs{ b } = b + 2(-b)_+$
with $b = f(w) - \frac{ 1 } { k \epsilon_k(x) }$,
we decompose the integral in \eqref{eq:absolute-coefficient-estimate} as
\[
    \int_0^{\epsilon_k(x)}
    \frac{1}{ (x - w)^{3/4} \epsilon_k(x - w)^{1/4} }
    \abs{ \frac{1}{k \epsilon_k(x)} - f(w) }
    \, dw
    =
    A + B + C,
\]
where
\begin{align*}
    A
    &=
    \int_0^{\epsilon_k(x)}
    \frac{1}{ (x - w)^{3/4} \epsilon_k(x - w)^{1/4} }
    f(w) \, dw, \\
    B
    &=
    -\frac{1}{k \epsilon_k(x)}
    \int_0^{\epsilon_k(x)}
    \frac{1}{ (x - w)^{3/4} \epsilon_k(x - w)^{1/4} }
    \, dw, \\
    C
    &=
    2 \int_0^{\epsilon_k(x)}
    \frac{1}{ (x - w)^{3/4} \epsilon_k(x - w)^{1/4} }
    \left(
        \frac{1}{k \epsilon_k(x)} - f(w)
    \right)_+
    \, dw,
\end{align*}
and we study each term in order.

First, we add and subtract
$x^{-3/4} \epsilon_k(x)^{-1/4}$
under the integral in $A$, so that
\[    
    A
    =
    \frac{ F(\epsilon_k(x)) }{ x^{3/4} \epsilon_k(x)^{1/4} } 
    +
    \int_0^{\epsilon_k(x)}
    \left(
        \frac{1}{ (x - w)^{3/4} \epsilon_k(x - w)^{1/4} }
        -
        \frac{1}{ x^{3/4} \epsilon_k(x)^{1/4} }
    \right)
    f(w) \, dw.
\]
Since $\epsilon_k(x) \leq \thetaub x$
and $0 \leq w \leq \epsilon_k(x)$, 
we also have that 
\(  \frac{1}{x - w}
    \leq
    \frac{1}{(1 - \thetaub)x}.
\)
This last bound, combined with the first inequality of
\eqref{eq:ratio-eps-x-to-x-minus-eps-3/4}
and the integral representation
\eqref{eq:consumption-function-integral-representation},
then implies
\begin{align*}
    \int_0^{\epsilon_k(x)}
    \left(
        \frac{1}{ (x - w)^{3/4} \epsilon_k(x - w)^{1/4} }
        -
        \frac{1}{ x^{3/4} \epsilon_k(x)^{1/4} }
    \right)
    f(w) \, dw 
    & \leq
    \frac{1}{ x^{3/4} \epsilon_k(x)^{1/4} }
    \int_0^{\epsilon_k(x)}
    \frac{w}{x - w} f(w) \, dw \\
    &\leq
    \frac{1}{ x^{3/4} \epsilon_k(x)^{1/4} }
    \frac{1}{(1 - \thetaub)x}
    \int_0^{\epsilon_k(x)} w f(w) \, dw  \\
    & =
    \frac{1}{1 - \thetaub}
    \frac{1}{ k x^{3/4} \epsilon_k(x)^{1/4} }.
\end{align*}
Thus,
\[
    A
    \leq
    \frac{  F(\epsilon_k(x)) }{  x^{3/4} \epsilon_k(x)^{1/4} }
    +
    \frac{1}{1 - \thetaub}
    \frac{1}{ k x^{3/4} \epsilon_k(x)^{1/4} }.
\]

Next, since $0 \leq w \leq \epsilon_k(x)$
and the map $x \mapsto \epsilon_k(x)$
is increasing, 
we have
\(
    (x - w)^{3/4} \epsilon_k(x - w)^{1/4}
    \leq
    x^{3/4} \epsilon_k(x)^{1/4}.
\)
This inequality, the rightmost bound
in \eqref{eq:ratio-eps-x-to-x-minus-eps-3/4}
and $\epsilon_k(x) / x \leq \thetaub$, 
then imply
\begin{equation}
\label{eq:H(x-w)-sandwich-bound}
    \frac{1}{ x^{3/4} \epsilon_k(x)^{1/4} }
    \leq
    \frac{1}{ (x - w)^{3/4} \epsilon_k(x - w)^{1/4} } 
    \leq
    \frac{x}{x - \epsilon_k(x)}
    \frac{1}{ x^{3/4} \epsilon_k(x)^{1/4} } 
    \leq
    \frac{1}{1 - \thetaub}
    \frac{1}{ x^{3/4} \epsilon_k(x)^{1/4} }.
\end{equation}
Multiplying the first inequality by $-1 / (k \epsilon_k(x))$
and integrating over $w \in [0, \epsilon_k(x)]$ gives
\[
    B
    =
    -\frac{1}{k \epsilon_k(x)}
    \int_0^{\epsilon_k(x)}
    \frac{1}{ (x - w)^{3/4} \epsilon_k(x - w)^{1/4} }
    \, dw
    \leq
    - \frac{1}{ k x^{3/4} \epsilon_k(x)^{1/4} }.
\]

Finally, we use the rightmost bound in
\eqref{eq:H(x-w)-sandwich-bound}
and combine it with
\eqref{eq:positive-part-integral-bound} to find
\begin{equation*}
    C
    \leq
    \frac{2}{ (1 - \thetaub) x^{3/4} \epsilon_k(x)^{1/4} }
    \int_0^{\epsilon_k(x)}
    \left(
        \frac{1}{k \epsilon_k(x)} - f(w)
    \right)_+
    \, dw 
    \leq
    \frac{ 2 \lambda_\Theta }{1 - \thetaub}
    \frac{1}{ k x^{3/4} \epsilon_k(x)^{1/4} }.
\end{equation*}

Combining the bounds for $A$, $B$, and $C$, we obtain
\[
    A + B + C
    \leq
    \frac{ F(\epsilon_k(x)) }{ x^{3/4} \epsilon_k(x)^{1/4} } 
    +
    \left(
        \frac{1}{1 - \thetaub}
        - 1
        + \frac{ 2 \lambda_\Theta }{1 - \thetaub}
    \right)
    \frac{1}{ k x^{3/4} \epsilon_k(x)^{1/4} } 
    =
    \left( F(\epsilon_k(x)) + \frac{C_\Theta}{k} \right)
    \frac{1}{ x^{3/4} \epsilon_k(x)^{1/4} },
\]
which proves \eqref{eq:absolute-coefficient-estimate}.
\halmos
\endproof

\proof{Proof of Lemma~\ref{lemma:vhat-upperbound-small-capacity}.}
Fix $k \geq K_1$.
The case $x = 0$ is immediate.
For $0 < x \leq \xhat_k$, we have $x \leq \wbar$,
and every item selected must have weight at most $x$.
The lower bound on the second ratio in
\eqref{eq:karamata-ratios-pointwise-bounds} and
Proposition~\ref{prop:properties-of-xhat}--\ref{property:consumption-bound-below-xhat}
therefore give
\[
    0
    \leq \vhat_k(x)
    \leq k F(x)
    \leq
    \frac{ k C_F }{ x } \int_0^x w f(w) \, dw
    \leq
    C_F \Theta
    =
    M_v,
\]
just as needed.
\halmos\endproof

\section{Proofs of the results in Section~\ref{sec:vhat-derivative-bounds}}
\label{apx:vhat-derivative-bounds-proofs}

\proof{Proof of Lemma~\ref{lem:vhat-continuity&differentiability}.}
For $k \leq K$, all assertions follow from $\vhat_k \equiv 0$.
For $k > K$, define
\[
    x_k
    =
    \max\{ x \in [0,c]: x \leq \epsilon_k(x) \}.
\]
Since
$\int_0^x w f(w) \, dw \leq x F(x) = o(x)$ as $x \to 0^+$,
the integral representation
\eqref{eq:consumption-function-integral-representation}
gives $\epsilon_k(x) > x$ for all sufficiently small $x > 0$.
Moreover, $\epsilon_k(c) < c$ by \eqref{eq:K}.
The concavity of $\epsilon_k$
(Proposition~\ref{prop:basic-properties-epsilon}--\ref{property:consumption-function-concavity})
therefore implies that $0 < x_k < c$, that $\epsilon_k(x_k) = x_k$,
and that
\(
    \epsilon_k(x) > x \text{ for } 0 < x < x_k,
    \text{ and }
    \epsilon_k(x) < x  \text{for } x_k < x \leq c.
\)
Consequently, $\hhat_k(x) = x$ below $x_k$ and
$\hhat_k(x) = \epsilon_k(x)$ above $x_k$.
By \eqref{eq:consumption-function-derivative},
$\hhat_k$ is continuously differentiable on
$(0,c) \setminus \{x_k\}$.
Set $\cD_k = \varnothing$ for $0 \leq k \leq K$, and define
\[
    \cD_k = \cD_{k-1} \cup \{x_k\}
    \qquad \text{for } k > K.
\]

We prove absolute continuity and continuous differentiability
simultaneously by induction.
Fix $k > K$ and suppose that $\vhat_{k-1}$ is absolutely continuous
on $[0,c]$ and continuously differentiable on
$(0,c) \setminus \cD_{k-1}$.
Since $\hhat_k$ and $\vhat_{k-1}$ are continuous and
$0 \leq \vhat_{k-1} \leq k-1$, dominated convergence applied to
\eqref{eq:value-function-recursion}, with dominating function
$k f(w)$ for the integral term, shows that $\vhat_k$ is continuous
on $[0,c]$.

For the following calculation, extend $\vhat_{k-1}$ by zero to negative
arguments, which preserves continuity and absolute continuity because
$\vhat_{k-1}(0) = 0$.
Assign its derivative the value zero wherever it is undefined.
For $x \in (0,c) \setminus \cD_{k-1}$ and $t \in (0,\wbar)$, let
\[
    G(x,t)
    =
    \int_0^t \{ 1 + \vhat_{k-1}(x-w) \} f(w) \, dw.
\]
This function is continuously differentiable, with
\[
    \partial_x G(x,t)
    =
    \int_0^t \vhat_{k-1}'(x-w) f(w) \, dw,
    \qquad
    \partial_t G(x,t)
    =
    \{ 1 + \vhat_{k-1}(x-t) \} f(t).
\]
To justify the first identity and its continuity, split the integral
at a sufficiently small $\delta > 0$.
Near $w = 0$, use uniform continuity of $\vhat_{k-1}'$
on a compact interval avoiding $\cD_{k-1}$.
On the remaining interval, $f$ is bounded, so continuity of translations
in $L^1$ and absolute continuity of the Lebesgue integral give joint
continuity as $x$ and $t$ vary.
Absolute continuity of $\vhat_{k-1}$, Fubini's theorem, and the fundamental
theorem of calculus then give the first identity.
The second identity and its continuity follow directly from continuity
of $\vhat_{k-1}$ and $f$.

By the induction hypothesis and the regularity of $G$ and $\hhat_k$,
the right-hand side of \eqref{eq:value-function-recursion} is
continuously differentiable on $(0,c) \setminus \cD_k$.
Consequently, so is $\vhat_k$.
Differentiating this expression using the chain and product rules
gives
\begin{align}
\label{eq:value-function-diff-explicit-formula}
    \vhat_k'(x)
    & =
    \{ 1 - F(\hhat_k(x)) \} \vhat_{k-1}'(x)
    + \int_0^{\hhat_k(x)}
        \vhat_{k-1}'(x-w) f(w) \, dw
    \notag \\
    & \quad
    + \hhat_k'(x) f(\hhat_k(x))
    \{ 1 + \vhat_{k-1}(x-\hhat_k(x)) - \vhat_{k-1}(x) \}.
\end{align}

It remains to establish absolute continuity of $\vhat_k$.
Since $0 \leq \vhat_{k-1} \leq k-1$ and $\hhat_k' \geq 0$,
the preceding identity implies, almost everywhere on $(0,c)$, that
\[
    \abs{ \vhat_k'(x) }
    \leq
    \abs{ \vhat_{k-1}'(x) }
    + \int_0^{\hhat_k(x)}
        \abs{ \vhat_{k-1}'(x-w) } f(w) \, dw
    + k \hhat_k'(x) f(\hhat_k(x)).
\]
Since $\hhat_k(x) = x$ for $x \leq x_k$ and
$\hhat_k(x) = \epsilon_k(x)$ for $x \geq x_k$, we have
\begin{equation*}
    \int_0^c \hhat_k'(x) f(\hhat_k(x)) \, dx
    =
    \int_0^{x_k} f(x) \, dx
    + \int_{x_k}^c \epsilon_k'(x) f(\epsilon_k(x)) \, dx 
    =
    F(x_k) + F(\epsilon_k(c)) - F(\epsilon_k(x_k)) 
    \leq 1.
\end{equation*}
Here, the second equality follows from the fundamental theorem of calculus applied to $x \mapsto F(\epsilon_k(x))$ on $[x_k,c]$, and the last
uses $\epsilon_k(x_k) = x_k$.
Furthermore, Tonelli's theorem and $\hhat_k(x) \leq x$ give
\[
    \int_0^c \int_0^{\hhat_k(x)}
        \abs{ \vhat_{k-1}'(x-w) } f(w) \, dw \, dx
    \leq
    \int_0^c \abs{ \vhat_{k-1}'(y) } \, dy.
\]
Consequently,
\[
    \int_0^c \abs{ \vhat_k'(x) } \, dx
    \leq
    2 \int_0^c \abs{ \vhat_{k-1}'(x) } \, dx + k
    < \infty.
\]
Since $\vhat_k$ is continuous on $[0,c]$ and continuously differentiable
outside the finite set $\cD_k$, the fundamental theorem of calculus
on the complementary intervals, followed by passage to their endpoints,
shows that $\vhat_k$ is absolutely continuous on $[0,c]$.
This completes the induction.

Finally, below $x_k$ we have $\hhat_k(x) = x$,
$\hhat_k'(x) = 1$, and $\vhat_{k-1}(0) = 0$.
Above $x_k$, we have $\hhat_k(x) = \epsilon_k(x)$ and
\[
    \hhat_k'(x) f(\hhat_k(x))
    =
    \frac{ 1 }{ k \epsilon_k(x) }.
\]
Substitution into \eqref{eq:value-function-diff-explicit-formula}
gives the two branches of \eqref{eq:vhat-derivative-recursion}.
\halmos
\endproof

% The proofs of Propositions~\ref{prop:vhat-intermediate-derivative-bounds-small-capacity} and~\ref{prop:vhat-derivative-bounds-small-capacity} follow next.
% We begin with the first one. 

\proof{Proof of Proposition~\ref{prop:vhat-intermediate-derivative-bounds-small-capacity}.}
We begin by showing that \eqref{eq:vhat-derivative-bounds-small-capacity} is an immediate consequence of \eqref{eq:vhat-derivative-bounds-small-capacity-induction}. 
In fact, we have that 
$\{1+F(x)\}^{k} \leq e^{k F(x)}$ 
and $ k F(x) \leq k F(\xhat_k) \leq M_v$ by \eqref{eq:vhat-upperbound-small-capacity}
for all $0 < x < \xhat_k$.
As such, one can just choose $C_v = (1 + M_v) e^{M_v}$ and \eqref{eq:vhat-derivative-bounds-small-capacity} holds.

Next, we fix \(x\in(0,\xhat_k)\) 
at which the derivative recursion \eqref{eq:vhat-derivative-recursion} holds
and establish that 
\begin{equation}
\label{eq:vhat-derivative-abs-recursion-bound}
\abs{ \vhat_k'(x) }
\leq 
\abs{ \vhat_{k-1}'(x) }
+\int_0^x \abs{ \vhat_{k-1}'(x-w) } f(w)\,dw
+(1 + M_v)f(x),
\qquad \text{for all } x \in (0, \xhat_k).
\end{equation}

For \(x<x_k\), the greedy branch of the derivative recursion \eqref{eq:vhat-derivative-recursion} gives
\[
\abs{ \vhat_k'(x) }
\leq 
 ( 1-F(x) ) \abs{ \vhat_{k-1}'(x) }
+ \int_0^x \abs{ \vhat_{k-1}'(x-w) } f(w)\,dw
+ f(x) \abs{ 1 - \vhat_{k-1}(x) },
\]
and \eqref{eq:vhat-derivative-abs-recursion-bound} holds after one drops $-F(x) \abs{ \vhat'_{k-1}(x)}$ and notices that the absolute value in the last summand can be bounded by $1+M_v$ from \eqref{eq:vhat-upperbound-small-capacity}.

For \(x\in(x_k,\xhat_k)\), the conservative branch of the derivative recursion \eqref{eq:vhat-derivative-recursion} gives
\[
\abs{ \vhat_k'(x) }
\leq 
( 1-F(\epsilon_k(x)) ) \abs{ \vhat_{k-1}'(x) }
+\int_0^{\epsilon_k(x)} \abs{ \vhat_{k-1}'(x-w) } f(w)\,dw 
+ \frac{\abs{ 1 + \vhat_{k-1}(x-\epsilon_k(x))- \vhat_{k-1}(x) }}{k\epsilon_k(x)}.
\]
Here, the last summand is bounded above by $( 1 + M_v ) f(x)$. 
To see this, note that, by \eqref{eq:vhat-upperbound-small-capacity},  the absolute value in the numerator is at most \(1+M_v\). 
Furthermore, the upper bound in the left estimate of
\eqref{eq:karamata-ratios-pointwise-bounds} at
$\epsilon = \epsilon_k(x) \leq \wbar$, the bounds
$\epsilon_k(x) \leq x \leq \xhat_k \leq \wbar$,
and the monotonicity of $w \mapsto w f(w)$ imply
\[
    \frac{ 1 }{ k \epsilon_k(x) }
    =
    \frac{ \int_0^{\epsilon_k(x)} w f(w) \, dw }{ x \epsilon_k(x) }
    \leq
    \frac{ \epsilon_k(x) f(\epsilon_k(x)) }{ x }
    \leq f(x).
\]
The bound \eqref{eq:vhat-derivative-abs-recursion-bound} follows after one drops $-F(\epsilon_k(x)) \abs{ \vhat'_{k-1}(x)}$ and enlarges the integration interval to \([0,x]\).

Next, we prove \eqref{eq:vhat-derivative-bounds-small-capacity-induction} by induction.
For \(k=K\), the result follows from \(\vhat_K \equiv 0\). 
Next, fix \(k>K\) and suppose that
\[
\abs{ \vhat_{k-1}'(x) }
\leq( 1 + M_v )(k-1)f(x)\{1+F(x)\}^{k-2},
\qquad \text{for all } 0 < x < \xhat_{k-1}.
\]
Since \(\xhat_k\leq\xhat_{k-1}\), the induction hypothesis applies to all of $(0, \xhat_k)$.
Substituting the induction hypothesis into \eqref{eq:vhat-derivative-abs-recursion-bound} yields
\begin{align}
\abs{ \vhat_k'(x) }
\leq (1 + M_v ) \bigg[&
(k-1)f(x)\{1+F(x)\}^{k-2}\notag \\
&+(k-1)\int_0^x
f(x-w)\{1+F(x-w)\}^{k-2}f(w)\,dw
+f(x)\bigg].
\label{eq:vhat-derivative-abs-w-inductionhyp}
\end{align}
To bound the integral in \eqref{eq:vhat-derivative-abs-w-inductionhyp}, 
we use the inequalities  $ (x - w) f(x - w) \leq x f(x)$, $F(x - w) \leq F(x)$, and $k > K$
to compute
\[
\int_0^x \frac{x-w}{x} f(x-w)\{1+F(x-w)\}^{k-2}f(w)\,dw 
\leq 
 f(x) F(x) \{1+F(x)\}^{k-2}.
\]
Similarly, from the bound $w f(w) \leq x f(x)$ and a change of variable, we obtain
\[
\int_0^x \frac{w}{x} f(x-w)\{1+F(x-w)\}^{k-2}f(w)\,dw 
\leq 
 f(x) \int_0^x  f( y )\{1+F( y )\}^{k-2} \,dy 
= f(x) \frac{\{ 1 + F(x) \}^{k-1} - 1}{k-1}.
\]
Adding up the last two integrals and their corresponding bounds, we obtain 
\begin{align*}
\int_0^x f(x-w)\{1+F(x-w)\}^{k-2}f(w)\,dw 
=& \int_0^x \left(\frac{x-w}{x} + \frac{w}{x}\right) f(x-w)\{1+F(x-w)\}^{k-2}f(w)\,dw\\
\leq &f(x) \left[ F(x) \{1+F(x)\}^{k-2} + \frac{\{ 1 + F(x) \}^{k-1} - 1}{k-1} \right].
\end{align*}
When we substitute this last bound into \eqref{eq:vhat-derivative-abs-w-inductionhyp} and rearrange, we find 
\[
\abs{ \vhat_k'(x) }
\leq(1+M_v)f(x)\Big[
(k-1)\{1+F(x)\}^{k-1}
+\{1+F(x)\}^{k-1}
\Big]
=(1+M_v)k f(x)\{1+F(x)\}^{k-1},
\]
just as needed to complete the induction proof. 
\Halmos\endproof

% We now prove Proposition~\ref{prop:vhat-derivative-bounds-small-capacity}.
% The argument combines \eqref{eq:vhat-derivative-bounds-small-capacity}
% with bounds valid on $(0, \xhat_k)$ to obtain the required derivative
% approximation in terms of $x$ and $\epsilon_{k+1}(x)$
% and the stated choice of $\Delta$.

\proof{Proof of Proposition~\ref{prop:vhat-derivative-bounds-small-capacity}.}
Fix $k \geq K$ and $x \in (0, \xhat_k)$ at which
$\vhat_k'(x)$ exists.
Let $C_v$ be the constant in
Proposition~\ref{prop:vhat-intermediate-derivative-bounds-small-capacity}.
Because both $\vhat'_k(x), 1/ \epsilon_{k+1}(x) \geq 0$, we bound each of them separately,
and we begin with $1 / \epsilon_{k+1}(x)$.
By Proposition~\ref{prop:properties-of-xhat}--\ref{property:epsilon/x-ratio-bounds}, for all $x < \xhat_k$, we have $\epsilon_k(x)/x \geq \thetalb$. Together with the space-time scaling equivalence \eqref{eq:consumption-function-space-time-scaling}, this implies  
\[
  \frac{x }{\epsilon_{k+1}(x)} 
  = 
  \frac{k+1}{k}\frac{\frac{k}{k+1}x}{\epsilon_{k}(\frac{k}{k+1}x)} 
  \leq 
  2 \cdot \frac{1}{\thetalb} = 2 \Theta.
\]
Consequently, since $\Delta \geq (2 \Theta)^{3/4}$, we also have
\begin{equation}
\label{eq:small-capacity-reciprocal-bound}
        \frac{ 1 }{ \epsilon_{k+1}(x) }
         =
        \left( \frac{ x }{ \epsilon_{k+1}(x) } \right)^{3/4}
        \frac{ 1 }{ x^{3/4} \epsilon_{k+1}(x)^{1/4} }
         \leq
        \frac{ (2 \Theta)^{3/4} }{ x^{3/4} \epsilon_{k+1}(x)^{1/4} }
        \leq
        \frac{ \Delta }{ x^{3/4} \epsilon_{k+1}(x)^{1/4} }.
\end{equation}

Next, we recall from
\eqref{eq:vhat-derivative-bounds-small-capacity} that
$\abs{ \vhat_k'(x) } \leq C_v k f(x)$.
We seek to bound $C_v k f(x)$ on the same interval.
If we use Proposition~\ref{prop:properties-of-xhat}--\ref{property:consumption-bound-below-xhat} one more time and recall 
\eqref{eq:epsilon-k-leq-wbar}, we have that 
\(
    x \leq \epsilon_k(\Theta x) \leq \wbar.
\)
Moreover, the monotonicity of the consumption function and
$\Theta \geq 1$ imply
\(
    \epsilon_{k+1}(x)
    \leq \epsilon_k(x)
    \leq \epsilon_k(\Theta x).
\)
It then follows that
\[
    x^{3/4} \epsilon_{k+1}(x)^{1/4}
    \leq
    \epsilon_k(\Theta x)
    \qquad \text{and} \qquad
    k f(x) x^{3/4} \epsilon_{k+1}(x)^{1/4}
    \leq
    k f(x) \epsilon_k(\Theta x).
\]
If we now multiply and divide the right-hand side of the last inequality by $x$, 
use the monotonicity of $w \mapsto w f(w)$,
and use the lower bound of the left estimate of \eqref{eq:karamata-ratios-pointwise-bounds} at
$\epsilon = \epsilon_k(\Theta x) \leq \wbar$,
we find 
\begin{equation*}
    k f(x) x^{3/4} \epsilon_{k+1}(x)^{1/4}
    \leq
    \frac{ k \epsilon_k(\Theta x)^2
        f( \epsilon_k(\Theta x) ) }{ x } 
    =
    \Theta \cdot \frac{ k \epsilon_k(\Theta x)^2
        f( \epsilon_k(\Theta x) ) }{ \Theta x } 
    \leq
    \Theta C_F
    \leq \frac{\Delta}{C_v}.
\end{equation*}
Finally, from \eqref{eq:vhat-derivative-bounds-small-capacity},
the bound above, and
\eqref{eq:small-capacity-reciprocal-bound}, we find
\[
    \abs{ \vhat_k'(x) - \frac{ 1 }{ \epsilon_{k+1}(x) } }
    \leq
    \max\left\{ C_v k f(x), \frac{ 1 }{ \epsilon_{k+1}(x) } \right\}
    \leq
    \frac{ \Delta }{ x^{3/4} \epsilon_{k+1}(x)^{1/4} },
\]
concluding the proof of the proposition.
\halmos\endproof

\proof{Proof of Proposition~\ref{prop:Pk-absolute-upper-bound}.}
Fix $k > K$ and $x \in (\xhat_k, c) \setminus \cD_k$.
Since $x \in [\xhat_k, c)$, the definition \eqref{eq:xk-definition}
implies that $\epsilon_k(x) \leq x$
and the integrand in \eqref{eq:Pk} is defined almost everywhere on $(0, \epsilon_k(x))$.
Furthermore, we have that $\cD_{k-1} \subseteq \cD_k$ and that $\cD_{k-1}$ is finite, so \eqref{eq:vhat-derivative-approximation-error-assumption}
applies at the chosen $x$ and almost everywhere on $(0,c)$.
If we now recall \eqref{eq:Pk},
apply the triangle inequality, and use the assumption \eqref{eq:vhat-derivative-approximation-error-assumption},
we find 
\[
\abs{ P_k(x) }
    \leq
    \{ 1 - F( \epsilon_k(x) ) \}
    \frac{ M }{ x^{3/4} \epsilon_k(x)^{1/4} }
    +
    \int_0^{\epsilon_k(x)}
    \frac{ M }{ (x - w)^{3/4} \epsilon_k(x - w)^{1/4} }
    \abs{ \frac{ 1 }{ k \epsilon_k(x) } - f(w) }
    \, dw.
\]
Lemma~\ref{lemma:absolute-coefficient-estimate} applied to the last summand then gives 
\[
    \abs{ P_k(x) }
    \leq 
    \left(
        1 - F( \epsilon_k(x) )
        + F( \epsilon_k(x) )
        + \frac{ C_\Theta }{ k }
    \right)
    \frac{ M }{ x^{3/4} \epsilon_k(x)^{1/4} }
    =
    \left( 1 + \frac{ C_\Theta }{ k } \right)
    \frac{ M }{ x^{3/4} \epsilon_k(x)^{1/4} },
\]
just as needed.
\halmos
\endproof

\proof{Proof of Proposition~\ref{prop:Rk-absolute-upper-bound}.}
Fix $k \geq K$ and $x \in (\xhat_k, c)$.
We recall the decomposition \eqref{eq:Rk} and estimate each difference separately. 
First, we estimate the time-shift difference.
The space-time scaling equivalence
\eqref{eq:consumption-function-space-time-scaling} gives us that
\[
    \epsilon_{k+1}(x)
    =
    \epsilon_k\left( x - \frac{ x }{ k + 1 } \right),
\]
so when we apply \eqref{eq:difference-of-consumption-function-reciprocals-linear-approximation}
with $w = x / (k + 1)$ and 
weaken the lower bound by using $1 / (k + 1) \leq 1 / k$ in the $ C_\epsilon$ term, 
we find
\begin{equation}
\label{eq:consumption-function-reciprocals-time-difference}
    - \frac{ \epsilon_k'(x) }{ \epsilon_k(x)^2 }
        \frac{ x }{ k + 1 }
    - \frac{  C_\epsilon }{ k x^{3/4} \epsilon_k(x)^{1/4} }
    \leq
    \frac{ 1 }{ \epsilon_k(x) }
    - \frac{ 1 }{ \epsilon_{k+1}(x) } 
    \leq
    - \frac{ \epsilon_k'(x) }{ \epsilon_k(x)^2 }
        \frac{ x }{ k + 1 }.
\end{equation}

Second, we recall from \eqref{eq:difference-of-consumption-function-reciprocals-dFintegral} that the difference of the consumption function reciprocals integrated against the weight distribution
satisfies the bounds 
\begin{equation}
\label{eq:consumption-function-reciprocals-space-difference-dFintegral}
    \frac{ \epsilon_k'(x) }{ \epsilon_k(x)^2 }
        \frac{x}{k}
    \leq
    \int_0^{\epsilon_k(x)}
    \left(
        \frac{ 1 }{ \epsilon_k(x - w) }
        -
        \frac{ 1 }{ \epsilon_k(x) }
    \right)
    f(w) \, dw 
     \leq
    \frac{ \epsilon_k'(x) }{ \epsilon_k(x)^2 }
        \frac{x}{k}
    +
    \frac{  C_\epsilon }{ k x^{3/4} \epsilon_k(x)^{1/4} }.
\end{equation}

Third, we estimate the integral of the difference of the reciprocals of the consumption function by integrating the upper and lower bounds in
\eqref{eq:difference-of-consumption-function-reciprocals-linear-approximation} to obtain 
\begin{equation}
\label{eq:consumption-function-reciprocals-space-difference-integral}
    - \frac{ \epsilon_k'(x) }{ 2k \epsilon_k(x) }
    - \frac{  C_\epsilon }{ 2k x^{3/4} \epsilon_k(x)^{1/4} }
      \frac{ \epsilon_k(x) }{ x }
    \leq 
    - \frac{ 1 }{ k \epsilon_k(x) }
    \int_0^{\epsilon_k(x)}
    \left(
        \frac{ 1 }{ \epsilon_k(x - w) }
        -
        \frac{ 1 }{ \epsilon_k(x) }
    \right)
    \, dw
    \leq
    - \frac{ \epsilon_k'(x) }{ 2k \epsilon_k(x) }.
\end{equation}

If we now factor $\frac{ \epsilon_k'(x) }{ \epsilon_k(x)^2 }$ in the bounds \eqref{eq:consumption-function-reciprocals-time-difference}, \eqref{eq:consumption-function-reciprocals-space-difference-dFintegral}, and \eqref{eq:consumption-function-reciprocals-space-difference-integral} and sum, 
we find 
\begin{equation*}
    - \left( 1 + \frac{ \epsilon_k(x) }{ 2x } \right)
      \frac{  C_\epsilon }{ k x^{3/4} \epsilon_k(x)^{1/4} } 
    \leq
    R_k(x)
    - \frac{ \epsilon_k'(x) }{ \epsilon_k(x)^2 }
    \left(
        \frac{ x }{ k(k + 1) }
        -
        \frac{ \epsilon_k(x) }{ 2k }
    \right) 
    \leq
    \frac{  C_\epsilon }{ k x^{3/4} \epsilon_k(x)^{1/4} }.
\end{equation*}
If we now recall that $\epsilon_k(x) \leq x$ and that $\epsilon_k'(x) \geq 0$, 
take absolute values,
use the triangle inequality, 
and rearrange,
we find
\[
    \abs{ R_k(x) }
    \leq 
    \frac{ \epsilon_k'(x) }{ k \epsilon_k(x)^2 }
        \abs{
            \frac{ x }{ k + 1 }
            -
            \frac{ \epsilon_k(x) }{ 2 }
        }
    +
    \frac{ 3 C_\epsilon }{ 2k x^{3/4} \epsilon_k(x)^{1/4} },
\]
and we just need to study the first summand on the right-hand side. 
From the upper bound of the right estimate in 
\eqref{eq:karamata-ratios-pointwise-bounds} at
$\epsilon = \epsilon_{k+1}(x) \leq \wbar$, we find $0 \leq x / (k + 1) \leq \epsilon_{k+1} (x) F(\epsilon_{k+1}(x))$, and then $0 \leq x / (k + 1) \leq \epsilon_{k} (x)$ since $\epsilon_k(x)$ is decreasing in $k$, and $F$ is bounded by one.
This shows that
\[
    -  \frac{ \epsilon_k(x) }{ 2 }
    \leq 
        \frac{ x }{ k + 1 }
        - \frac{ \epsilon_k(x) }{ 2 }
    \leq 
    \frac{ \epsilon_k(x) }{ 2 }.
\]
Furthermore, from 
Proposition~\ref{prop:basic-properties-epsilon}--\ref{property:consumption-function-concavity} 
we know that the map $x \mapsto \epsilon_k(x)/x$ is decreasing so 
differentiating gives 
\[
\left(\frac{\epsilon_k(x)}{x}\right)'
= 
\frac{1}{x^2}
\left( 
    x \epsilon_k'(x) - \epsilon_k(x)
\right)
\leq 
0,
\qquad \text{hence} \qquad
\epsilon_k'(x) 
\leq \frac{\epsilon_k(x)}{x}.
\]
When we use these last two observations in the upper bound for $\abs{ R_k(x) }$,
and note that $\epsilon_k(x) \leq x$ implies $x^{3/4} \epsilon_k(x)^{1/4} \leq x$,
we find 
\[
    \abs{ R_k(x) }
    \leq
    \frac{ 1 }{ 2kx } 
    +    
    \frac{ 3 C_\epsilon }{ 2k x^{3/4} \epsilon_k(x)^{1/4} }
    \leq 
    \frac{ 1 }{ 2 k x^{3/4} \epsilon_k(x)^{1/4}} 
    +    
    \frac{ 3 C_\epsilon }{ 2k x^{3/4} \epsilon_k(x)^{1/4} }
    =
    \left(\frac{ 3 C_\epsilon + 1 }{ 2 }\right)
    \frac{ 1 }{ k x^{3/4} \epsilon_k(x)^{1/4} },
\]
just as needed.
\halmos
\endproof

% \subsection{Final induction argument}
% \label{subsec:vhat-derivative-bound-final-induction-step}

% We are now ready to prove Proposition~\ref{prop:vhat-derivative-bounds} by induction on
% $k \geq K$.
% The proof combines the small-capacity derivative bound of
% Proposition~\ref{prop:vhat-derivative-bounds-small-capacity} with the propagated-error and
% local-error bounds of Propositions~\ref{prop:Pk-absolute-upper-bound} and
% \ref{prop:Rk-absolute-upper-bound}.

\proof{Proof of Proposition~\ref{prop:vhat-derivative-bounds}.}
We prove the desired bound by induction on $k \geq K$.
For the base case $k = K$, the policy definition gives
$\vhat_K \equiv 0$, and hence $\vhat_K' \equiv 0$ on $(0, c)$
and 
\(
    \abs{ \vhat_K'(x) - \frac{ 1 }{ \epsilon_{K+1}(x) } }
     =
    \frac{ 1 }{ \epsilon_{K+1}(x) } 
\).
If we now multiply and divide this last right-hand side by $x^{3/4}$, rearrange,
and use the increasing property of the map  $x \mapsto x / \epsilon_{K+1}(x)$ (Proposition~\ref{prop:basic-properties-epsilon}--\ref{property:consumption-function-concavity}),  
we obtain for every $x \in (0, c)$, 
\begin{align*}
    \abs{ \vhat_K'(x) - \frac{ 1 }{ \epsilon_{K+1}(x) } }
    =
    \frac{ (  x / \epsilon_{K+1}(x) )^{3/4} }{ x^{3/4} \epsilon_{K+1}(x)^{1/4} }
    \leq
    \frac{ ( c / \epsilon_{K+1}(c) )^{3/4} }{ x^{3/4} \epsilon_{K+1}(x)^{1/4} } 
     \leq
    \frac{ \Delta }{ x^{3/4} \epsilon_{K+1}(x)^{1/4} }.
\end{align*}
This establishes the base case throughout $(0, c)$.

For the induction step, fix $k > K$ and assume that
\[
    \abs{ \vhat_{k-1}'(x) - \frac{ 1 }{ \epsilon_k(x) } }
    \leq
    \frac{ \Delta }{ x^{3/4} \epsilon_k(x)^{1/4} }
    \qquad
    \text{for every } x \in (0,c) \setminus \cD_{k-1}.
\]
Fix any $x \in (0,c) \setminus \cD_k$.
If $x < \xhat_k$, then
Proposition~\ref{prop:vhat-derivative-bounds-small-capacity}
and the choice of $\Delta$ immediately give
\[
    \abs{ \vhat_k'(x) - \frac{ 1 }{ \epsilon_{k+1}(x) } }
    \leq
    \frac{ \Delta }{ x^{3/4} \epsilon_{k+1}(x)^{1/4} }.
\]
%If $x = \xhat_k$, the same bound follows by taking a limit from below, since $\vhat_k'$ is continuous at $x$ by Lemma~\ref{lem:vhat-continuity&differentiability} and $\epsilon_{k+1}$ is continuous.
If $x > \xhat_k$, the decomposition \eqref{eq:vhat-derivative-decomposition-large-capacity},
the triangle inequality, the bound in
Proposition~\ref{prop:Pk-absolute-upper-bound} with $M = \Delta$ and the inequality in Proposition~\ref{prop:Rk-absolute-upper-bound},
and a rearrangement yield
\begin{equation}
\label{eq:vhat-derivative-induction-intermediate-bound}
    \abs{ \vhat_k'(x) - \frac{ 1 }{ \epsilon_{k+1}(x) } }
    \leq
    \abs{ P_k(x) } + \abs{ R_k(x) } 
    \leq
    \frac{ \Delta }{ x^{3/4} \epsilon_k(x)^{1/4} }
    +
    \frac{ \Delta C_\Theta + (3 C_\epsilon + 1)/2 }{
        k x^{3/4} \epsilon_k(x)^{1/4}
    }.
\end{equation}
The right-hand side of \eqref{eq:vhat-derivative-induction-intermediate-bound}
depends on $\epsilon_{k}(x)$, so we use the space-time scaling identity
\eqref{eq:consumption-function-space-time-scaling}
to compare it with the desired bound involving $\epsilon_{k+1}(x)$.
We have 
\(
    \epsilon_{k+1}(x)
    =
    \epsilon_k\left( x - \frac{ x }{ k + 1 } \right),
\)
so the lower bound in
\eqref{eq:taylor-first-order-eps-a}
with $a = 1/4$ and $w = x / (k + 1)$ implies
\[
    \frac{ 1 }{ \epsilon_k(x)^{1/4} }
    \leq
    \frac{ 1 }{ \epsilon_{k+1}(x)^{1/4} }
    -
    \frac{ x }{ 4(k + 1) }
    \frac{ \epsilon_k'(x) }{ \epsilon_k(x)^{5/4} }.
\]
We now upper bound the right-hand side by using the lower bound in
Proposition~\ref{prop:joint-properties-f-epsilon}--\ref{property:epsilon'/epsilon^2-delta-bounds} with $a = 1/4$.
We have
\[
    \frac{ \epsilon_k'(x) }{ \epsilon_k(x)^{5/4} }
    \geq
    \frac{1}{C_F} \frac{ 1 }{x \epsilon_k(x)^{1/4} },
\]
so we find 
\[
    \frac{ 1 }{ \epsilon_k(x)^{1/4} }
    \leq
    \frac{ 1 }{ \epsilon_{k+1}(x)^{1/4} }
    -
    \frac{ x }{ 4(k + 1) }
    \frac{ \epsilon_k'(x) }{ \epsilon_k(x)^{5/4} }
    \leq 
    \frac{ 1 }{ \epsilon_{k+1}(x)^{1/4} }
    -
    \frac{1}{4 C_F}\frac{ 1 }{(k + 1) \epsilon_k(x)^{1/4} }.
\]
If we now divide through by $x^{3/4}$ and use the inequality $k + 1 \leq 2k$,
we obtain
\[
    \frac{ 1 }{ x^{3/4} \epsilon_k(x)^{1/4} }
    \leq
    \frac{ 1 }{ x^{3/4} \epsilon_{k+1}(x)^{1/4} }
    -
    \frac{1}{8 C_F}\frac{ 1 }{k  x^{3/4} \epsilon_k(x)^{1/4} }.
\]
When we substitute this last bound into the first summand on the right-hand side of \eqref{eq:vhat-derivative-induction-intermediate-bound}, we find 
\begin{align*}
    \abs{ \vhat_k'(x) - \frac{ 1 }{ \epsilon_{k+1}(x) } }
    & \leq
    \frac{ \Delta }{ x^{3/4} \epsilon_{k+1}(x)^{1/4} }
    +
    \frac{
        \Delta C_\Theta + (3 C_\epsilon + 1)/2 - \Delta / (8 C_F)
    }{ k x^{3/4} \epsilon_k(x)^{1/4} }
     \\
    & =
    \frac{ \Delta }{ x^{3/4} \epsilon_{k+1}(x)^{1/4} }
    +
    \frac{ (3 C_\epsilon + 1)/2 - \Delta d_\Theta }{
        k x^{3/4} \epsilon_k(x)^{1/4}
    }
    \\
    & \leq
    \frac{ \Delta }{ x^{3/4} \epsilon_{k+1}(x)^{1/4} }.
\end{align*}
Here, the equality uses the definition 
$d_\Theta = 1 / (8 C_F) - C_\Theta$,
and the last inequality follows from $\Delta d_\Theta \geq (3 C_\epsilon + 1)/2$.
This then establishes the desired bound for every
$x \in (0, c) \setminus \cD_k$, 
completing the induction argument. 

To prove \eqref{eq:intermediate-integrand-overall-bound}, we recall from 
Lemma~\ref{lem:vhat-continuity&differentiability} that
$\vhat_k$ is absolutely continuous on $[0,c]$.
Hence, the fundamental theorem of calculus gives
\begin{align*}
    \vhat_k(x) - \vhat_k(x-w)
    - \frac{ w }{ \epsilon_{k+1}(x) }
    & =
    \int_{x-w}^x
        \left(
            \vhat_k'(y) - \frac{ 1 }{ \epsilon_{k+1}(x) }
        \right)
    \, dy    \\
    & =
    \int_{x-w}^x
        \left(
            \vhat_k'(y) - \frac{ 1 }{ \epsilon_{k+1}(y) }
        \right)
    \, dy
    +
    \int_{x-w}^x
        \left(
            \frac{ 1 }{ \epsilon_{k+1}(y) }
            - \frac{ 1 }{ \epsilon_{k+1}(x) }
        \right)
    \, dy.
\end{align*}
We consider the two terms on this last right-hand side separately.
For the first term, since $k \geq K$,
Proposition~\ref{prop:vhat-derivative-bounds} applies for almost every $y \in [x-w, x]$ at index $k$. 
Thus, by the triangle inequality,
\[
    \abs{
        \int_{x-w}^x
            \left(
                \vhat_k'(y) - \frac{ 1 }{ \epsilon_{k+1}(y) }
            \right)
        \, dy
    }
    \leq
    \int_{x-w}^x
        \frac{ \Delta }{ y^{3/4} \epsilon_{k+1}(y)^{1/4} }
    \, dy.
\]
For the second term, since $w \leq \hhat_{k+1}(x)$,
\eqref{eq:1/epsilon-integral-bounds} gives
\[
    0
    \leq
    \int_{x-w}^x
        \left(
            \frac{ 1 }{ \epsilon_{k+1}(y) }
            - \frac{ 1 }{ \epsilon_{k+1}(x) }
        \right)
    \, dy
    \leq
    \frac{ w^2 }{ (x-w) \epsilon_{k+1}(x) }.
\]
Combining these two bounds and applying the triangle inequality
proves the desired result.
\halmos
\endproof

% \subsection{Proofs of Corollaries of Proposition~\ref{prop:vhat-derivative-bounds}} 
% \label{sec:corollaries-of-non-asymptotic-derivative-bounds}

% We now prove the two corollaries of
% Proposition~\ref{prop:vhat-derivative-bounds}.
% We begin with Corollary~\ref{cor:uniform-bound-on-difference-in-value-function},
% which establishes a uniform absolute bound on the increments of the value function
% over each selection interval.

\proof{Proof of Corollary~\ref{cor:uniform-bound-on-difference-in-value-function}.}
Let $M_v$ be the constant in
Lemma~\ref{lemma:vhat-upperbound-small-capacity}.
Fix $k \geq 1$, $x \in [0,c]$, and
$w \in [0, \hhat_{k+1}(x)]$.
First, if $k < K$ or $x = 0$, then $\hhat_{k+1}(x) = 0$.
Therefore, the only possible value of $w$ is zero,
and the result follows trivially.

Suppose that $k \geq K$ and $x < \xhat_{k+1}$.
By Proposition~\ref{prop:properties-of-xhat}--\ref{property:fixed-points-decreasing},
we have $\xhat_{k+1} \leq \xhat_k$, so the bound \eqref{eq:vhat-upperbound-small-capacity}
and the inequalities $ 0 \leq \vhat_k(x-w) \leq \vhat_k(x)$ imply
\[
    \abs{ \vhat_k(x) - \vhat_k(x-w) }
    \leq \vhat_k(x) 
    \leq M_v.
\]

Suppose instead that $k \geq K$ and $x \geq \xhat_{k+1}$, so that $\hhat_{k+1}(x) = \epsilon_{k+1}(x) \leq x$.
By Proposition~\ref{prop:properties-of-xhat}--\ref{property:epsilon/x-ratio-bounds}
and the choice $\thetaub < 1/2$, 
$\hhat_{k+1}(x) = \epsilon_{k+1}(x) \leq \thetaub x < x/2$.
The estimate \eqref{eq:intermediate-integrand-overall-bound},
applied at index $k$ and extended to the endpoints
by continuity, gives
\[
    \abs{ \vhat_k(x) - \vhat_k(x-w) }
    \leq
    \int_{x-w}^x
        \frac{ \Delta }{ y^{3/4} \epsilon_{k+1}(y)^{1/4} }
    \, dy
    +
    \frac{ w^2 }{ (x-w) \epsilon_{k+1}(x) }
    +
    \frac{ w }{ \epsilon_{k+1}(x) }.
\]
Since $x-w \leq y$, we can use the monotonicity of
$\epsilon_{k+1}$ to bound the integrand by its value at $x-w$.
Moreover, since $w \leq \epsilon_{k+1}(x)$, we have
\(
    \frac{ w }{ \epsilon_{k+1}(x) } \leq 1,
\)
and
\(
    \frac{ w^2 }{ (x-w) \epsilon_{k+1}(x) }
    \leq
    \frac{ \epsilon_{k+1}(x) }{ x - \epsilon_{k+1}(x) }.
\)
As a result,
\begin{equation*}
    \abs{ \vhat_k(x) - \vhat_k(x-w) }
    \leq
    \int_{x-w}^x
        \frac{ \Delta }{ (x-w)^{3/4} \epsilon_{k+1}(x-w)^{1/4} }
    \, dy
    +
    \frac{ \epsilon_{k+1}(x) }{ x - \epsilon_{k+1}(x) }
    + 1.
\end{equation*}
Applying \eqref{eq:ratio-eps-x-to-x-minus-eps-3/4}
to the integrand of the first term, integrating for $y \in [x-w, x]$, using $w \leq \epsilon_{k+1}(x)$, 
and rearranging, we obtain
\begin{align*}
    \abs{ \vhat_k(x) - \vhat_k(x-w) }
    & \leq
    \int_{x-w}^x
        \frac{ \Delta x }{ (x - \epsilon_{k+1}(x)) x^{3/4} \epsilon_{k+1}(x)^{1/4} }
    \, dy
    +
    \frac{ \epsilon_{k+1}(x) }{ x - \epsilon_{k+1}(x) }
    + 1 \\
    & \leq
    \frac{ \Delta x \epsilon_{k+1}(x) }{ ( x - \epsilon_{k+1}(x) ) x^{3/4} \epsilon_{k+1}(x)^{1/4} }
    +
    \frac{ \epsilon_{k+1}(x) }{ x - \epsilon_{k+1}(x) }
    + 1 \\
    & = 
    \frac{ \Delta }{ 1 - \epsilon_{k+1}(x) / x }
    \left( \frac{ \epsilon_{k+1}(x) }{ x } \right)^{3/4}
    +
    \frac{ \epsilon_{k+1}(x)/x }{ 1 - \epsilon_{k+1}(x)/x }
    + 1.
\end{align*}
If we now recall that $\epsilon_{k+1}(x)/x \leq \thetaub < 1/2$,
we conclude that
\[
    \abs{ \vhat_k(x) - \vhat_k(x-w) } 
    \leq
    2^{1/4} \Delta + 2.
\]
Choosing
\(
    \Delta_v
    =
    \max\left\{ M_v, 2^{1/4} \Delta + 2 \right\}
\)
therefore proves the desired result.
\halmos
\endproof

% We next prove Corollary~\ref{cor:error-bound-of-linear-integrand-approximation-k},
% which bounds the error obtained when replacing
% $\left( 1 + \vhat_k(x-w) - \vhat_k(x) \right)^2$
% by its squared linear approximation
% $\left( 1 - \frac{ w }{ \hhat_{k+1}(x) } \right)^2$
% and integrating against the weight distribution over the selection interval.

\proof{Proof of Corollary~\ref{cor:error-bound-of-linear-integrand-approximation-k}.}
Fix $k \geq K$ and $x \in (0,c]$.
For each $w \in [0, \hhat_{k+1}(x)]$, write
\[
    a = 1 + \vhat_k(x-w) - \vhat_k(x)
    \qquad \text{and} \qquad
    b = 1 - \frac{ w }{ \hhat_{k+1}(x) }.
\]
Using $a^2 - b^2 = (a+b)(a-b)$, the integrand on the left-hand side
of \eqref{eq:error-bound-of-linear-integrand-approximation-square-k}
can be written as $(a+b)(a-b)f(w)$.
Since $0 \leq w \leq \hhat_{k+1}(x)$, we have $0 \leq b \leq 1$.
Moreover, Corollary~\ref{cor:uniform-bound-on-difference-in-value-function},
applied at index $k$, gives $\abs{ a } \leq 1 + \Delta_v$.
Thus, $\abs{ a+b } \leq \Delta_v + 2$.
Together with the triangle inequality, this bounds the left-hand side
of \eqref{eq:error-bound-of-linear-integrand-approximation-square-k} by
\begin{equation}
\label{eq:error-bound-of-linear-integrand-approximation-square-k-intermediate}
    \int_0^{\hhat_{k+1}(x)}
        \abs{ (a+b)(a-b) } f(w) \, dw
    \leq
    (\Delta_v + 2)
    \int_0^{\hhat_{k+1}(x)}
        \abs{
            \vhat_k(x) - \vhat_k(x-w)
            - \frac{ w }{ \hhat_{k+1}(x) }
        }
        f(w) \, dw.
\end{equation}

Suppose first that $x < \xhat_{k+1}$.
Corollary~\ref{cor:uniform-bound-on-difference-in-value-function}
and $w/\hhat_{k+1}(x) \leq 1$ bound the absolute value
in this last integral by $\Delta_v + 1$.
Thus, the right-hand side of
\eqref{eq:error-bound-of-linear-integrand-approximation-square-k-intermediate}
is at most
\[
    (\Delta_v + 2)(\Delta_v + 1) F(\hhat_{k+1}(x))
    \leq
    (\Delta_v + 2)(\Delta_v + 1) F(x)
    \leq
    \frac{ (\Delta_v + 2)(\Delta_v + 1) M_v }{ k+1 },
\]
where the last inequality follows from
\eqref{eq:vhat-upperbound-small-capacity}, with $M_v$ as in
Lemma~\ref{lemma:vhat-upperbound-small-capacity}.

Suppose instead that $x \geq \xhat_{k+1}$.
By Proposition~\ref{prop:properties-of-xhat}--\ref{property:epsilon/x-ratio-bounds}
and the choice $\thetaub < 1/2$, 
we have
$\hhat_{k+1}(x) = \epsilon_{k+1}(x) \leq \thetaub x < x/2$.
For $w \in (0, \epsilon_{k+1}(x))$ and $y \in [x-w, x]$,
monotonicity of $\epsilon_{k+1}$ and
\eqref{eq:ratio-eps-x-to-x-minus-eps-3/4} give
\[
    \frac{ 1 }{ y^{3/4} \epsilon_{k+1}(y)^{1/4} }
    \leq
    \frac{ x }{ x - \epsilon_{k+1}(x) }
    \frac{ 1 }{ x^{3/4} \epsilon_{k+1}(x)^{1/4} }
    \leq
    \frac{ 2 }{ x^{3/4} \epsilon_{k+1}(x)^{1/4} }.
\]
Also, since $w \leq \epsilon_{k+1}(x) \leq x/2$,
\[
    \frac{ w^2 }{ (x-w) \epsilon_{k+1}(x) }
    \leq
    \frac{ w }{ x-w }
    \leq
    \frac{ w }{ x - \epsilon_{k+1}(x) }
    \leq
    \frac{ 2w }{ x }.
\]
Hence, by \eqref{eq:intermediate-integrand-overall-bound},
\[
    \abs{
        \vhat_k(x) - \vhat_k(x-w)
        - \frac{ w }{ \epsilon_{k+1}(x) }
    }
    \leq
    \int_{x-w}^x
        \frac{ \Delta }{ y^{3/4} \epsilon_{k+1}(y)^{1/4} }
    \, dy
    +
    \frac{ w^2 }{ (x-w) \epsilon_{k+1}(x) }
    \leq
    2w \left\{
        \frac{ \Delta }{ x^{3/4} \epsilon_{k+1}(x)^{1/4} }
        + \frac{ 1 }{ x }
    \right\}.
\]
Integrating against $f(w) \, dw$ over $[0, \epsilon_{k+1}(x)]$
and using \eqref{eq:consumption-function-integral-representation},
we see that the right-hand side of
\eqref{eq:error-bound-of-linear-integrand-approximation-square-k-intermediate}
is bounded by
\[
    2(\Delta_v + 2)
    \left\{
        \frac{ \Delta x^{1/4} }{ (k+1) \epsilon_{k+1}(x)^{1/4} }
        + \frac{ 1 }{ k+1 }
    \right\}.
\]
Finally, the upper bound on the second ratio in
\eqref{eq:karamata-ratios-pointwise-bounds} at
$\epsilon = \epsilon_{k+1}(x) \leq \wbar$ gives
\[
    \frac{ x^{1/4} }{ (k+1) \epsilon_{k+1}(x)^{1/4} }
    =
    \frac{ F^{1/4}(\epsilon_{k+1}(x)) }{ (k+1)^{3/4} }
    \left(
        \frac{ x }{ (k+1) \epsilon_{k+1}(x) F(\epsilon_{k+1}(x)) }
    \right)^{1/4}
    \leq
    \frac{ F^{1/4}(\epsilon_{k+1}(x)) }{ (k+1)^{3/4} }.
\]
Since $\hhat_{k+1}(x) = \epsilon_{k+1}(x)$ in this case, choosing
\(
    \Delta_2
    =
    (\Delta_v + 2)
    \max\left\{
        (\Delta_v + 1) M_v, 2 \Delta, 2
    \right\}
\)
establishes the desired bound in both cases.
\halmos
\endproof

\section{Proofs of the results in Section~\ref{sec:variance-asymptotics}}
\label{apx:variance-asymptotics-proofs}

\proof{Proof of Lemma~\ref{lem:series-upper-bound}.}
Fix $t\in\{0,\ldots,T-1\}$ and set 
\(
p_j:=F(\hhat_{T-j+1}(X_{j-1})),
\)
for $j=t+1,\ldots,T$.
Since $p_j=\E{ \1{G_j}\midbar\cF_{j-1} }$, the tower property and linearity of conditional expectation give
\[
\sum_{j=t+1}^T \E{ p_j\midbar\cF_t }
=
\E{
\sum_{j=t+1}^T
\1{W_j\leq\hhat_{T-j+1}(X_{j-1})}
\midbar\cF_t
}
=
\vhat_{T-t}(X_t),
\]
proving the result for  $a=1$.
For $0<a<1$, H\"older's inequality with weights $a$ and $1-a$ applied respectively to $p_j$ and $1/(T-j+1)$ gives
\[
\sum_{j=t+1}^T
p_j^a \left(\frac{ 1 }{ T-j+1 }\right)^{1-a}
\leq
\left(\sum_{j=t+1}^T p_j\right)^a
\left(\sum_{j=1}^{T-t}\frac{ 1 }{ j }\right)^{1-a}.
\]
If we now take conditional expectations and apply conditional Jensen's inequality for the concave function $z \mapsto z^a$, we obtain
\begin{align*}
    \sum_{j=t+1}^T
    \E{
    \frac{ p_j^a }{ (T-j+1)^{1-a} }
    \midbar\cF_t
    }
    & \leq
     \E{ \left(\sum_{j=t+1}^T p_j\right)^a \midbar \cF_t }
      \left(\sum_{j=1}^{T-t}\frac{ 1 }{ j }\right)^{1-a}\\
    & \leq 
    \left( \sum_{j=t+1}^T\E{ p_j\midbar\cF_t }\right)^a \left(\sum_{j=1}^{T-t}\frac{ 1 }{ j }\right)^{1-a} \\
    & =
    \vhat_{T-t}(X_t)^a \left(\sum_{j=1}^{T-t}\frac{ 1 }{ j }\right)^{1-a}.
    \end{align*}
The claimed bound then follows from the harmonic-sum estimate $\sum_{j=1}^{T-t}j^{-1}\leq 1+\log(T-t)$.
\halmos
\endproof

% We now use Lemma~\ref{lem:series-upper-bound} with $a = 1/4$ and $t = 0$
% to sum the one-period approximation errors and show that their
% accumulated contribution is asymptotically negligible on the scale
% of $\vhat_T(c)$.

\proof{Proof of Proposition~\ref{prop:A_t-square-accumulated-approximation-error-negligible}.}
Fix $t \in [T]$ with $T-t+1 > K$, and suppose that $X_{t-1} > 0$.
Let 
\[
I_t =   \int_0^{\hhat_{T-t+1}(X_{t-1})}
            \left(
                1 - \frac{ w }{ \hhat_{T-t+1}(X_{t-1}) }
            \right)^2
            f(w) \, dw,
\]
be the integral of the squared linear approximation of the value-function difference on the selection interval against the weight distribution. 
Adding and subtracting $I_t$ and applying the triangle inequality, we obtain
\begin{equation*}
    \abs{
        \E{ A_t^2 \midbar \cF_{t-1} }
        - \sigma^2 F(\hhat_{T-t+1}(X_{t-1}))
    }
    \leq
        \abs{
        \E{ A_t^2 \midbar \cF_{t-1} }
        -
        I_t
    } 
    +
    \abs{
        I_t
        - \sigma^2 F(\hhat_{T-t+1}(X_{t-1}))
    }.
\end{equation*}

To bound the first term, the integral representation
\eqref{eq:At-conditional-second-moment} and
Corollary~\ref{cor:error-bound-of-linear-integrand-approximation-k},
applied at index $T-t \geq K$ and capacity $X_{t-1}$, give
\[
    \abs{
        \E{ A_t^2 \midbar \cF_{t-1} }
        -
        I_t
    } 
    \leq
    \Delta_2
    \left\{
        \frac{ F^{1/4}(\hhat_{T-t+1}(X_{t-1})) }{ (T-t+1)^{3/4} }
        + \frac{ 1 }{ T-t+1 }
    \right\}.
\]

To bound the second term, first note that
$\hhat_{T-t+1}(X_{t-1}) \leq \epsilon_{T-t+1}(X_{t-1})$.
Since $T-t+1 > K \geq K_1$, the monotonicity of $F$ and the lower
bound on the second ratio in
\eqref{eq:karamata-ratios-pointwise-bounds} at
$\epsilon = \epsilon_{T-t+1}(X_{t-1}) \leq \wbar$
give
\begin{equation}
\label{eq:ub-hhatF(hhat)}
        \hhat_{T-t+1}(X_{t-1}) F(\hhat_{T-t+1}(X_{t-1}))
        \leq
        \epsilon_{T-t+1}(X_{t-1}) F(\epsilon_{T-t+1}(X_{t-1})) 
        \leq
        C_F \frac{ X_{t-1} }{ T-t+1 }.
\end{equation}
Moreover, $0 < \hhat_{T-t+1}(X_{t-1}) \leq \wbar$, so
Lemma~\ref{lem:linear-approximation-second-moment-bounds},
applied with $h = \hhat_{T-t+1}(X_{t-1})$, gives
\begin{align*}
    \abs{
        I_t
        - \sigma^2 F(\hhat_{T-t+1}(X_{t-1}))
    }
    & \leq
    M_F \hhat_{T-t+1}(X_{t-1})^{3/4}
    F(\hhat_{T-t+1}(X_{t-1})) \\
    & \leq
    M_F (C_F c)^{3/4}
    \frac{ F^{1/4}(\hhat_{T-t+1}(X_{t-1})) }{ (T-t+1)^{3/4} },
\end{align*}
where the last inequality follows from
\eqref{eq:ub-hhatF(hhat)} and $X_{t-1} \leq c$.

Combining these two estimates gives a finite constant $C$,
independent of $T$, $t$, and the remaining capacity, such that,
for every $t \in [T]$,
\[
    \abs{
        \E{ A_t^2 \midbar \cF_{t-1} }
        - \sigma^2 F(\hhat_{T-t+1}(X_{t-1}))
    }
    \leq
    C \left\{
        \frac{ F^{1/4}(\hhat_{T-t+1}(X_{t-1})) }{ (T-t+1)^{3/4} }
        + \frac{ 1 }{ T-t+1 }
    \right\}.
\]
Taking expectations, summing over $t$, and applying
Lemma~\ref{lem:series-upper-bound} with $a = 1/4$ at time zero gives
\[
    \sum_{t=1}^T
    \E{
        \abs{
            \E{ A_t^2 \midbar \cF_{t-1} }
            - \sigma^2 F(\hhat_{T-t+1}(X_{t-1}))
        }
    }
    \leq
    C \left\{
        \vhat_T(c)^{1/4} (1 + \log T)^{3/4}
        + 1 + \log T
    \right\},
\]
where we also used $\sum_{t=1}^T \frac{1}{T-t+1} \leq 1 + \log T$.
For $T > K$, dividing both sides by $\vhat_T(c)$ and using the
little-$o$ estimate in
Proposition~\ref{prop:joint-properties-f-epsilon}--\ref{property:fluid-benchmark-growth-relative-to-logT},
together with the logarithmic regret bound in
Proposition~\ref{prop:log-regret-bound}, shows that the resulting ratio
converges to zero as $T \to \infty$, completing the proof.
\halmos
\endproof

\proof{Proof of Proposition~\ref{prop:B_t-square-total-contribution-negligible}.}
We begin by showing that the second limit follows from the first.
For $t = 1, \ldots, T$ and $x \in [0,c]$, the bounded-difference
estimate of Corollary~\ref{cor:uniform-bound-on-difference-in-value-function}
implies that
\[
    \abs{ 1 + \vhat_{T-t}(x-w) - \vhat_{T-t}(x) }
    \leq 1 + \Delta_v
    \qquad \text{for all } w \in [0, \hhat_{T-t+1}(x)].
\]
Taking absolute values in the integral representation
\eqref{eq:Bt-integral-representation}, applying the triangle inequality,
and squaring gives, for every $t \in [T]$, that
\[
    \abs{ B_t }^2
    \leq (1 + \Delta_v)^2 F(\hhat_{T-t+1}(X_{t-1}))^2.
\]
Taking expectations, summing over $t$, and dividing by $\vhat_T(c)$
gives, for all $T > K$,
\[
    0
    \leq
    \frac{ 1 }{ \vhat_T(c) } \sum_{t=1}^T \E{ B_t^2 }
    \leq
    (1 + \Delta_v)^2
    \frac{ 1 }{ \vhat_T(c) }
    \sum_{t=1}^T \E{ F(\hhat_{T-t+1}(X_{t-1}))^2 }.
\]
Hence, it remains to show that this last right-hand side converges to zero
as $T \to \infty$.

Let
\(
    m_T = \lfloor \sqrt{ \vhat_T(c) } \rfloor
\)
be the largest integer less than or equal to $\sqrt{ \vhat_T(c) }$.
By Proposition~\ref{prop:log-regret-bound} and
Proposition~\ref{prop:joint-properties-f-epsilon}--\ref{property:fluid-benchmark-growth-relative-to-logT},
we have $\vhat_T(c) \to \infty$, 
so we may take $T$ sufficiently large that $K < m_T \leq T$.
For $t \leq T-m_T$, at least $m_T + 1$ periods remain.
The definition of the threshold function \eqref{eq:threshold-function-definition}, the bound $X_{t-1} \leq c$, and the monotonicity
of the consumption function give
\[
    F(\hhat_{T-t+1}(X_{t-1}))
    \leq
    F(\epsilon_{T-t+1}(X_{t-1}))
    \leq
    F(\epsilon_{m_T}(c)).
\]
Consequently,
\[
    F(\hhat_{T-t+1}(X_{t-1}))^2
    \leq
    F(\epsilon_{m_T}(c))
    F(\hhat_{T-t+1}(X_{t-1}))
    \qquad
    \text{for } t \leq T-m_T.
\]
In the last $m_T $ periods, each squared probability is at most one.
Using the identity
\[
    \sum_{t=1}^T
        \E{ F(\hhat_{T-t+1}(X_{t-1})) }
    =
    \E{ \widehat N_T(c) }
    =
    \vhat_T(c),
\]
we obtain
\[
        0
        \leq
        \frac{ 1 }{ \vhat_T(c) }
        \sum_{t=1}^T
            \E{ F(\hhat_{T-t+1}(X_{t-1}))^2 }
        \leq
        F(\epsilon_{m_T}(c))
        + \frac{ m_T }{ \vhat_T(c) } 
        \leq
        F(\epsilon_{m_T}(c))
        + \frac{ 1 }{ \sqrt{ \vhat_T(c) } }.
\]
As $T \to \infty$, we have $m_T \to \infty$.
The integral representation
\eqref{eq:consumption-function-integral-representation}
and the positivity of $f$ on $(0,\wbar]$ then imply that
$\epsilon_{m_T}(c) \to 0$.
Thus, the first term on the last right-hand side tends to $F(0) = 0$
by continuity, while the second tends to zero because
$\vhat_T(c) \to \infty$.
This proves the first limit and, by the initial bound, the second.
\halmos
\endproof

% \subsection{Proof of the variance asymptotics}
% \label{subsec:variance-asymptotics}
% We now combine the estimates from
% Sections~\ref{sec:At-squared-contribution} and~\ref{sec:Bt-square-contribution}
% to complete the proof of Proposition~\ref{prop:variance-asymptotics}.

\proof{Proof of Proposition~\ref{prop:variance-asymptotics}.}
For $T > K$, using \eqref{eq:conditional-variance-decomposition}
at time zero, the tower property, and the identity
\[
    \sum_{t=1}^T
        \E{ F(\hhat_{T-t+1}(X_{t-1})) }
    =
    \E{ \widehat N_T(c) }
    =
    \vhat_T(c),
\]
we obtain
\begin{align*}
    \abs{ \frac{ \omega_T(c) }{ \vhat_T(c) } - \sigma^2 }
    & =
    \frac{ 1 }{ \vhat_T(c) }
    \abs{
        \sum_{t=1}^T \E{ A_t^2 }
        - \sum_{t=1}^T \E{ B_t^2 }
        - \sigma^2 \sum_{t=1}^T
            \E{ F(\hhat_{T-t+1}(X_{t-1})) }
    } \\
    & \leq
    \frac{ 1 }{ \vhat_T(c) }
    \sum_{t=1}^T
    \E{
        \abs{
            \E{ A_t^2 \midbar \cF_{t-1} }
            - \sigma^2 F(\hhat_{T-t+1}(X_{t-1}))
        }
    }
    +
    \frac{ 1 }{ \vhat_T(c) }
    \sum_{t=1}^T \E{ B_t^2 }.
\end{align*}
The two terms on the last right-hand side converge to zero as
$T \to \infty$ by
Propositions~\ref{prop:A_t-square-accumulated-approximation-error-negligible}
and~\ref{prop:B_t-square-total-contribution-negligible}, respectively.
Since $\omega_T(c) = \Var{ \widehat N_T(c) }$ and
$\vhat_T(c) = \E{ \widehat N_T(c) }$,
this proves the proposition.
\halmos
\endproof

\section{Proofs of the results in Section~\ref{sec:clt}}
\label{apx:clt-proofs}

\proof{Proof of Lemma~\ref{lemma:compensator-L2-convergence}.}
Since
\(
    \E{ \1{ G_t } \midbar \cF_{t-1} }
    =
    F(\hhat_{T-t+1}(X_{t-1})),
\)
the differences
\[
    D_t
    =
    F(\hhat_{T-t+1}(X_{t-1})) - \1{ G_t },
    \qquad t = 1, \ldots, T,
\]
are $\cF_t$-measurable and satisfy
$\E{ D_t \midbar \cF_{t-1} } = 0$.
Thus, $\{D_t: 1 \leq t \leq T\}$ is a martingale difference sequence
with respect to $\{\cF_t: 0 \leq t \leq T\}$, and
\[
    \sum_{t=1}^T D_t = \Lambda_T - \widehat N_T(c).
\]
Moreover, these differences are also orthogonal: for $s < t$, $D_s$ is
$\cF_{t-1}$-measurable, so
\(
    \E{ D_s D_t }
    =
    \E{ D_s \E{ D_t \midbar \cF_{t-1} } }
    = 0.
\)
Thus, 
\[
    \E{ (\Lambda_T - \widehat N_T(c))^2 }
    =
    \E{ \left( \sum_{t=1}^T D_t \right)^2 }
    =
    \sum_{t=1}^T \E{ D_t^2 }
    =
    \sum_{t=1}^T
    \E{ \Var{ \1{ G_t } \midbar \cF_{t-1} } }.
\]
Conditional on $\cF_{t-1}$, $\1{ G_t }$ is Bernoulli with success
probability $F(\hhat_{T-t+1}(X_{t-1}))$.
Therefore,
\begin{align*}
    \sum_{t=1}^T
    \E{ \Var{ \1{ G_t } \midbar \cF_{t-1} } }
    & =
    \sum_{t=1}^T
    \E{
        F(\hhat_{T-t+1}(X_{t-1}))
        \{ 1 - F(\hhat_{T-t+1}(X_{t-1})) \}
    } \\
    & \leq
    \sum_{t=1}^T
    \E{ F(\hhat_{T-t+1}(X_{t-1})) }
    = \vhat_T(c).
\end{align*}
If we now write $\Lambda_T  - \vhat_T(c) = ( \Lambda_T - \widehat N_T(c) ) + ( \widehat N_T(c) -  \vhat_T(c))$, square both sides, use the inequality 
$(a+b)^2 \leq 2a^2 + 2b^2$, and recall that 
$\E{ \widehat N_T(c) } = \vhat_T(c)$,
we obtain
\begin{equation*}
    \E{ (\Lambda_T - \vhat_T(c))^2 }
    \leq
    2 \E{ (\Lambda_T - \widehat N_T(c))^2 }
    + 2 \E{ (\widehat N_T(c) - \vhat_T(c))^2 } 
    \leq 2 \vhat_T(c) + 2 \omega_T(c).
\end{equation*}
Consequently,
\[
    \E{
        \left( \frac{ \Lambda_T }{ \vhat_T(c) } - 1 \right)^2
    }
    \leq
    \frac{ 2 }{ \vhat_T(c) }
    \left( 1 + \frac{ \omega_T(c) }{ \vhat_T(c) } \right)
    \longrightarrow 0
    \qquad \mbox{as } T \to \infty,
\]
because Proposition~\ref{prop:variance-asymptotics} gives
$\omega_T(c)/\vhat_T(c) \to \sigma^2$, while
$\vhat_T(c) \to \infty$ follows from
Proposition~\ref{prop:log-regret-bound} and
Proposition~\ref{prop:joint-properties-f-epsilon}--\ref{property:fluid-benchmark-growth-relative-to-logT}.
\halmos
\endproof

% The mean-square convergence in
% Lemma~\ref{lemma:compensator-L2-convergence} represents a crucial step towards proving Proposition~\ref{prop:quadratic-variation-convergence}. 
% The rest of the argument follows next, and it uses the bounds in
% Propositions~\ref{prop:A_t-square-accumulated-approximation-error-negligible}
% and~\ref{prop:B_t-square-total-contribution-negligible}
% one more time. 

\proof{Proof of Proposition~\ref{prop:quadratic-variation-convergence}.}
For all sufficiently large $T$, write
\[
    \frac{ Q_T }{ \omega_T(c) } - 1
    =
    \frac{ \vhat_T(c) }{ \omega_T(c) }
    \left( \frac{ Q_T }{ \vhat_T(c) } - \sigma^2 \right)
    + \sigma^2 \frac{ \vhat_T(c) }{ \omega_T(c) } - 1.
\]
The triangle inequality therefore gives
\[
    \E{ \abs{ \frac{ Q_T }{ \omega_T(c) } - 1 } }
    \leq
    \frac{ \vhat_T(c) }{ \omega_T(c) }
    \E{ \abs{ \frac{ Q_T }{ \vhat_T(c) } - \sigma^2 } }
    +
    \abs{ \sigma^2 \frac{ \vhat_T(c) }{ \omega_T(c) } - 1 }.
\]
By Proposition~\ref{prop:variance-asymptotics},
$\vhat_T(c)/\omega_T(c) \to 1/\sigma^2$ as $T \to \infty$.
Thus, the second term tends to zero, and it remains to show that
\begin{equation}\label{eq:QT/vT-convergence}
    \E{ \abs{ \frac{ Q_T }{ \vhat_T(c) } - \sigma^2 } }
    \longrightarrow 0
    \qquad \mbox{as } T \to \infty.
\end{equation}
Adding and subtracting $\sigma^2 \Lambda_T/\vhat_T(c)$ and applying
the triangle inequality gives
\[
    \E{ \abs{ \frac{ Q_T }{ \vhat_T(c) } - \sigma^2 } }
    \leq
    \E{
        \abs{
            \frac{ Q_T }{ \vhat_T(c) }
            - \sigma^2 \frac{ \Lambda_T }{ \vhat_T(c) }
        }
    }
    + \sigma^2
    \E{ \abs{ \frac{ \Lambda_T }{ \vhat_T(c) } - 1 } }.
\]
Here, the second term tends to zero by
Lemma~\ref{lemma:compensator-L2-convergence} and the
Cauchy--Schwarz inequality.
To bound the first term, the definition of the predictable quadratic
variation \eqref{eq:QT} and the decomposition
\eqref{eq:dt-squared-decomposition} give
\[
    Q_T
    =
    \sum_{t=1}^T \E{ A_t^2 \midbar \cF_{t-1} }
    - \sum_{t=1}^T B_t^2.
\]
Hence, the definition of $\Lambda_T$ in \eqref{eq:LambdaT} and
another application of the triangle inequality give
\begin{align*}
    \E{
        \abs{
            \frac{ Q_T }{ \vhat_T(c) }
            - \sigma^2 \frac{ \Lambda_T }{ \vhat_T(c) }
        }
    }
    & =
    \frac{ 1 }{ \vhat_T(c) }
    \E{
        \abs{
            \sum_{t=1}^T
            \left\{
                \E{ A_t^2 \midbar \cF_{t-1} }
                - \sigma^2 F(\hhat_{T-t+1}(X_{t-1}))
            \right\}
            - \sum_{t=1}^T B_t^2
        }
    } \\
    & \leq
    \frac{ 1 }{ \vhat_T(c) }
    \sum_{t=1}^T
    \E{
        \abs{
            \E{ A_t^2 \midbar \cF_{t-1} }
            - \sigma^2 F(\hhat_{T-t+1}(X_{t-1}))
        }
    }
    + \frac{ \sum_{t=1}^T \E{ B_t^2 } }{ \vhat_T(c) }.
\end{align*}
As $T \to \infty$, the last right-hand side converges to zero by
Propositions~\ref{prop:A_t-square-accumulated-approximation-error-negligible}
and~\ref{prop:B_t-square-total-contribution-negligible},
concluding the proof of \eqref{eq:QT/vT-convergence} and of the $L^1$ convergence in the proposition. 
The asserted convergence in probability then follows by Markov's inequality.
\halmos
\endproof

% \section{Proof of Theorem~\ref{thm:main-theorem}}
% \label{subsec:proof-of-clt}

% We now apply Proposition~\ref{prop:martingale-clt} to the normalized
% martingale differences $Z_t$, recenter and rescale using the fluid benchmark, 
% and transfer the resulting limit to the optimal policy.

\proof{Proof of Theorem~\ref{thm:main-theorem}.}
We verify that the two conditions \eqref{eq:negligibility-condition-of-martingale-clt} and \eqref{eq:convergence-in-probability-condition-of-martingale-clt} in Proposition~\ref{prop:martingale-clt}
hold for the normalized martingale differences 
\(
Z_t = \frac{d_t}{\sqrt{\omega_T(c)}},
\)
$1 \leq t \leq T$, in \eqref{eq:Zt}.
We recall that $d_t = A_t + B_t$.
For $t > T-K$, the threshold $\hhat_{T-t+1}$ is zero, and
$A_t = B_t = 0$ almost surely.
For $t \leq T-K$, the index $T-t$ is at least $K$, so \eqref{eq:At} and
Corollary~\ref{cor:uniform-bound-on-difference-in-value-function},
applied at index $T-t$, give
$\abs{ A_t } \leq 1 + \Delta_v$.
Furthermore, we also have $B_t = -\E{ A_t \midbar \cF_{t-1} }$, so
$\abs{ B_t } \leq 1 + \Delta_v$.
It then follows that 
\[
    \max_{1 \leq t \leq T} \norm{ Z_t }_\infty
    \leq
    \frac{ 2(1 + \Delta_v) }{ \sqrt{ \omega_T(c) } }
    \longrightarrow 0,
    \quad \text{as } T \to \infty,
\]
by the asymptotic equivalence $\omega_T(c) \sim \sigma^2 \vhat_T(c)$
and the divergence of $\vhat_T(c)$ from 
the logarithmic regret bound in
Proposition~\ref{prop:log-regret-bound} and 
Proposition~\ref{prop:joint-properties-f-epsilon}--\ref{property:fluid-benchmark-growth-relative-to-logT}.

For the normalized predictable quadratic variation condition \eqref{eq:convergence-in-probability-condition-of-martingale-clt},
Proposition~\ref{prop:quadratic-variation-convergence} immediately gives
\[
    \sum_{t=1}^T \E{ Z_t^2 \midbar \cF_{t-1} }
    =
    \frac{ Q_T }{ \omega_T(c) }
    \xrightarrow{p} 1,
    \qquad \mbox{as } T \to \infty.
\]
Both conditions of Proposition~\ref{prop:martingale-clt} are
therefore satisfied, so
\[
    \frac{ \widehat N_T(c) - \vhat_T(c) }{ \sqrt{ \omega_T(c) } }
    =
    \sum_{t=1}^T Z_t
    \xrightarrow{d} N(0,1),
    \qquad \mbox{as } T \to \infty.
\]

To center and scale by the fluid benchmark, write
\[
    \frac{ \widehat N_T(c) - T F(\epsilon_T(c)) }
         { \sqrt{ T F(\epsilon_T(c)) } }
    =
    \frac{ \widehat N_T(c) - \vhat_T(c) }{ \sqrt{ \omega_T(c) } }
    \sqrt{ \frac{ \omega_T(c) }{ T F(\epsilon_T(c)) } }
    +
    \frac{ \vhat_T(c) - T F(\epsilon_T(c)) }
         { \sqrt{ T F(\epsilon_T(c)) } }.
\]
By Proposition~\ref{prop:log-regret-bound} and
Proposition~\ref{prop:joint-properties-f-epsilon}--\ref{property:fluid-benchmark-growth-relative-to-logT},
the second term tends to zero and
$\vhat_T(c)/\{ T F(\epsilon_T(c)) \} \to 1$.
Together with
$\omega_T(c)/\vhat_T(c) \to \sigma^2$, this implies that the
square-root factor in the first term tends to $\sigma$.
Slutsky's theorem therefore gives
\[
    \frac{ \widehat N_T(c) - T F(\epsilon_T(c)) }
         { \sqrt{ T F(\epsilon_T(c)) } }
    \xrightarrow{d} N(0,\sigma^2),
    \qquad \mbox{as } T \to \infty,
\]
which proves the asserted central limit theorem for the reoptimized
heuristic.

To prove the central limit theorem for the number of item selections made by the optimal policy $\pi^*$, run both
policies on the same sequence of item weights.
Let $N_T^{\mathrm{off}}(c)$ denote the maximum number of items that
can be selected from this sequence with capacity $c$ in hindsight.
On each realization, both online selection counts are bounded above
by $N_T^{\mathrm{off}}(c)$, so Proposition~\ref{prop:log-regret-bound} and the optimality of $\pi^*$
imply
\[
    \E{ \abs{ N_T^*(c) - \widehat N_T(c) } }
    \leq
    \E{ N_T^{\mathrm{off}}(c) - N_T^*(c) }
    + \E{ N_T^{\mathrm{off}}(c) - \widehat N_T(c) } 
    \leq
    2 C (1 + \log T).
\]
By Proposition~\ref{prop:joint-properties-f-epsilon}--\ref{property:fluid-benchmark-growth-relative-to-logT}, dividing the preceding bound by $\sqrt{ T F(\epsilon_T(c)) }$ gives
\[
    \frac{ N_T^*(c) - \widehat N_T(c) }
         { \sqrt{ T F(\epsilon_T(c)) } }
    \longrightarrow 0,
    \qquad \mbox{in } L^1 \mbox{ as } T \to \infty.
\]
Since convergence in $L^1$ implies convergence in probability,
Slutsky's theorem and the central limit theorem for the reoptimized
heuristic give
\[
    \frac{ N_T^*(c) - T F(\epsilon_T(c)) }
         { \sqrt{ T F(\epsilon_T(c)) } }
    =
    \frac{ N_T^*(c) - \widehat N_T(c) }
         { \sqrt{ T F(\epsilon_T(c)) } }
    +
    \frac{ \widehat N_T(c) - T F(\epsilon_T(c)) }
         { \sqrt{ T F(\epsilon_T(c)) } }
    \xrightarrow{d} N(0,\sigma^2)
    \qquad \mbox{as } T \to \infty,
\]
just as needed.
\halmos
\endproof

\end{APPENDICES}

%%%%%%%%%%%%%%%%%
\end{document}